\documentclass[12pt]{amsart}
\usepackage[]{amscd,amsfonts,amsmath,amsxtra,amssymb}
\usepackage{graphics}
\usepackage[all]{xy}
\usepackage{hyperref}
\newtheorem{sub}{}[section]
\newtheorem{subsub}{}[sub]
\def\ov#1{\overline{#1}}

\def\coker{\mathop{\rm coker}\nolimits}
\def\Hom{\mathop{\rm Hom}\nolimits}
\def\HHom{\mathop{\mathcal Hom}\nolimits}
\def\Ext{\mathop{\rm Ext}\nolimits}
\def\EExt{\mathop{\mathcal Ext}\nolimits}
\def\Tor{\mathop{\rm Tor}\nolimits}
\def\Hilb{\mathop{\rm Hilb}\nolimits}
\def\Pic{\mathop{\rm Pic}\nolimits}
\def\Aut{\mathop{\rm Aut}\nolimits}

\def\End{\mathop{\rm End}\nolimits}
\def\EEnd{\mathop{\mathcal End}\nolimits}
\def\GL{\mathop{\rm GL}\nolimits}

\def\imm{\mathop{\rm im}\nolimits}
\def\deg{\mathop{\rm deg}\nolimits}
\def\Deg{\mathop{\rm Deg}\nolimits}
\def\spec{\mathop{\rm spec}\nolimits}
\def\lra{\longrightarrow}

\def\sigg{\mathop{\hbox{$\displaystyle\sum$}}\limits}
\def\psigg{\mathop{\hbox{$\displaystyle\bigotimes$}}\limits}

\def\hfl#1#2{\smash{\mathop{\ \hbox to 12mm{\rightarrowfill}}
\limits^{\scriptstyle#1}_{\scriptstyle#2} \ }}
\def\hflb#1#2{\smash{\mathop{\hbox to 12mm{\leftarrowfill}}
\limits^{\scriptstyle#1}_{\scriptstyle#2}}}

\def\m#1{{\hbox{$#1$}}}

\def\ot{\otimes}
\def\og{\leavevmode\raise.3ex\hbox{$\scriptscriptstyle\langle\!\langle$}}
\def\fg{\leavevmode\raise.3ex\hbox{$\scriptscriptstyle\,\rangle\!\rangle$}}

\def\nsp{\lbrace 0\rbrace}

\def\Ssect#1#2{\pagebreak[3]\begin{sub}\label{#2}{\sc\small\small  #1}\rm\medskip}

\def\sepsec{\vskip 1.5cm}
\def\sepsub{\vskip 0.6cm}
\def\sepsubsub{\vskip 0.4cm}
\def\sepprop{\vskip 0.4cm}

\def\xmat#1{\[\xymatrix{#1}\]}

\def\flinc{\ar@{^{(}->}}
\def\fleq{\ar@{=}}
\def\flon{\ar@{->>}}
\def\fmaps{\ar@{|-{>}}}
\def\wT{{\widetilde T}}
\def\Nligne{\hfil\break}
\newcommand{\N}{{\mathbb N}}

\newcommand{\Z}{{\mathbb Z}}

\newcommand{\C}{{\mathbb C}}

\renewcommand{\P}{{\mathbb P}}

\newcommand{\F}{{\mathbb F}}
\newcommand{\E}{{\mathbb E}}
\newcommand{\G}{{\mathbb G}}

\renewcommand{\L}{{\mathbb L}}

\newcommand{\ka}{{\mathcal A}}

\newcommand{\kc}{{\mathcal C}}
\newcommand{\kd}{{\mathcal D}}
\newcommand{\ke}{{\mathcal E}}
\newcommand{\kf}{{\mathcal F}}
\newcommand{\kg}{{\mathcal G}}
\newcommand{\kh}{{\mathcal H}}
\newcommand{\ki}{{\mathcal I}}
\newcommand{\kj}{{\mathcal J}}
\newcommand{\kk}{{\mathcal K}}
\newcommand{\kl}{{\mathcal L}}
\newcommand{\km}{{\mathcal M}}
\newcommand{\kn}{{\mathcal N}}
\newcommand{\ko}{{\mathcal O}}

\newcommand{\kt}{{\mathcal T}}
\newcommand{\ku}{{\mathcal U}}
\newcommand{\kv}{{\mathcal V}}

\newcommand{\kx}{{\mathcal X}}
\newcommand{\ky}{{\mathcal Y}}

\begin{document}

\def\refname{R\'ef\'erences}
\def\contentsname{Sommaire}
\def\proofname{D\'emonstration}

\author{Jean--Marc Dr\'{e}zet}
\address{
Institut de Math\'ematiques de Jussieu,
Case 247,
4 place Jussieu,
F-75252 Paris, France}
\email{drezet@math.jussieu.fr}
\urladdr{http://www.math.jussieu.fr/\~{}drezet}

\begin{abstract}
This paper is devoted to the study of coherent sheaves on non reduced curves
that can be locally embedded in smooth surfaces. If $Y$ is such a curve then
there is a filtration \m{C\subset C_2\subset\cdots\subset C_n=Y} such that $C$
is the reduced curve associated to $Y$, and for every \m{P\in C} there exists
\m{z\in\ko_{Y,P}} such that \m{(z^i)} is the ideal of \m{C_i} in \m{\ko_{Y,P}}.
We define, using canonical filtrations, new invariants of coherent sheaves on
$Y$ : the {\em generalized rank} and {\em degree}, and use them to state a {\em
Riemann-Roch theorem} for sheaves on $Y$. We define {\em quasi locally free
sheaves}, which are locally isomorphic to direct sums of \m{\ko_{C_i}}, and
prove that every coherent sheaf on $Y$ is quasi locally free on some nonempty
open subset of $Y$. We give also a simple criterion of quasi locally freeness.
We study the ideal sheaves \m{\ki_{n,Z}} in $Y$ of finite subschemes $Z$ of $C$.
When $Y$ is embedded in a smooth surface we deduce some results on deformations
of \m{\ki_{n,Z}} (as sheaves on $S$). When \m{n=2}, i.e. when $Y$ is a {\em
double curve}, we can completely describe the torsion free sheaves on $Y$. In
particular we show that these sheaves are reflexive. The torsion free sheaves of
generalized rank 2 on \m{C_2} are of the form \m{\ki_{2,Z}\ot\kl}, where $Z$ is
a finite subscheme of $C$ and $\kl$ is a line bundle on $Y$. We begin the study
of moduli spaces of stable sheaves on a double curve, of generalized rank 3 and
generalized degree $d$. These moduli spaces have many components. Sometimes
one of them is a multiple structure on the moduli space of stable vector bundles
on $C$ of rank 3 and degree $d$.
\end{abstract}

\title[{\tiny Faisceaux coh\'erents sur les courbes multiples}]
{Faisceaux coh\'erents sur les courbes multiples}
\maketitle
\tableofcontents

\section{Introduction}
\label{intro}

Les courbes projectives multiples {\em primitives} ont \'et\'e d\'efinies et
\'etudi\'ees par C.~B\u anic\u a et O.~Forster dans \cite{ba_fo}. Les cas les
plus simples et qu'on \'etudiera en d\'etail ici sont les courbes non r\'eduites
plong\'ees dans une surface lisse, et dont la courbe r\'eduite associ\'ee est
projective lisse.

Les faisceaux semi-stables sur des courbes projectives r\'eduites non lisses
on \'et\'e \'etudi\'es par de nombreux auteurs, notamment par C.S.~Seshadri
dans \cite{ses}, et U.N. Bhosle dans \cite{bho}, \cite{bho2} et d'autres
articles. Les faisceaux semi-stables sur des vari\'et\'es non r\'eduites 
semblables \`a celles qui sont consid\'er\'ees ici sont le sujet de l'article
\cite{in} de M.-A. Inaba. L'article \cite{in2} du m\^eme auteur traite des
faisceaux stables sur une vari\'et\'e non irr\'eductible ayant deux composantes
qui se coupent. Il est possible qu'on puisse obtenir par l'\'etude
des faisceaux coh\'erents sur les courbes non r\'eduites des r\'esultats sur la
d\'eg\'en\'eration des fibr\'es vectoriels ou des vari\'et\'es des modules de
fibr\'es semi-stables sur les courbes lisses. Certains r\'esultats ont d\'ej\`a
\'et\'e obtenus en utilisant des courbes r\'eduites mais singuli\`eres (cf.
\cite{sun0}, \cite{sun1}).

Les faisceaux coh\'erents sur les courbes non r\'eduites interviennent aussi
lorsqu'on veut \'etudier les faisceaux de dimension 1 sur les surfaces.
Les faisceaux sur des courbes non r\'eduites apparaissent
(cf. \cite{lp}, \cite{lp4}) comme limites de fibr\'es vectoriels sur des
courbes lisses. Leur r\^ole est sans doute plus important si on cherche \`a
obtenir d'autres vari\'et\'es de modules fins de faisceaux de dimension 1 que
les classiques vari\'et\'es de modules de faisceaux semi-stables (cf.
\cite{dr}).

Le but du pr\'esent article est de donner les bases de l'\'etude des faisceaux
coh\'erents sur une courbe multiple primitive $Y$ et de leurs vari\'et\'es de
modules. On introduit deux nouveaux invariants des faisceaux coh\'erents : le
{\em rang} et le {\em degr\'e g\'en\'eralis\'es}, avec lesquels on peut
\'enoncer un {\em th\'eor\`eme de Riemann-Roch} sur les courbes primitives.
On s'int\'eressera aux faisceaux g\'en\'eriques qui sont ici les faisceaux {\em
quasi localement libres} jouant le m\^eme r\^ole que les faisceaux localement
libres sur les vari\'et\'es lisses. On \'etudiera ensuite les faisceaux
d'id\'eaux de sous-sch\'emas finis de la courbe r\'eduite associ\'ee $C$, qui
sont les premiers exemples non triviaux de faisceaux sur $Y$. On s'int\'eressera
enfin aux {\em courbes doubles}. Dans ce cas on peut d\'ecrire pr\'ecis\'ement
les faisceaux sans torsion sur $Y$ et prouver en particulier qu'ils sont
r\'eflexifs. Pour finir on s'int\'eressera aux vari\'et\'es de modules de
faisceaux stables de rang g\'en\'eralis\'e 3 et de degr\'e g\'en\'eralis\'e $d$
sur une courbe double et on mettra en \'evidence de multiples composantes. Une
d'elles est une structure multiple sur la vari\'et\'e de modules des fibr\'es
vectoriels de rang 3 et de degr\'e $d$ sur $C$.

\sepsub

\Ssect{Faisceaux coh\'erents sur les courbes multiples primitives}{intro_1}

{\bf Courbes multiples primitives - }
Soit $C$ une courbe alg\'ebrique projective irr\'eductible lisse, de genre
\m{g_C}, plong\'ee dans une vari\'et\'e projective lisse $X$ de dimension 3 sur
$\C$. Soit \m{Y\subset X} une sous-vari\'et\'e ferm\'ee
de Cohen-Macaulay dont la
sous-vari\'et\'e r\'eduite associ\'ee est $C$. On dit que $Y$ est {\em
primitive} si elle peut localement \^etre plong\'ee dans une surface~: pour
tout point $P$ de $C$ il existe une surface \m{S\subset X} et un ouvert $U$ de
$X$ contenant $P$ tels que \m{U\cap S} soit lisse et \ \m{Y\cap U\subset S} .
Si tel est le cas il existe un entier \m{n>0} tel que pour tout point
\m{P\in C} on peut choisir l'ouvert $U$ pr\'ec\'edent et des coordonn\'ees
locales en $P$, $x$, $z$, $t$ de telle sorte que l'id\'eal de $S$ dans $U$ soit
\m{(x)} et celui de $Y$ \m{(x,z^n)}. L'entier $n$ s'appelle la {\em
multiplicit\'e} de $Y$ et $C$ son {\em support}.

Il existe une filtration canonique
\[C=C_1\subset\cdots\subset C_n=Y\ ,\]
o\`u \m{C_i} est l'intersection de $Y$ et du \m{i^{\rm eme}} voisinage
infinit\'esimal de $C$ dans $X$. Au voisinage de $P$ l'id\'eal de \m{C_i} dans
$U$ est donc \m{(x,z^i)}. On note, pour \m{1\leq i\leq n}, \m{\ko_i} le faisceau
structural de \m{C_i} (en particulier, \m{\ko_1=\ko_C} et \m{\ko_n=\ko_Y}), et
\m{\ki_C} le faisceau d'id\'eaux de $C$ dans $Y$.

\sepsubsub

{\bf Rang et degr\'e g\'en\'eralis\'es - Th\'eor\`eme de Riemann-Roch - }
Soit $\ke$ un faisceau coh\'erent sur $Y$. On d\'efinit en \ref{QLL-def} la
{\em premi\`ere filtration canonique} de $\ke$ :
\[\ke_{n+1}=0\subset\ke_n\subset\cdots\subset\ke_1=\ke \ .\]
Pour \m{1\leq i\leq n}, \m{\ke_{i+1}=\ki_C^i\ke} est le noyau de la restriction \
\m{\ke_{i}\to\ke_{i\mid C}} .
On a donc \ \m{\ke_i/\ke_{i+1}=\ke_{i\mid C}}, \m{\ke/\ke_{i+1}=\ke_{\mid
C_i}} . Le gradu\'e \
\m{{\rm Gr}(\ke)=\bigoplus_{i=1}^n\ke_i/\ke_{i+1}} \
est un faisceau de \m{\ko_C}-modules.

Les entiers \
\m{R(\ke)=rg({\rm Gr}(\ke))} , \m{\Deg(\ke)=\deg({\rm Gr}(\ke))} \
s'appellent respectivement le {\em rang g\'en\'eralis\'e} et le {\em degr\'e
g\'en\'eralis\'e} de $\ke$. On montre que le rang et le degr\'e
g\'en\'eralis\'es sont {\em additifs} : si \
\m{0\to\ke'\to\ke\to\ke''\to 0} \
est une suite exacte de faisceaux coh\'erents sur \m{C_n} alors on a
\[R(\ke) \ = \ R(\ke')+R(\ke") , \ \ \ \
\Deg(\ke) \ = \ \Deg(\ke')+\Deg(\ke") .\]
On d\'eduit des d\'efinitions le th\'eor\`eme de Riemann-Roch pour les courbes
primitives (cf. \ref{RR}) :

\sepprop

{\bf Th\'eor\`eme : } {\em Si $\ke$ est un faisceau coh\'erent sur \m{C_n},
alors on a
\[\chi(\ke) \ = \ \Deg(\ke)+R(\ke)(1-g_C) .\]}

\sepprop

On montre aussi que le rang et le degr\'e g\'en\'eralis\'es sont {\em
invariants par d\'eformation}, c'est-\`a-dire que dans une famille plate de
faisceaux coh\'erents sur \m{C_n} param\'etr\'ee par une vari\'et\'e
irr\'eductible, tous les faisceaux ont les m\^emes rang et degr\'e
g\'en\'eralis\'es.

On utilisera aussi la {\em seconde filtration canonique} de $\ke$
:  \Nligne \m{\ke^{(n+1)}=\nsp\subset \ke^{(n)}\subset\cdots\subset
\ke^{(2)}\subset \ke^{(1)}=\ke} , o\`u $\ke^{(i)}$ est le sous-faisceau de $\ke$
annulateur de \m{\ki_C^{n+1-i}}. Cette filtration a \'et\'e d\'efinie et
utilis\'ee par M.A. Inaba dans \cite{in}.

\sepsubsub

{\bf Faisceaux semi-stables - }\label{i_h}
On d\'eduit du th\'eor\`eme de Riemann-Roch
pr\'ec\'edent le calcul des polyn\^omes de Hilbert des faisceaux coh\'erents
sur $Y$(cf. \ref{hilb_pol}). Il en d\'ecoule que la d\'efinition des faisceaux
de dimension 1 sur $Y$
(semi-)stables au sens de C.~Simpson (cf. \cite{si}) est analogue \`a celle des
fibr\'es (semi-)stables sur les courbes lisses : un faisceau $\ke$ de dimension
1 sur $Y$ est {\em semi-stable} (resp. {\em stable}) si et seulement si il est
pur et si pour tout sous-faisceau propre \m{\kf\subset\ke} on a
\[\frac{\Deg(\kf)}{R(\kf)} \ \leq \ \frac{\Deg(\ke)}{R(\ke)}
\quad\quad\quad {\rm (resp. }\quad < \quad{\rm )} .\]

\sepsubsub

{\bf Fibr\'es vectoriels et groupe de Picard - }
On d\'emontre (th\'eor\`eme \ref{pr5}) que
si \m{1\leq i\leq n}, tout fibr\'e vectoriel sur
\m{C_i} peut se prolonger en fibr\'e vectoriel sur \m{C_n}.

On en d\'eduit bri\`evement en \ref{piccn} et \ref{piccn_desc} la description du
groupe de Picard de \m{C_n}. Le morphisme de restriction
\m{\Pic(C_n)\to\Pic(C)}
est surjectif et son noyau, le groupe des fibr\'es en droites sur \m{C_n} dont
la restriction \`a $C$ est triviale est une somme directe de groupes \m{\G_a}.

\sepsubsub

{\bf Faisceaux quasi localement libres - }
Soient \m{P\in C} et $M$ un \m{\ko_{n,P}}-module de type fini. On dit que $M$
est {\em quasi libre} s'il existe des entiers \m{m_1,\ldots,m_n} non n\'egatifs
tels que
\[M \ \simeq \  \ \bigoplus_{i=1}^nm_i\ko_{i,P} .\]
Si tel est le cas la suite \m{(m_1,\ldots,m_n)} est unique, on l'appelle le
{\em type} de $M$.

Soient $\ke$ un faisceau coh\'erent sur \m{C_n} et \m{P\in C}. On dit que $\ke$
est {\em quasi localement libre en} $P$ s'il existe des entiers
\m{m_1,\ldots,m_n} non n\'egatifs et un ouvert \m{U\subset C} contenant $P$
tels que
\[\ke_{\mid U} \ \simeq \ \bigoplus_{i=1}^nm_i\ko_{i\mid U} .\]
Dans ce cas le \m{\ko_{n,P}}-module \m{\ke_P} est quasi libre.
On dit que $\ke$ est {\em quasi localement libre} s'il l'est en tout point de
$C$. On d\'emontre en \ref{str_gen} le

\sepprop

{\bf Th\'eor\`eme : }{\em
Soit $\ke$ un faisceau coh\'erent sur \m{C_n}. Alors il existe un ouvert non
vide $U$ tel que $\ke$ soit quasi localement libre en tout point de $U$.}

\sepprop

Rappelons que tout faisceau coh\'erent sur une vari\'et\'e r\'eduite est
localement libre sur un ouvert dense. On peut donc dire que les les faisceaux
quasi localement libres sur les courbes primitives jouent le m\^eme r\^ole que
les faisceaux localement libres sur les vari\'et\'es r\'eduites.

On donne en \ref{str_gen0b} une caract\'erisation des faisceaux quasi
localement libres :

\sepprop

{\bf Th\'eor\`eme : }{\em
Soit $\ke$ un faisceau coh\'erent sur \m{C_n}. Alors $\ke$ est quasi localement
libre si et seulement si tous les termes de \m{{\rm Gr}(\ke)} sont localement
libres sur $C$.}

\sepprop

Les faisceaux quasi localement libres ont des propri\'et\'es semblables \`a
celles des faisceaux localement libres. En particulier on montre (th\'eor\`eme
\ref{str_mor}) que le noyau d'un morphisme surjectif de faisceaux coh\'erents
quasi localement libres l'est.
\end{sub}

\sepsub

\Ssect{Faisceaux sans torsion sur les courbes primitives doubles}{intro_5}

On suppose que \m{n=2}. On dit alors que \m{C_2} est une courbe primitive
{\em double}. Soit \m{L=\ki_C} le faisceau d'id\'eaux de $C$ dans
\m{C_2}. On peut dans ce cas d\'ecrire compl\`etement les faisceaux quasi
localement libres et les faisceaux sans torsion sur \m{C_2} (par d\'efinition
un faisceau est dit {\em sans torsion} s'il ne poss\`ede pas de sous-faisceau
dont le support est de dimension 0).

\sepsubsub

{\bf Faisceaux quasi localement libres - } 
Soient $\ke$ un faisceau quasi localement libre sur \m{C_2} et \m{E\subset\ke}
sa premi\`ere filtration canonique. Si \m{F=\ke/E} on a donc une suite exacte
\[0\lra E\lra\ke\lra F\lra 0 ,\]
et $E$, $F$ sont des fibr\'es vectoriels sur $C$. Le morphisme canonique \
\m{\ke\ot\ki_C\to\ke} \ induit un morphisme surjectif
\ \m{\Phi_\ke:F\ot L\to E .}
On a un isomorphisme canonique
\[\EExt^1_{\ko_2}(F,E) \ \simeq \ \HHom(F\ot L,E)\]
et \m{\Phi_\ke} n'est autre que l'image dans \m{H^0(\EExt^1_{\ko_2}(F,E))} de
l'\'el\'ement de \m{\Ext^1_{\ko_2}(F,E)} associ\'e \`a la suite exacte
pr\'ec\'edente. Le faisceau $\ke$ est localement libre si et seulement si
\m{\Phi_\ke} est un isomorphisme. En g\'en\'eral le noyau du morphisme
compos\'e \
\m{\ke\to F\to E\ot L^*} \
(o\`u \m{L^*} est le fibr\'e dual de $L$ sur $C$ et le second morphisme
provient de \m{\Phi_\ke}) est un fibr\'e vectoriel sur $C$ et c'est le plus
grand sous-faisceau de $\ke$ de support $C$.

R\'eciproquement, si \ \m{\sigma\in\Ext^1_{\ko_2}(F,E)} \ est tel le morphisme
associ\'e \ \m{\Phi:E\ot L\to F} \ soit surjectif, le faisceau $\ke$ extension
de $F$ par $E$ d\'efini par $\sigma$ est quasi localement libre,
\m{E\subset\ke} est sa premi\`ere filtration canonique et \ \m{\Phi_\ke=\Phi}.

Soit \m{P\in C}. On donne en \ref{Def_1} une condition n\'ecessaire et
suffisante pour qu'un \m{\ko_{2,P}}-module quasi libre de type \m{(m_1,m_2)} se
d\'eforme en modules quasi libres de type \m{(n_1,n_2)} : on doit avoir \
\m{m_1+2m_2=n_1+2n_2} \ et \ \m{n_1\geq m_1}.

\sepsubsub

{\bf Faisceaux sans torsion - }
Soient $\ke$ un faisceau coh\'erent sans torsion sur \m{C_2} et \m{E\subset\ke}
sa premi\`ere filtration canonique. Alors $E$ est un fibr\'e vectoriel mais
\m{\ke/E} a de la torsion si $\ke$ n'est pas quasi localement libre. Posons
\m{\ke/E=F\oplus T}  ,
o\`u $F$ est un fibr\'e vectoriel sur $C$ et $T$ un faisceau de torsion sur
$C$. On appelle {\em index} de $\ke$ l'entier
\ \m{i(\ke)=h^0(T)} .
Soient \m{P\in C}, \m{z\in \ko_{2,P}} une \'equation locale de $C$ et \m{x\in
\ko_{2,P}} un \'el\'ement au dessus d'un g\'en\'erateur de l'id\'eal maximal de
\m{\ko_{C,P}}. L'id\'eal maximal de \m{\ko_{2,P}} est donc \m{(x,z)}. On note
\m{I_{k,P}=(x^k,z)} , pour tout entier positif $k$.

On donne dans \ref{DUAL_ST} la structure des
\m{\ko_{2,P}}-modules de type fini sans torsion. Si $M$ en est un, alors il
existe des entiers $m$, $p$ et une suite d'entiers \m{n_1,\ldots,n_p} tels que
\[M \ \simeq \ \biggl(\bigoplus_{i=1}^pI_{n_i,P}\biggr)\oplus m\ko_{2,P} .\]
Les deux r\'esultats suivants sont d\'emontr\'es dans \ref{DUAL_ST} :

\sepprop

{\bf Th\'eor\`eme : }{\em Tout faisceau coh\'erent sans torsion sur \m{C_2} est
r\'eflexif.}

\sepprop

On a donc le m\^eme r\'esultat que sur les courbes lisses.

\sepprop

{\bf Th\'eor\`eme : }{\em Tout faisceau coh\'erent sans torsion sur \m{C_2} est
isomorphe au noyau d'un morphisme surjectif \ \m{\kf\to T}, o\`u $\kf$ est un
faisceau quasi localement libre sur \m{C_2} et $T$ un faisceau de torsion sur
$C$.}

\sepprop

On a donc un r\'esultat analogue \`a ce qu'on a sur les surfaces lisses. Le
fibr\'e $\kf$ n'est pas en g\'en\'eral unique, mais on peut les d\'ecrire tous.
\end{sub}

\sepsub

\Ssect{Faisceaux d'id\'eaux de points sur les courbes multiples
primitives}{XXX20}

Soit \m{Z\subset C} un sous-ensemble fini. On note \m{\ki_{n,Z}} le faisceau
d'id\'eaux de $Z$ sur \m{C_n}.

\sepsubsub

\begin{subsub}{\bf Limites de fibr\'es en droites sur les courbes lisses - }\rm 
Pour montrer l'int\'er\^et de l'\'etude des faisceaux d'id\'eaux de points sur
les courbes multiples, examinons les cas du plan projectif \m{\P_2}. Supposons
que \m{C_n\subset\P_2}, $C$ \'etant une courbe de degr\'e $d$. Soit $\km$ la
vari\'et\'e de modules de faisceaux de dimension 1 contenant tous les fibr\'es
en droites de degr\'e 0 sur les courbes lisses de degr\'e $nd$ de \m{\P_2}. On a
\m{\dim(\km)=n^2d^2+1} (cf. \cite{lp}, \cite{lp4}). Les faisceaux du type
\m{\ki_{n,Z}(k)}, o\`u $Z$ a \m{knd} points sont limites de faisceaux de $\km$.
\'Etant donn\'e que si \m{\ki_{n,Z}\simeq\ki_{n,Z'}} alors \m{Z=Z'} (car $Z$ est
pr\'ecis\'ement le lieu des points de $C$ o\`u \m{\ki_{n,Z}} n'est pas libre),
les faisceaux du type \m{\ki_{n,Z}(k)} (avec $k$ fix\'e et \m{\#Z=knd})
constituent une famille de dimension \m{knd} de faisceaux limites de $\km$. On
obtient donc des familles de dimension arbitrairement grande de faisceaux
limites de $\km$. 
\end{subsub}

\sepsubsub

\begin{subsub}{\bf R\'esultats g\'en\'eraux - }\rm 
On calcule en \ref{def_id} les dimensions de \m{\End(\ki_{n,Z})} et\Nligne
\m{\Ext^1_{\ko_n}(\ki_{n,Z},\ki_{n,Z})} . 
Cela permet d'obtenir, si \m{C_n} est plong\'ee dans $S$ (o\`u $S$ est le plan
projectif \m{\P_2} ou une surface K3), des informations sur les d\'eformations
\m{\ki_{n,Z}} en tant que faisceau sur $S$. Soit \m{(\ke_t)_{t\in T}} une telle
d\'eformation, \m{t_0\in T} l'origine, tel que \m{\ke_{t_0}=\ki_{n,Z}}. On a
alors une application canonique
\[\Theta_{t_0}^\ke:T_{t_0}T\lra H^0(\ko_S(nC))/\langle\sigma^n\rangle \]
(o\`u \m{\sigma\in H^0(\ko_S(C))} est une \'equation de $C$)
qui est en gros l'application tangente du morphisme associant \`a un point $t$
de $T$ la courbe support de \m{\ke_t}. On d\'emontre en \ref{def_id2} le

\sepprop

{\bf Th\'eor\`eme : } {\em Si $\ke$ est une d\'eformation compl\`ete de
\m{\ki_{n,Z}}, l'application \m{\Theta^\ke_{t_0}} induit
un isomorphisme
\[\Ext^1_{\ko_S}(\ki_{n,Z},\ki_{n,Z})/
\Ext^1_{\ko_{C_n}}(\ki_{n,Z},\ki_{n,Z})\ \simeq \ V/\langle\sigma^n\rangle ,\]
o\`u \m{V\subset H^0(\ko_S(nC))} est l'espace des courbes passant par tous
les points de $Z$.}

\sepprop

Il est bien entendu possible d'obtenir \m{\ki_{n,Z}} comme limite de
fibr\'es en droites sur des courbes lisses de $S$ ne passant pas par $Z$. Le
r\'esultat pr\'ec\'edent indique que dans ce cas les courbes lisses en question
convergeront vers \m{C_n} avec une multiplicit\'e $>1$ de mani\`ere \`a annuler
l'application tangente \m{\Theta_{t_0}^\ke}.
\end{subsub}

\sepsubsub

\begin{subsub}{\bf Le cas des courbes doubles - }\rm On suppose maintenant que
\m{n=2}. Les seuls faisceaux sans torsion sur \m{C_2} de rang g\'en\'eralis\'e
1 sont les fibr\'es en droites sur $C$.
Les exemples les plus simples de faisceaux sans torsion de rang
g\'en\'eralis\'e 2 sur \m{C_n} sont ceux qui sont quasi localement libres. Ils
sont de deux sortes :
les fibr\'es vectoriels de rang 2 sur $C$ et les fibr\'es en droites sur
\m{C_2} .

On montre en \ref{XXX2} qu'un faisceau sans torsion de rang g\'en\'eralis\'e
2 qui n'est pas quasi localement libre est isomorphe \`a un faisceau du type
\m{\ki_Z\ot\kl}, o\`u \m{Z\subset C} est un sous-sch\'ema fini, \m{\ki_Z} son
faisceau d'id\'eaux sur \m{C_2} et $\kl$ un fibr\'e en droites sur \m{C_n}. Le
sous sch\'ema $Z$ est unique, mais pas le fibr\'e en droites $\kl$.
\end{subsub}

\sepsubsub

{\bf D\'eformations des faisceaux de rang g\'en\'eralis\'e 2 - }

Les fibr\'es en droites sur \m{C_2} ne peuvent \'evidemment se d\'eformer qu'en
fibr\'es en droites. Seuls les fibr\'es vectoriels de rang 2 sur $C$ peuvent se
d\'eformer en faisceaux sans torsion de rang g\'en\'eralis\'e 2 d'un autre
type. On d\'emontre dans le th\'eor\`eme \ref{DUAL_0c} qu'un fibr\'e vectoriel
de rang 2 sur $C$ se d\'eforme en faisceau non concentr\'e sur $C$ si et
seulement si il contient un fibr\'e en droites de degr\'e suffisamment
\'elev\'e.

Un fibr\'e vectoriel de rang 2 sur $C$ est donc la limite d'un nombre fini
(\'eventuellement nul) de familles de faisceaux d'index positif.

Dans ce qui suit on suppose pour simplifier les notations que \m{C_2} est
plong\'ee dans une surface lisse $S$. Soit $E$ un fibr\'e vectoriel sur $C$. On 
montre en \ref{DUAL_0} que l'application canonique
\[\Ext^1_{\ko_2}(E,E)\lra\Ext^1_{\ko_S}(E,E)\]
est un isomorphisme.

Soit \m{Z\subset S} un sous-sch\'ema fini. On montre de m\^eme en \ref{C2_x6_1}
que l'application canonique
\[\Ext^1_{\ko_2}(\ko_Z,\ko_Z)\lra\Ext^1_{\ko_S}(\ko_Z,\ko_Z)\]
est un isomorphisme. Les d\'eformations de \m{\ko_Z} en tant que faisceau sur
\m{C_2} sont les m\^emes que ses d\'eformations en tant que faisceau sur $S$
contenues dans \m{C_2}.

Soit \m{\ke=\ki_Z\ot\kl} un faisceau sans torsion de rang g\'en\'eralis\'e
2 sur \m{C_2}. Alors toutes les d\'eformations de $\ke$ proviennent de
d\'eformations de \m{\ko_Z} et de $\kl$ (cf. \ref{C2_x7d}).
\end{sub}

\sepsub

\Ssect{Faisceaux quasi localement libres de rang g\'en\'eralis\'e 3 sur les
courbes doubles}{intro_7}

Soient $S$ une surface projective lisse irr\'eductible et \m{C\subset S} une
courbe projective lisse irr\'eductible. Soient \m{C_2\subset S} la courbe
double associ\'ee, \m{L=\ko_C(-C)} et \m{l=-\deg(L)}. On suppose que
\m{l=C^2\geq 1}. Le genre de $C$ est \ \m{g=\frac{1}{2}(C^2+K_SC)+1} .
On \'etudie dans \ref{RGEN3_2} les vari\'et\'es de modules de faisceaux quasi
localement libres de rang g\'en\'eralis\'e 3 sur \m{C_2}.

Si $\ke$ est un faiseau quasi localement libre de rang g\'en\'eralis\'e 3 sur
$S$, la premi\`ere filtration canonique de $\ke$ donne une suite exacte \
\m{0\to E_\ke\to\ke\to F_\ke\to 0}, o\`u \m{E_\ke}, \m{F_\ke} sont des fibr\'es
vectoriels sur $C$, \m{E_\ke} \'etant de rang 1 et \m{F_\ke} de rang 2. Soient
\m{\Phi_\ke:F_\ke\ot L\to E_\ke} le morphisme canonique surjectif, et
\m{\Gamma_\ke=\ker(\Phi_\ke)\ot L^*} . On a
aussi une autre suite exacte \ \m{0\to G_\ke\to\ke\to E_\ke\ot L^*\to 0}, o\`u
\m{G_\ke} est le plus grand sous-faisceau de $\ke$ de support $C$, qui est
localement libre de rang 2 sur $C$. On montre que les degr\'es de \m{E_\ke},
\m{F_\ke}, \m{G_\ke}, \m{\Gamma_\ke} sont invariants par d\'eformation. 

Soient $\gamma$, $\epsilon$ des entiers tels que \m{\gamma-l<\epsilon<\gamma} .
Soit \m{\km_s(\epsilon,\gamma)} l'ensemble des faisceaux $\ke$ tels que
\m{E_\ke} soit de degr\'e $\epsilon$, \m{\Gamma_\ke} de degr\'e $\gamma$, et que
\m{F_\ke}, \m{G_\ke} soient stables. On montre que \m{\km_s(\epsilon,\gamma)}
est un ouvert non vide de la vari\'et\'e de modules \m{M(3,2\epsilon+\gamma+l)}
des faisceaux semi-stables (au sens de Simpson) de rang g\'en\'eralis\'e 3 et
de degr\'e g\'en\'eralis\'e \m{2\epsilon+\gamma+l} sur \m{C_2}. On a \
\m{\dim(\km_s(\epsilon,\gamma))=5g+2l-4} . 

L'adh\'erence \m{\ov{\km_s}(\epsilon,\gamma)} de \m{\km_s(\epsilon,\gamma)} dans
\m{M(3,2\epsilon+\gamma+l)} est une composante irr\'eductible de cette
derni\`ere. Il existe une autre composante \m{M_C(3,2\epsilon+\gamma+l)}, celle
qui est constitu\'ee des fibr\'es semi-stables sur $C$. Les vari\'et\'es
\m{\ov{\km_s}(\epsilon,\gamma)} rencontrent toutes
\m{M_C(3,2\epsilon+\gamma+l)}. J'ignore si elles sont disjointes.

La vari\'et\'e \m{\km_s(\epsilon,\gamma)} est non r\'eduite. La vari\'et\'e
r\'eduite sous-jacente \m{\km_s^{red}(\epsilon,\gamma)} est lisse. En tout point
$\ke$ le conoyau de l'application canonique \
\m{T\km_s^{red}(\epsilon,\gamma)_\ke\to T\km_s(\epsilon,\gamma)_\ke} \ est
canoniquement isomorphe \`a \m{H^0(L^*)}. 

Pour des faisceaux de rang plus \'elev\'e, la situation est plus compliqu\'ee,
car les rangs et degr\'es des gradu\'es de la premi\`ere filtration canonique ne
sont pas invariants par d\'eformation.
\end{sub}

\newpage

\section{Pr\'eliminaires}\label{prelim}

\Ssect{Courbes multiples}{cour_mul}

(cf. \cite{ba_fo}, \cite{va}, \cite{man}, \cite{be_fr}).

Soit $X$ une vari\'et\'e alg\'ebrique lisse connexe de dimension 3, et
\m{C\subset X} une courbe lisse connexe. On appelle {\em courbe
multiple de support $C$} un sous-sch\'ema de Cohen-Macaulay
\m{ Y\subset X} tel que l'ensemble des points ferm\'es de $Y$ soit $C$.
Autrement dit, \m{Y_{red}=C}.

Soit $n$ le plus petit entier tel que \m{Y\subset C^{(n-1)}}, \m{C^{(k-1)}}
d\'esignant le $k$-i\`eme voisinage infinit\'esimal de $C$, c'est-\`a-dire \
\m{\ki_{C^{(k-1)}}=\ki_C^{k}} .
On a une filtration \ \m{C=C_1\subset C_2\subset\cdots\subset C_{n}=Y} \
o\`u $C_i$ est le plus grand sous-sch\'ema de Cohen-Macaulay contenu dans
\m{Y\cap C^{(i-1)}}. On appelle $n$ la {\em multiplicit\'e} de $Y$.

On dit que $Y$ est {\em primitive} si  pour tout point ferm\'e $x$ de $C$,
il existe une surface $S$ de $X$ contenant un voisinage de $x$ dans $Y$ et
lisse en $x$. Dans ce cas, \m{L=\ki_C/\ki_{C_2}} est un fibr\'e en droites sur
$C$ et on a \ \m{\ki_{C_{j}}/\ki_{C_{j+1}}=L^j} \ pour \m{1\leq j<n}.
Soit \m{P\in C}. Alors il existe des \'el\'ements $x$, $y$, $t$ de
\m{m_{X,P}} (l'id\'eal maximal de \m{\ko_{X,P}}) dont les images dans
 \m{m_{X,P}/m_{X,P}^2} forment une base, et que pour \m{1\leq i<n} on ait
\ \m{\ki_{C_i,P}=(x,y^{i})} .

Le cas le plus simple est celui o\`u $Y$ est contenue dans une surface lisse de
$X$. Dans ce cas il est m\^eme inutile de mentionner la vari\'et\'e ambiente
$X$, et on peut voir une courbe primitive de Cohen-Macaulay comme une courbe
multiple qui est une sous-vari\'et\'e ferm\'ee d'une surface lisse.

Soient $S$ une surface lisse, \m{Y\subset S} une courbe primitive de
multiplicit\'e $n$ et $C$ la courbe r\'eduite associ\'ee. Soient \m{P\in C}
et \m{f\in\ko_{S,P}} une \'equation locale de $C$. Alors on a
\ \m{\ki_{C_i,P}=(f^{i})} \
pour \m{0\leq j<n}, en particulier \m{I_{Y,P}=(f^n)}, et
\ \m{L=\ko_C(-C)} .

\end{sub}

\sepsub

\Ssect{Les $\Ext$ de faisceaux d\'efinis sur des sous-vari\'et\'es}{ext_sv}

Soient $X$ une vari\'et\'e projective et \m{Y\subset X} une sous-vari\'et\'e
ferm\'ee. Si \ \m{j:Y\to X} \ est l'inclusion et $E$ un faisceau coh\'erent sur
$Y$, on notera aussi souvent $E$ le faisceau \m{j_*(E)} sur $X$.

\sepprop

\begin{subsub}{\bf Proposition : }\label{pr3} Soient $E$, $F$ des faisceaux
coh\'erents sur $Y$. Alors on a une suite exacte canonique
\[0\lra\Ext^1_{\ko_Y}(F,E)\lra\Ext^1_{\ko_X}(F,E)\lra\Hom(
\Tor^1_{\ko_X}(F,\ko_Y),E)\lra\Ext^2_{\ko_Y}(F,E)  .\]
\end{subsub}

\begin{proof} Soient $\ko_X(1)$ un fibr\'e en droites tr\`es ample sur $X$ et
\m{\ko_Y(1)} sa restriction \`a $Y$. Soient \m{n_0>0} un entier tel que
\m{F(n_0)} soit engendr\'e par ses sections globales et \m{M_0=H^0(F(n_0))}. 
Soient \m{F_0} le noyau du morphisme canonique surjectif \ \m{\ko_X(-n_0)\ot M_0
\to F} , \m{n_1} un entier tel que \m{F_1(n_1)} soit engendr\'e par ses sections
globales et \m{M_1=H^0(F_0(n_1))}. En continuant ce proc\'ed\'e
on obtient la r\'esolution localement libre de $F$ 
\xmat{
\cdots\ko_X(-n_2)\ot M_2\ar[r]^-{f_2} & \ko_X(-n_1)\ot M_1\ar[r]^-{f_1} &
\ko_X(-n_0)\ot M_0\ar[r]^-{f_0} & F\ar[r] & 0
}
En restreignant cette r\'esolution \`a $Y$ on obtient le complexe
\xmat{
\cdots \ko_Y(-n_2)\ot M_2\ar[r]^-{f_{Y2}} &\ko_Y(-n_1)\ot M_1\ar[r]^-{f_{Y1}} &
\ko_Y(-n_0)\ot M_0\ar[r]^-{f_{Y0}} & F\ar[r] & 0
}
Pour tout $i\geq 0$, soient \m{\kx_i=\ker(f_i)}, \m{\ky_i=\ker(f_{Yi})} . On a
une suite exacte
\xmat{
0\ar[r] & \kx_0\ar[r] & \ko_X(-n_0)\ot M_0\ar[r]^-{f_0} & F\ar[r] & 0
}
d'o\`u on d\'eduit la suivante
\xmat{
0\ar[r] & \Tor^1_{\ko_X}(F,\ko_Y)\ar[r] & \kx_{0\mid Y}\ar[r]^-\alpha &
\ko_Y(-n_0)\ot M_0\ar[r]^-{f_{0Y}} & F\ar[r] & 0
}
On a donc une suite exacte
\[ (*) \ \ \ \ \ \
\xymatrix{
0\ar[r] & \Tor^1_{\ko_X}(F,\ko_Y)\ar[r] & \kx_{0\mid Y}\ar[r] & \ky_0\ar[r] & 0
}\]
Posons \ \m{M'_2=H^0(\ky_1(n_2))} . Rappelons que \
\m{M_2=H^0(\kx_1(n_2))} .
Le morphisme \Nligne
 \m{\kx_{1\mid Y}\to\ko_Y(-n_1)\ot M_1} \ d\'eduit de l'inclusion
\ \m{\kx_1\subset\ko_X(-n_1)\ot M_1} \ est \`a valeurs dans \m{\ky_1}.
On en d\'eduit une application lin\'eaire naturelle \
\m{\Phi : M_2\to M'_2} .
On note \Nligne
 \m{f'_2:\ko_Y(-n_2)\ot M'_2\to\ko_Y(-n_1)\ot M_1} \ la compos\'ee
\xmat{
\ko_Y(-n_2)\ot M'_2\ar[r] & \ky_1\flinc[r] & \ko_Y(-n_1)\ot M_1
}
de l'inclusion et de l'\'evaluation. Alors on a un diagramme commutatif avec
lignes exactes
\xmat{
\ko_X(-n_2)\ot M_2\ar[r]^-{f_2}\ar[d]^{r_2\ot\Phi} &
\ko_X(-n_1)\ot M_1\ar[r]^-{f_1}\ar[d]^{r_1\ot I_{M_1}} &
\ko_X(-n_0)\ot M_0\ar[r]^-{f_0}\ar[d]^{r_0\ot I_{M_0}} & F\ar[r]\fleq[d] & 0 \\
\ko_Y(-n_2)\ot M'_2\ar[r]^-{f'_2} & \ko_Y(-n_1)\ot M_1\ar[r]^-{f_{1Y}} &
\ko_Y(-n_0)\ot M_0\ar[r]^-{f_{0Y}} & F\ar[r] & 0
}
o\`u pour $i=0,1,2$, \m{r_i:\ko_X(-n_i)\to\ko_Y(-n_i)} est le morphisme
canonique.
d'o\`u le diagramme commutatif
\xmat{
\Hom(\ko_Y(-n_0)\ot M_0,E)\ar[r]^-{F_{0}}\fleq[d] &
\Hom(\ko_Y(-n_1)\ot M_1,E)\ar[r]^-{F_{1Y}}\fleq[d] &
\Hom(\ko_Y(-n_2)\ot M'_2,E)\ar[d]^{{}^t\Phi\ot I} \\
\Hom(\ko_X(-n_0)\ot M_0,E)\ar[r]^-{F_{0}} &
\Hom(\ko_X(-n_1)\ot M_1,E)\ar[r]^-{F_{1}} & \Hom(\ko_X(-n_2)\ot M_2,E)
}
On peut supposer que $n_0$ est suffisamment grand pour que
\m{\Ext^1_{\ko_Y}(F,E)} puisse se calculer \`a l'aide de la suite exacte du
haut, c'est \`a dire que \
\m{\Ext^1_{\ko_Y}(F,E)\simeq\ker(F_{1Y})/\imm(F_0)} .
On d\'eduit du diagramme pr\'ec\'edent une application lin\'eaire
\[\Theta : \Ext^1_{\ko_Y}(F,E)=\ker(F_{1Y})/\imm(F_0)\lra
\Ext^1_{\ko_X}(F,E)=\ker(F_{1})/\imm(F_0) . \]
Montrons que $\Theta$ est injective. Soit \ \m{\alpha:\ko_Y(-n_1)\ot M_1\to E}
\ s'annulant sur \m{\imm(f'_2)}, tel que le morphisme induit \
\m{\ov{\alpha}:\ko_X(-n_1)\ot M_1\to E} \ ait pour image $0$ dans
\m{\Ext^1_{\ko_X}(F,E)} . Cela signifie que $\ov{\alpha}$ se factorise par
$f_1$ :
\xmat{
\ov{\alpha}:\ko_X(-n_1)\ot M_1\ar[r]^-{f_1} &
\ko_X(-n_0)\ot M_0\ar[r]^-\gamma & E
}
Par restriction \`a $Y$ on obtient la factorisation
\xmat{
\alpha:\ko_Y(-n_1)\ot M_1\ar[r]^-{f_{Y1}} &
\ko_Y(-n_0)\ot M_0\ar[r]^-{\gamma_{\mid Y}} & E
}
qui prouve que l'image de $\alpha$ dans \m{\Ext^1_{\ko_Y}(F,E)} est $0$.

Par d\'efinition de $\Theta$ on a \
\m{\coker(\Theta)\simeq\ker(F_1)/\ker(F_{1Y})} .
On a
\ \m{\ker(F_{1Y})=\Hom(\ky_0,E)} , \m{\ker(F_{1})=\Hom(\kx_{0\mid Y},E)} .
D'apr\`es la suite exacte $(*)$ on a donc une suite exacte
\[0\lra\coker(\Theta)\lra\Hom(\Tor^1_{\ko_X}(F,\ko_Y),E)\lra
\Ext^1_{\ko_Y}(\ky_0,E)  .\]
D'apr\`es le choix de $n_0$ et la suite exacte \
\m{0\to\ky_0\to\ko_Y(-n_0)\ot M_0\to F\to 0} \
on voit que \ \m{\Ext^1(\ky_0,E)\simeq\Ext^2_{\ko_Y}(F,E)} .
Le fait que la suite exacte de la proposition \ref{pr3} ne d\'epend pas du
choix de la r\'esolution localement libre de $F$ est laiss\'e au lecteur.
\end{proof}
\end{sub}

\sepsub

\Ssect{Vari\'et\'es de Brill-Noether}{BN}

Les r\'esultats qui suivent seront utilis\'es en \ref{GEN_FIB2}. Soit $C$ une
courbe projective irr\'eductible lisse de genre \m{g\geq 2}. Pour tout entier
$d$ on note \m{J^d} la jacobienne des fibr\'es en droites de degr\'e $d$ sur $C$
et \m{M_s(2,d)} la vari\'et\'e de modules des fibr\'es stables de rang 2 et de
degr\'e $d$ sur $C$, qui est une vari\'et\'e alg\'ebrique irr\'eductible lisse
de dimension \m{4g-3}. Soit
\[W_0^{2,d} \ = \ \lbrace E\in M_s(2,d) ; h^0(E)>0\rbrace ,\]
qui est une sous-vari\'et\'e ferm\'ee de \m{M_s(2,d)}. Si \m{d\leq 0} on a
\m{W_0^{2,d}=\emptyset}, et si \m{d\geq 2g-2} on a \m{W_0^{2,d}=M_s(2,d)} .
Si \m{0<d\leq 2g-2}, \m{W_0^{2,d}} est non vide, irr\'eductible et de dimension
\m{2g-2+d}. Soit \m{G_0^{2,d}} la vari\'et\'e des {\em paires de Brill-Noether}
\m{(E,s)}, o\`u \m{E\in M_s(2,d)} et $s$ est une droite de \m{H^0(E)}. Si
\m{0<d\leq 2g-2}, c'est une vari\'et\'e lisse irr\'eductible de dimension
\m{2g-2+d} (cf. \cite{BMNO}, \cite{ra-vi}, \cite{sun}).

On note \ \m{det_d:M_s(2,d)\to J^d} \ le morphisme associant \`a un fibr\'e
stable $E$ son d\'eterminant.

\sepprop

\begin{subsub}\label{BN1}{\bf Proposition : } Si \m{d>0}, on a \
\m{det_d(W_0^{2,d})=J^d} .
\end{subsub}

\begin{proof}
Il suffit de traiter les cas \m{d=1,2}. En effet, il existe un entier
\m{\delta>0} tel que \m{d-2\delta} soit \'egal \`a 1 ou 2. Soit \m{U\in
J^\delta} tel que \m{h^0(U)>0} . Alors \m{W_0^{2,d}} contient tous les
fibr\'es du type \m{E\ot U}, avec \m{E\in W_0^{2,d-2\delta}} . Donc
\m{det_d(W_0^{2,d})} contient tous les \m{det_{d-2\delta}(E)\ot U^2}. Si
\m{det_{d-2\delta}(W_0^{2,d-2\delta})=J^{d-2\delta}} on a donc
\m{det_d(W_0^{2,d})=J^d} .

{\em Le cas \m{d=1} - } Soit \m{D\in J^1} . On va montrer que \m{D\in
\det_1(W_0^{2,1})}. On a \m{h^1(D^*)=g>0}, donc
il existe une extension non triviale \ \m{0\to\ko_C\to E\to D\to 0} . Alors $E$
est stable : il suffit de montrer que si \m{U\in J^1}, alors \m{\Hom(U,E)=\nsp}.
Supposons que \m{\Hom(U,E)\not=\nsp}. On a \ \m{\Hom(U,E)\subset\Hom(U,D)}, donc
\m{\Hom(U,D)\not=\nsp}, donc \m{U=D}, et il existe une section du morphisme
\m{E\to D}, ce qui est absurde car l'extension est non triviale. On a
\m{det(E)=D}, \m{E\in W_0^{2,1}}, donc \m{D\in\det_1(W_0^{2,1})}.

{\em Le cas \m{d=2} - } Soit \m{D\in J^2} . On va montrer que \m{D\in
\det_2(W_0^{2,2})}. On a \m{h^1(D^*)=g+1>2}. Soit \ \m{0\to\ko_C\to E\to D\to 0}
\ une extension non triviale. Supposons que $E$ ne soit pas stable. Alors il
existe \m{U\in J^1} tel que \m{\Hom(U,E)\not=\nsp}. On a donc
\m{\Hom(U,D)\not=\nsp}, et il existe donc un point $x$ de $C$ tel que
\m{U=D(-x)}. De plus, l'application canonique \ \m{\Hom(D(-x),D)=\C\to
\Ext^1_{\ko_C}(D(-x),\ko_C)} \ est nulle. Soit \m{\sigma\in H^1(D^*)}
correspondant \`a l'extension pr\'ec\'edente. On a un carr\'e commutatif
\xmat{\End(D)\ar[rr]\ar[d]^\simeq & & H^1(D^*)\ar[d]^{\alpha_x}\\
\Hom(D(-x),D)\ar[rr] & & H^1(D^*(x))}
o\`u les fl\`eches horizontales proviennent de $\sigma$ et les verticales de la
section de \m{\ko_C(x)}. On en d\'eduit que \m{\sigma\in\ker(\alpha_x)}. Puisque
\m{\dim(\ker(\alpha_x))=1} et \m{h^1(D^*)>2} il existe  \m{\sigma_0\in H^1(D^*)}
tel que pour tout \m{x\in C}, \m{\sigma_0\not\in\ker(\alpha_x)}. Si \
\m{0\to\ko_C\to E_0\to D\to 0} \ est l'extension correspondante, le fibr\'e
\m{E_0} est stable d'apr\`es ce qui pr\'ec\`ede. On a \m{det(E_0)=D},
\m{E_0\in W_0^{2,2}}, donc \m{D\in\det_2(W_0^{2,2})}.
\end{proof}

\end{sub}

\sepsec

\section{Fibr\'es vectoriels sur les courbes multiples}

Soient $C_n$ une courbe projective multiple primitive de multiplicit\'e $n>1$,
\Nligne \m{C=C_1\subset C_2\subset\cdots\subset C_n=Y} \
la filtration canonique et  $L$ le fibr\'e en droites sur $C$ associ\'e
(cf. \ref{cour_mul}).
On pose, pour \m{1\leq i\leq n}, \m{\ko_i=\ko_{C_i}} .

\sepsub

\Ssect{Prolongements de fibr\'es vectoriels}{Extens}

\sepprop

\begin{subsub}\label{pr5}{\bf Th\'eor\`eme :} Si \m{1\leq i\leq n},
tout fibr\'e vectoriel sur \m{C_i} peut se prolonger en fibr\'e vectoriel sur
\m{C_n}.
\end{subsub}

\begin{proof}
On peut en raisonnant par r\'ecurrence supposer que \m{i=n-1}. Soit $\ki$
le faisceau d'id\'eaux de \m{C_{n-1}} dans \m{C_n}. On a un
isomorphisme de \m{\ko_{n}}-modules :
\ \m{\ki\simeq\ki_{\mid C_{n-1}}}.

Soit $E$ un fibr\'e vectoriel sur \m{C_{n-1}}, de rang $r$. Il existe un
recouvrement ouvert \m{(U_i)} de $C_{n-1}$ tel que chaque restriction
\m{E_{\mid U_i}} soit un fibr\'e trivial.
Soient \ \m{\lambda_i:E_{\mid U_i}\simeq\ko_{n-1}(U_i)\ot\C^r} \ des
trivialisations, \
\m{\lambda_{ij}=\lambda_j\circ\lambda_i^{-1}\in\GL(r,\ko_{n-1}(U_{ij}))} .
Soient \m{\Lambda_{ij}\in\GL(r,\ko_{n}(U_{ij}))} une extension de
\m{\lambda_{ij}} et \
\m{\rho_{ijk}=\Lambda_{jk}\Lambda_{ij}-\Lambda_{ik}} \
(\'el\'ement de \m{\ko_{n}(U_{ijk})\ot\End(\C^r)}). Alors les
\m{\Lambda_{ij}} d\'efinissent un fibr\'e vectoriel sur \m{C_n} si et
seulement si les \m{\rho_{ijk}} sont nuls. Leurs restrictions \`a \m{C_{n-1}}
sont nulles, donc on peut les consid\'erer comme des \'el\'ements de
\[\End(\ko_{n}(U_{ijk})\ot\C^r)\ot\ki = \
 \HHom(\ko_{n-1}\ot\C^r,\ki\ot\C^r)(U_{ijk}) \ .\]
Soit \
\m{\mu_{ijk}=(\lambda_k)^{-1}\rho_{ijk}\lambda_i} ,
qui est un \'el\'ement de \m{\HHom(E,E\ot\ki)(U_{ijk})} .
Pour obtenir une extension de $E$ \`a \m{C_n} on peut remplacer les
\m{\Lambda_{ij}} par \
\m{\Lambda'_{ij}=\Lambda_{ij}-\beta_{ij}} ,
avec \m{\beta_{ij}} nul sur \m{C_{n-1}}~. On peut donc consid\'erer les
\m{\beta_{ij}} comme des \'el\'ements de 
\m{\HHom(\ko_i\ot\C^r,\ki\ot\C^r)(U_{ij})}.
Soit \
\m{\rho'_{ijk}=\Lambda'_{jk}\Lambda'_{ij}-\Lambda'_{ik}} .
Alors on a \
\m{\rho'_{ijk}=\rho_{ijk}-\beta_{jk}\Lambda_{ij}-
\Lambda_{jk}\beta_{ij}+\beta_{ik}} .
Posons \ \m{\alpha_{ij}=(\lambda_j)^{-1}\beta_{ij}\lambda_i}, qui est un
\'el\'ement de \m{\HHom(E,E\ot\ki)(U_{ij})} . Alors on a
\m{\rho'_{ijk}=0} si et seulement si
\[
(*) \ \ \ \ \ \ \mu_{ijk} \ = \ \alpha_{ij}+\alpha_{jk}-\alpha_{ik} .\]
On a \
\m{\Lambda_{kl}\rho_{ijk}-\rho_{ijl}+\rho_{ikl}-\rho_{jkl}\Lambda_{ij}=0} ,
d'o\`u il d\'ecoule que
\m{\mu_{ijk}-\mu_{ijl}+\mu_{ikl}-\mu_{jkl}=0} ,
c'est-\`a-dire que \m{(\mu_{ijk})} est un cocycle associ\'e au fibr\'e
vectoriel \m{\HHom(E,E\ot\ki)} sur $C_{n-1}$ et au recouvrement
\m{(U_i)}. Comme $C_{n-1}$ est une courbe, on a \
\m{H^2(\HHom(E,E\ot\ki))=\nsp} ,
d'o\`u l'existence des \m{\alpha_{ij}} satisfaisant l'\'egalit\'e $(*)$ et des
\m{\beta_{ij}} d\'efinissant le prolongement voulu de $E$.
\end{proof}

\sepsubsub

\begin{subsub}\label{app_ext0}R\'esolutions canoniques - \rm
On suppose que \m{C_n} s'\'etend en une courbe primitive \m{C_{n+1}} de
multiplicit\'e \m{n+1}. Soit \m{\ko_{n+1}=\ko_{C_{n+1}}} .
Le faisceau d'id\'eaux $\L_n$ de $C$ dans \m{C_{n+1}} est un \m{\ko_{n}}-module
localement libre de rang 1. D'apr\`es le th\'eor\`eme \ref{pr5}, on peut le
prolonger en un fibr\'e en droites sur \m{C_{n+1}}, not\'e $\L$.
Si \m{1\leq p\leq n}, \m{\L_n^p} est le faisceau d'id\'eaux de $C_p$ dans
\m{C_{n+1}}, c'est un \m{\ko_{n+1-p}}-module libre de rang 1, et on a
\m{\L_n^p=\L^p_{\mid C_{n+1-p}}} .
Soit \
\m{\mu_p : \L^p\to\ko_{n+1}} \
le morphisme compos\'e \ \m{\L^p\to\L^p_{\mid C_{n+1-p}}=\L_n^p
\subset\ko_{n+1}}. On a alors,
pour \m{1\leq i\leq n} une r\'esolution localement libre canonique de \m{\ko_i}
\xmat{\cdots\ar[r] & \L^{2(n+1)}\ar[rr]^-{\mu_{n+1-i}} & & \L^{n+1+i}
\ar[r]^-{\mu_i} & \L^{n+1}\ar[rr]^-{\mu_{n+1-i}} & & \L^{i}\ar[r]^-{\mu_i} &
\ko_{n+1}\ar[r] & \ko_i }
\end{subsub}

Soit $E_{n+1}$ un fibr\'e vectoriel sur $C_{n+1}$, $E_n$ sa restriction \`a
\m{C_{n}} et $E_C$ sa restriction \`a $C$. On a une suite exacte canonique sur
\m{C_{n+1}}
\[0\lra E_n\ot\L\lra E_{n+1}\lra E_C\lra 0 \ .\]
D'apr\`es la suite spectrale des Ext (cf. (cf. \cite{go}, 7.3), on a une suite
exacte
\[0\lra H^1(\HHom(E_C,E_n\ot\L))\lra\Ext^1_{\ko_{n+1}}(E_C,E_n\ot\L)\lra
H^0(\EExt^1_{\ko_{n+1}}(E_C,E_n\ot\L))\lra 0 .\]

\sepprop

\begin{subsub}\label{lem3_1}{\bf Lemme : }
Les faisceaux \m{\HHom(E_C,E_n\ot\L)} et
\m{\EExt^1_{\ko_{n+1}}(E_C,E_n\ot\L)} sont
de support $C$ et on a des isomorphismes canoniques
\[\HHom(E_C,E_n\ot\L) \ \simeq \ E_C^*\ot E_C\ot L^n , \ \ \ \ \ \
\EExt^1_{\ko_{n+1}}(E_C,E_n\ot\L) \ \simeq \ E_C^*\ot E_C .\]
\end{subsub}

\begin{proof}
De la r\'esolution pr\'ec\'edente de \m{\ko_C=\ko_1} on d\'eduit la
r\'esolution localement libre de \m{E_C}
\xmat{
\cdots\ar[r] & E_{n+1}\ot\L^{n+1}\ar[rr]^-{I_{E_{n+1}}
\ot\mu_{n}} & &
E_{n+1}\ot\L\ar[rr]^-{I_{E_{n+1}}\ot\mu_1} & & E_{n+1}\ar[r] & E_C }
d'o\`u d\'ecoule ais\'ement le lemme.
\end{proof}

\sepprop

On a donc une suite exacte
\xmat{0\ar[r] & H^1(E_C^*\ot E_C\ot L^n)\ar[r] &
\Ext^1_{\ko_{n+1}}(E_C,E_n\ot\L)\ar[r]^-\pi & \End(E_C)\ar[r] & 0 .}

\sepprop

\begin{subsub}\label{lem3}{\bf Corollaire : }
Soit \
\m{0\to E_n\ot\L\to F\to E_C\to 0} \
une extension correspondant \`a \
\m{\sigma\in\Ext^1_{\ko_{n+1}}(E_C,E_n\ot\L)}. Alors $F$ est localement libre
si et seulement si \m{\pi(\sigma)} est un automorphisme.
\end{subsub}

\begin{proof}
L'assertion est locale. Il suffit donc de montrer que si \m{P\in C} et si
\[0\lra\ko_{n-1,P}\ot\C^n\lra M\lra\ko_{C,P}\ot\C^n\lra 0\]
est une suite exacte de \m{\ko_{n,P}}-modules, alors $M$ est libre si et
seulement si l'\'el\'ement associ\'e de \
\m{\Ext^1_{\ko_{n,P}}(\ko_{n-1,P}\ot\C^n,\ko_{C,P}\ot\C^n)
\simeq\End(\C^n)\ot\ko_{C,P}} \ est un automorphisme, ce qui est ais\'e.
\end{proof}

\sepprop

\begin{subsub}\label{extens_c2} Le cas des fibr\'es en droites - \rm
On suppose ici que \m{\deg(L)<0}.
Soient \m{D_n} un fibr\'e en
droites sur \m{C_n}, \m{D_C=D_{n\mid C}}. On a alors une suite exacte
\xmat{0\ar[r] & H^1(L^n)\ar[r] & \Ext^1_{\ko_{n+1}}(D_C,D_n\ot\L)\ar[r]^-\pi
& \C\ar[r] & 0 \ .}

L'ensemble \m{P_{D_n}} des
prolongements de \m{D_n} \`a \m{C_{n+1}} s'identifie \`a \m{\pi^{-1}(1)},
espace affine isomorphe \`a \m{H^1(L^n)}.
En particulier, consid\'erons \m{P_{\ko_n}}, le groupe des fibr\'es en droites
sur \m{C_{n+1}} donc la restriction \`a \m{C_n} est le fibr\'e trivial. La
bijection \
\m{P_{\ko_n}\simeq H^1(L^n)} \
d\'efinie par le prolongement \m{\ko_{n+1}} de \m{\ko_{n}} est un morphisme de
groupes ab\'eliens (on peut le voir par exemple en utilisant la
d\'emonstration du th\'eor\`eme \ref{pr5}).
\end{subsub}

\end{sub}

\sepsub

\Ssect{Groupe de Picard de $C_n$}{piccn}

Soit $X$ une vari\'et\'e alg\'ebrique. Rappelons qu'on appelle {\em groupe de
Picard de $X$} une vari\'et\'e alg\'ebrique, not\'ee habituellement
\m{\Pic(X)},
munie d'un fibr\'e en droites $\kl$ sur \ \m{\Pic(X)\times X}, appel\'e {\em
fibr\'e de Poincar\'e},
tels que pour toute famille plate $\kd$ de fibr\'es en droites sur
$X$, param\'etr\'ee par une vari\'et\'e alg\'ebrique $T$, il existe un unique
morphisme \m{f_\kd:T\to\Pic(X)} et un fibr\'e en droites $L$ sur $T$ tel que
\[\kd \ \simeq \ (f_\kd\times I_X)^*(\kl)\ot p_T^*(L) , \]
($p_T$ d\'esignant la projection \m{T\times X\to T}). S'il existe, \m{\Pic(X)}
est unique \`a isomorphisme pr\`es, et $\kl$ est unique \`a un fibr\'e en
droites sur \m{\Pic(X)} pr\`es. L'ensemble des points ferm\'es de \m{\Pic(X)}
s'identifie avec celui des classes d'isomorphisme de fibr\'es en droites sur
$X$.

\sepprop

\begin{subsub}{\bf Th\'eor\`eme : }\label{theo_pic} On suppose que 
\m{\deg(L)<0}. Soit \m{n\geq 1} un
entier. Alors il existe un groupe de Picard pour \m{C_n}.
\end{subsub}

\begin{proof}
On construit \m{\Pic(C_n)} et le fibr\'e de Poincar\'e \m{\kl_n} sur \
\m{\Pic(C_n)\times C_n} \ par r\'ecurrence sur $n$. Ils sont \'evidemment bien
connus pour \m{n=1}. Supposons que \m{n>1} et que \m{\Pic(C_{n-1})} et
\m{\kl_{n-1}} existent. On peut
voir ce dernier comme un faisceau coh\'erent sur \m{\Pic(C_{n-1})\times C_n}.
On note \ \m{p:\Pic(C_{n-1})\times C_n\to C_n},\
\m{q:\Pic(C_{n-1})\times S\to\Pic(C_{n-1})} \ les projections.

Soient \m{\L_{n -1}} le faisceau d'id\'eaux de \m{C} dans \m{C_n},
et $\L$ un prolongement de \m{\L_{n-1}} \`a
\m{C_n}. Le faisceau \m{\Ext^1} relatif sur \m{\Pic(C_{n-1})},
\m{\ke=\Ext^1_q(\kl_{n-1}\ot p^*(\ko_C),\kl_{n-1}\ot p^*(\L))} \
est localement libre. Pour tout \m{x\in\Pic(C_{n-1})} on a \
\m{\ke_x=\Ext^1_{\ko_n}(\kl_{n-1,x\mid C},\kl_{n-1,x}\ot\L)} .
On a vu qu'on avait une suite exacte canonique
\xmat{0\ar[r] & H^1(L^{n-1})\ar[r] &
\Ext^1_{\ko_n}(\kl_{n-1,x\mid C},\kl_{n-1,x}\ot\L)\ar[r]^-{\pi_x} &
\End(\kl_{n-1,x})=\C\ar[r] & 0 \ .}
On a donc un morphisme surjectif \
\m{\ke\to\ko_{\Pic(C_{n-1})}} .
On prend \ \m{\Pic(C_n)=\pi^{-1}(1)} et \m{\kl_n} est une extension
universelle. Les d\'etails de la d\'emonstration sont laiss\'es au lecteur.
\end{proof}
\end{sub}

\sepsub

\Ssect{Description de $\Pic(C_n)$}{piccn_desc}

On suppose que \m{\deg(L)<0}.
Comme dans le cas des courbes lisses, les composantes irr\'eductibles de
\m{\Pic(C_n)} sont d\'efinies par le degr\'e : ce sont les vari\'et\'es
\m{\Pic^i(C_n)}, \m{i\in\Z} des fibr\'es en droites sur \m{C_n} dont la
restriction \`a $C$ est de degr\'e $i$.

Rappelons que pour tout fibr\'e en droites \m{D_{n-1}} sur \m{C_{n-1}},
\m{P_{D_{n-1}}} d\'esigne l'ensemble des prolongements de \m{D_{n-1}} \`a
\m{C_n}. On a une suite exacte de groupes ab\'eliens
\[0\lra P_{\ko_{n-1}}\simeq H^1(L^{n-1})\lra\Pic(C_n)\lra
\Pic(C_{n-1})\lra 0\]
(cf. \ref{extens_c2}). Soit \m{{\bf P}_n\subset \Pic(C_n)} le sous-groupe des
fibr\'es en droites dont la restriction \`a $C$ est le fibr\'e trivial. On a
une filtration \
\m{0=G_0\subset G_1\subset\cdots\subset G_{n-1}={\bf P}_n} \
telle que pour \m{1\leq i\leq n-1}, on ait \m{G_i/G_{i-1}\simeq
H^1(L^i)} (\m{G_i} est le groupe des fibr\'es en droites dont la restriction
\`a \m{C_i} est triviale). Il d\'ecoule de \cite{se}, chap. VII, 7, corollaire,
que \m{{\bf P}_n} est isomorphe \`a un produit de groupes \m{\G_a},
c'est-\`a-dire un $\C$-espace vectoriel de dimension finie.

On a donc une suite exacte
\[0\lra{\bf P}_n\lra\Pic(C_n)\lra\Pic(C)\lra 0\]
et on obtient un r\'esultat analogue \`a ce que l'on observe lorsqu'on calcule
le groupe de Picard d'une courbe int\`egre en fonction de celui de sa
normalis\'ee (cf. \cite{ro}, \cite{se}).
\end{sub}

\sepsec

\section{Invariants des faisceaux coh\'erents sur les courbes multiples}
\label{QLL0}

Soient $C_n$ une courbe projective multiple primitive de multiplicit\'e $n>1$
et \Nligne \m{C=C_1\subset C_2\subset\cdots\subset C_n} \
sa filtration canonique. Soient $L$ le fibr\'e en droites associ\'e sur $C$,
\m{\ki_C} le faisceau d'id\'eaux de $C$ dans \m{C_n} et \m{g_C} le genre de $C$.

On pose, pour \m{1\leq i\leq n}, \m{\ko_i=\ko_{C_i}}, et \m{\ko_0=0}, ce sont
des faisceaux coh\'erents sur \m{C_n}.

Soient \m{P\in C}, \m{z\in\ko_{n,P}} une \'equation locale de $C$, et
\m{x\in\ko_{n,P}} au dessus d'un g\'en\'erateur de l'id\'eal maximal de
\m{\ko_{C,P}}.

\sepsub

\Ssect{Filtrations canoniques, rang, degr\'e et fonction caract\'eristique}
{QLL-def}

Soient $M$ un \m{\ko_{n,P}}-module de type fini et $\ke$ un faisceau coh\'erent
sur \m{C_n}.

\sepsubsub

\begin{subsub}\label{QLL-def1}Premi\`ere filtration canonique - \rm
On d\'efinit la {\em filtration canonique de $M$} : c'est la filtration
\[M_{n+1}=\nsp\subset M_n\subset\cdots\subset M_{2}\subset M_1=M\]
telle que pour \m{1\leq i\leq n}, \m{M_{i+1}} soit le noyau du morphisme
canonique surjectif \Nligne
\m{M_{i}\to M_{i}\ot_{\ko_{n,P}}\ko_{C,P}} .
On a donc
\[M_{i}/M_{i+1} \ = \ M_{i}\ot_{\ko_{n,P}}\ko_{C,P}, \ \ \ \
M/M_{i+1} \ \simeq \ M\ot_{\ko_{n,P}}\ko_{i,P}, \ \ \ \
M_{i+1} \ = \ z^iM .\]
Le gradu\'e
\[{\rm Gr}(M) \ = \ \bigoplus_{i=1}^nM_i/M_{i+1} \ = \ \bigoplus_{i=1}^n
z^{i-1}M/z^iM\]
est un \m{\ko_{C,P}}-module.
Les propri\'et\'es suivantes sont imm\'ediates : si \m{1<i\leq n}\Nligne
- on a $M_i=\nsp$ si et seulement si $M$ est un
$\ko_{i-1,P}$-module,\Nligne
- $M_i$ est un $\ko_{n+1-i,P}$-module, et sa filtration canonique
est \
\m{\nsp\subset M_n\subset\cdots\subset M_{i+1}\subset M_i} ,\Nligne
- tout morphisme de \m{\ko_{n,P}}-modules envoie la premi\`ere
filtration canonique du premier module sur celle du second. 

On d\'efinit de m\^eme la {\em premi\`ere filtration canonique de $\ke$} :
c'est la filtration
\[\ke_{n+1}=0\subset \ke_n\subset\cdots\subset \ke_{2}\subset \ke_1=\ke\]
telle que pour \m{1\leq i\leq n}, \m{\ke_{i+1}} soit le noyau du morphisme
canonique surjectif \ \m{\ke_i\to\ke_{i\mid C}}.
On a donc \
\m{\ke_{i}/\ke_{i+1}=\ke_{i\mid C}} ,
\m{\ke/\ke_{i+1}=\ke_{\mid C_i}} .
Le gradu\'e
\m{{\rm Gr}(\ke)} est un \m{\ko_{C}}-module.
Les propri\'et\'es suivantes sont imm\'ediates : si \m{1<i\leq n}\Nligne
- on a \m{\ke_i=\ki_C^{i-1}\ke}, et donc \ \m{{\rm
Gr}(\ke)=\bigoplus_{j=0}^{n-1}\ki_C^j\ke/\ki_C^{j+1}\ke} .\Nligne
- on a $\ke_i=0$ si et seulement si $\ke$ est un faisceau sur
$C_{i-1}$,\Nligne
- $\ke_i$ est un faisceau sur $C_{n+1-i}$, et sa filtration canonique
est \
\m{0\subset\ke_n\subset\cdots\subset\ke_{i+1}\subset\ke_i} .\Nligne
- tout morphisme de faisceaux coh\'erents sur \m{C_n} envoie la
premi\`ere filtration canonique du premier sur celle du second.

{\bf Exemples : } 1 - Si \m{\ke=\ko_m}, \m{1\leq m\leq n}, on a \
\m{\ke_i/\ke_{i+1}=\ko_C\ot L^{i-1}} \ pour \m{1\leq i\leq m}, et
\m{\ke_i=0} si \m{i>m}.

2 - Si $\ke$ est le faisceau d'id\'eaux du point \m{P\in C} on a \
\m{\ke_i/\ke_{i+1}=(\ko_C(-P)\ot L^{i-1})\oplus\C_P} \ si \m{1\leq i<n}, et \
\m{\ke_n=\ko_C(-P)\ot L^{n-1}} .
\end{subsub}

\sepsubsub

\begin{subsub}\label{2-fc}Seconde filtration canonique - \rm
On d\'efinit la {\em seconde filtration canonique de $M$} : c'est la filtration
\[M^{(n+1)}=\nsp\subset M^{(n)}\subset\cdots\subset M^{(2)}\subset M^{(1)}=M\]
avec \
\m{M^{(i)} \ = \ \big\lbrace u\in M ; z^{n+1-i}u=0\big\rbrace} .
Si \ \m{M_{n+1}=\nsp\subset M_n\subset\cdots\subset M_{2}\subset M_1=M} \ est
la (premi\`ere) filtration canonique de $M$ on a \ \m{M_i\subset M^{(i)}} \
pour \m{1\leq i\leq n}. Le gradu\'e
\[{\rm Gr}_2(M) \ = \ \bigoplus_{i=1}^nM^{(i)}/M^{(i+1)} \]
est un \m{\ko_{C,P}}-module.
Les propri\'et\'es suivantes sont imm\'ediates : si \m{1<i\leq n}\Nligne
- $M^{(i)}$ est un $\ko_{n+1-i,P}$-module, et sa filtration canonique
est \hfil\break
\m{\nsp\subset M^{(n)}\subset\cdots\subset M^{(i+1)}\subset M^{(i)}}~,
\Nligne
- tout morphisme de \m{\ko_{n,P}}-modules envoie la seconde
filtration canonique du premier module sur celle du second. 

On d\'efinit de m\^eme la {\em seconde filtration canonique de $\ke$} :
\[\ke^{(n+1)}=\nsp\subset \ke^{(n)}\subset\cdots\subset
\ke^{(2)}\subset \ke^{(1)}=\ke  .\]
Le gradu\'e
\m{{\rm Gr}_2(\ke)} est un \m{\ko_{C}}-module.
Les propri\'et\'es suivantes sont imm\'ediates : si \m{1<i\leq n}\Nligne
- $\ke^{(i)}$ est un faisceau sur $C_{n+1-i}$, et sa filtration
canonique est \
\m{0\subset\ke^{(n)}\subset\cdots\subset\ke^{(i+1)}\subset\ke^{(i)}} ,
\Nligne
- tout morphisme de faisceaux coh\'erents sur \m{C_n} envoie la
seconde filtration canonique du premier sur celle du second.

{\bf Exemples : } 1 - Si \m{\ke=\ko_m}, \m{1\leq m\leq n}, on a \
\m{\ke^{(i)}/\ke^{(i+1)}=\ko_C\ot L^{i+m-n-1}} \ pour \m{n-m+1\leq i\leq n}, et
\m{\ke_i=\ko_m} si \m{i\leq n-m+1}.

2 - Si $\ke$ est le faisceau d'id\'eaux du point \m{P\in C} on a \
\m{\ke^{(i)}/\ke^{(i+1)}=L^{i-1}} \ si \m{2\leq i\leq n}, et \
\m{\ke^{(1)}/\ke^{(2)}=\ko_C(-P)} .
\end{subsub}

\sepsubsub

\begin{subsub}Rang - \rm L'entier \
\m{R(M)=rg({\rm Gr}(M))} \
s'appelle le {\em rang g\'en\'eralis\'e} de $M$.

L'entier \
\m{R(\ke)=rg({\rm Gr}(\ke))} \
s'appelle le {\em rang g\'en\'eralis\'e} de $\ke$. On a donc
\m{R(\ke)=R(\ke_P)} pour tout \m{P\in C}.

\medskip

{\bf Exemple : } Si $\ke$ est localement libre, on a \ \m{R(\ke)=n.rg(\ke)} .
Si le support de $\ke$ est contenu dans $C$ (c'est-\`a-dire que $\ke$ est
un faisceau coh\'erent sur $C$ vu comme faisceau sur \m{C_n}), on a \
\m{R(\ke)=rg(\ke_{\mid C})} .
\end{subsub}

\sepsubsub

\begin{subsub}Degr\'e - \rm L'entier \
\m{\Deg(\ke)=\deg({\rm Gr}(\ke))} \
s'appelle le {\em degr\'e de } $\ke$.

Si \m{R(\ke)>0} on pose \m{\mu(\ke)=\Deg(\ke)/R(\ke)} et on appelle ce nombre la
{\em pente} de $\ke$.

{\bf Exemple : } Si $\ke$ est localement libre, on a
\[\Deg(\ke) \ = \ n.\deg(\ke_{\mid C})+\frac{n(n-1)}{2}rg(\ke_{\mid C})\deg(L)
.\]
Si le support de $\ke$ est contenu dans $C$ (c'est-\`a-dire que $\ke$ est
un faisceau coh\'erent sur $C$ vu comme faisceau sur \m{C_n}), on a \
\m{\Deg(\ke)=\deg(\ke_{\mid C})} .
\end{subsub}

\sepsubsub

\begin{subsub}Fonctions caract\'eristiques - \rm Soit $M$ un
\m{\ko_{n,P}}-module de type fini,
\[M_{n+1}=\nsp\subset M_n\subset\cdots\subset M_{2}\subset M_1=M\]
sa filtration canonique. On associe \`a
$M$ sa {\em fonction caract\'eristique}, c'est la fonction lin\'eaire par
morceaux
\[\xymatrix@R=-1pt{
F(M):[0,n]\ar[r] & \N \\
k\in\lbrace 0,\ldots,n\rbrace\fmaps[r] &
R(M_{n+1-k})=\sigg_{i=n+1-k}^{n}rg(M_i/M_{i+1})}\]
Si $\ke$ est un faisceau coh\'erent sur \m{C_n}, on d\'efinit de m\^eme sa
fonction caract\'eristique \m{F(\ke)}.

On d\'efinit de m\^eme la {\em seconde fonction caract\'eristique} de $M$ ou de
\m{\ke}, \`a partir de leur seconde filtration canonique.

Plus g\'en\'eralement \`a toute filtration $\sigma$ de $M$ dont les gradu\'es
sont des \m{\ko_{C,P}}-modules est associ\'ee de
mani\`ere \'evidente une fonction caract\'eristique \m{F_M(\sigma)} .
\end{subsub}

\sepsubsub

\begin{subsub}{\bf Propri\'et\'es des filtrations canoniques et des fonctions
caract\'eristiques - }\label{QLL-prop}\rm Les d\'emonstrations des r\'esultats
suivants sont faciles et laiss\'ees au lecteur.

1 - La fonction caract\'eristique d'un \m{\ko_{n,P}}-module de type fini ou
d'un faisceau coh\'erent sur $C_n$ est une fonction convexe.

2 - La seconde fonction caract\'eristique d'un \m{\ko_{n,P}}-module de type
fini ou d'un faisceau coh\'erent sur $C_n$ est une fonction concave.

3 - {\em Fonction caract\'eristique g\'en\'erique. }
Soient $N$ un entier positif, \m{N=pn+m}, avec \m{0\leq m<n}, et
\m{M_0=p\ko_{n,P}\oplus\ko_{m,P}} (resp.
\m{\ke_0=p\ko_n\oplus\ko_m}), qui est de rang g\'en\'eralis\'e \m{N}.
Soit $M$ un \m{\ko_{n,P}}-module (resp.
$\ke$ un faisceau coh\'erent sur \m{C_n}) de rang g\'en\'eralis\'e \m{N}.
Alors on a \ \m{F(M)\leq F(M_0)} (resp. \m{F(\ke)\leq F(\ke_0)}).

4 - Soient $M$ un \m{\ko_{n,P}}-module de type fini et $\sigma$ une filtration
de $M$ dont les gradu\'es sont des \m{\ko_{C,P}}-modules. Alors on a \
\m{F_M(\sigma)\geq F(M)} ,
avec \'egalit\'e si et seulement si $\sigma$ est la filtration canonique de
$M$.
\end{subsub}
\end{sub}

\sepsub

\Ssect{Th\'eor\`eme de Riemann-Roch g\'en\'eralis\'e}{RR}

\begin{subsub}{\bf Th\'eor\`eme : }\label{RR0}
Soit $\ke$ un faisceau coh\'erent sur \m{C_n}. Alors on a
\[\chi(\ke) \ = \ \Deg(\ke)+R(\ke)(1-g_C) .\]
\end{subsub}

Cela d\'ecoule \'evidemment du th\'eor\`eme de Riemann-Roch habituel sur $C$
et des d\'efinitions de \m{R(\ke)} et \m{\Deg(\ke)}.

\sepprop

\begin{subsub}\label{hilb_pol}Polyn\^ome de Hilbert - \rm
Soient \m{\ko_n(1)} un fibr\'e en droites tr\`es ample sur \m{C_n},
\m{\ko_C(1)} sa restriction \`a $C$
 et \m{\delta=deg(\ko_C(1))}. Si $\ke$ est un faisceau coh\'erent sur
\m{C_n} et $m$ un entier, on pose comme d'habitude \
\m{\ko_n(m)=\ko_n(1)^{\ot m}}, \m{\ke(m)=\ke\ot_{\ko_{C_n}}\ko_n(m)}, et si
$\kf$ est un faisceau coh\'erent sur $C$, \m{\kf(m)=\kf\ot_{\ko_C}\ko_C(m)}.
On a un isomorphisme canonique \
\m{Gr(\ke(m))\simeq Gr(\ke)(m)} ,
d'o\`u on d\'eduit imm\'ediatement \
\m{R(\ke(m)) = R(\ke)} , \m{\Deg(\ke(m))=\Deg(\ke)+R(\ke)\delta m} .
Le polyn\^ome
\[P_\ke(X) \ = \ \Deg(\ke)+R(\ke)(1-g_C) + R(\ke)\delta\cdot X \]
s'appelle le {\em polyn\^ome de Hilbert de $\ke$}. On a donc \ \m{P_\ke(m)=
\chi(\ke(m))} \ pour tout entier $m$.
\end{subsub}

\end{sub}

\sepsub

\Ssect{Propri\'et\'es du rang et du degr\'e g\'en\'eralis\'es}{prop_rang}

Soit $M$ un \m{\ko_{n,P}}-module de type fini. Si la limite
\[\lim_{p\to\infty}\biggl(\frac{1}{p}\dim_\C(M\ot_{\ko_{n,P}}\ko_{n,P}/(x^p))
\biggr)\]
existe et est finie, on note \m{R_0(M)} sa valeur, et on dit que \m{R_0(M)}
{\em est d\'efini}.

\sepprop

\begin{subsub}{\bf Th\'eor\`eme : }\label{R_R_0}
Soit $M$ un \m{\ko_{n,P}}-module de type fini. Alors \m{R_0(M)} est d\'efini
et on a \ \m{R_0(M)=R(M)} .
\end{subsub}

\begin{proof} Traitons d'abord le cas o\`u $M$ est un \m{\ko_{C,P}}-module. On
a alors \ \m{M\simeq T\oplus r\ko_{C,P}}, $T$ \'etant le sous-module de torsion
de $M$ et \ \m{r=rg(M)=R(M)}. Puisque $T$ est de type fini il existe un entier
\m{p_0>0} tel que $T$ soit annul\'e par \m{x^{p_0}}. On a donc si \m{p\geq p_0}
\Nligne\m{\dim_C(M\ot_{\ko_{n,P}}\ko_{n,P}/(x^p)) \ = \ \dim_C(T)+rp} ,
et le th\'eor\`eme en d\'ecoule imm\'ediatement.

Montrons maintenant que si \
\m{0\to M'\to M\to M''\to 0} \
est une suite exacte de \m{\ko_{n,P}}-modules de type fini et si le
th\'eor\`eme \ref{R_R_0} est vrai pour deux des modules \m{M'}, $M$, \m{M''},
alors il est vrai pour le troisi\`eme. Pour tout \m{p>0} on a une suite exacte
\[\Tor^1(M'',\ko_{n,P}/(x^p))\lra M'\ot\ko_{n,P}/(x^p)\lra M\ot\ko_{n,P}/(x^p)
\lra M''\ot\ko_{n,P}/(x^p)\lra 0 .\]

Il existe un entier $N$ tel que pour tout
\m{q>0} on ait \ \m{\dim_C(\Tor^1(M'',\ko_{n,P}/(x^q)))\leq N}.
Pour le voir on consid\`ere la r\'esolution de \m{\ko_{n,P}/(x^p)}
\xmat{
0\ar[r] & \ko_{n,P}\ar[r]^-{\times x^p} & \ko_{n,P}\ar[r] & \ko_{n,P}/(x^p) .}
Il en d\'ecoule que \m{\Tor^1(M'',\ko_{n,P}/(x^p)))} est isomorphe au
sous-module \m{M''_p} de \m{M''} des \'el\'ements annul\'es par \m{x^p}. On a
\m{M''_1\subset M''_2\subset\cdots\subset M''_p\subset M''_{p+1}\subset\cdots}
Puisque $M$ est no\'eth\'erien, cette suite stationnaire, d'o\`u le 
r\'esultat.

On a donc
\[0 \ \leq \ \dim_\C(M'\ot\ko_{n,P}/(x^p))+\dim_\C(M''\ot\ko_{n,P}/(x^p))
-\dim_\C(M\ot\ko_{n,P}/(x^p)) \ \leq \ N ,\]
et l'assertion en d\'ecoule imm\'ediatement. On en d\'eduit aussi que si
$M$ a une filtration finie dont tous les gradu\'es poss\`edent la
propri\'et\'e du th\'eor\`eme \ref{R_R_0}, alors il en est de m\^eme de $M$
(d\'emonstration par r\'ecurrence sur la longueur de la filtration).

Pour achever la d\'emonstration du th\'eor\`eme \ref{R_R_0}, on utilise la
filtration canonique de $M$. Ses gradu\'es sont des \m{\ko_{C,P}}-modules, donc
le th\'eor\`eme \ref{R_R_0} est vrai pour eux. Il est donc v\'erifi\'e aussi
pour $M$.
\end{proof}

\sepprop

\begin{subsub}{\bf Corollaire : }\label{str_gen2}
1 - Soit \
\m{0\to M'\to M\to M''\to 0} \
une suite exacte de \m{\ko_{n,P}}-modules de type fini. Alors on a \
\m{R(M)=R(M')+R(M'')} .

2 - Soit \
\m{0\to E\to F\to G\to 0} \
une suite exacte de faisceaux coh\'erents sur \m{C_n}. Alors on a \
\m{R(F)=R(E)+R(G)} , \m{\Deg(F)=\Deg(E)+\Deg(G)} .
\end{subsub}

\begin{proof}
La d\'emonstration de 1- est contenue dans celle du th\'eor\`eme \ref{R_R_0}.
On en d\'eduit l'assertion de 2- sur les rangs g\'en\'eralis\'es. Celle
concernant les degr\'es d\'ecoule du th\'eor\`eme de Riemann-Roch
g\'en\'eralis\'e et du fait que \ \m{\chi(F)=\chi(E)+\chi(G)} .
\end{proof}

\sepprop

\begin{subsub}{\bf Proposition : }\label{str_gen3}
1 - Le rang g\'en\'eralis\'e des \m{\ko_{n,P}}-modules de type fini est
invariant par d\'eformation.

2 - Le rang et le degr\'e g\'en\'eralis\'es des faisceaux coh\'erents sur
\m{C_n} sont invariants par d\'eformation.
\end{subsub}

La seconde assertion signifie que si $\ke$ est une famille plate de faisceaux
coh\'erents sur \m{C_n} param\'etr\'ee par
une vari\'et\'e alg\'ebrique $S$, et si \m{s_0\in S} est un point ferm\'e, il
existe un ouvert \m{U\subset S} contenant \m{s_0} tel que pour tout point
ferm\'e $s$ de $U$ ont ait \m{R(\ke_s)=R(\ke_{s_0})} et \m{\Deg(\ke_s)=
\Deg(\ke_{s_0})}. Pour d\'emontrer 2- on utilise le fait que le polyn\^ome de
Hilbert des faisceaux coh\'erents (cf. \ref{prop_rang}) est invariant par
d\'eformation, et 1- peut se d\'eduire de 2-.

\sepprop

\begin{subsub}{\bf Remarque : }\label{str_rem}\rm
Le rang habituel d'un faisceau coh\'erent $\ke$ sur \m{C_n} est d\'efini de la
fa\c con suivante : pour tout point $P$ de $C$, soit $r_P$ le nombre minimal de
g\'en\'erateurs du \m{\ko_{n,P}}-module \m{\ke_P}. Alors il existe un ouvert
non vide $U$ de $C$ sur lequel $r_P$ prend sa valeur minimale. La valeur de
$r_P$ en les points de $U$ est le {\em rang} de $\ke$. Il est not\'e
\m{rg(\ke)}. Pour ce rang, on va montrer par un exemple que
les propositions \ref{str_gen2} et \ref{str_gen3} sont fausses. Supposons pour
simplifier que $C_n$ soit plong\'ee dans une surface.
soient $p$, $q$ des entiers positifs tels que \m{p+q\leq n}.
Alors on a une suite exacte canonique \
\m{0\to\ko_{p}(-qC)\to\ko_{p+q}\to\ko_{q}\to 0} .
Mais on a \
\m{rg(\ko_{p}(-qC))=rg(\ko_{p+q})=rg(\ko_{q})=1} ,
ce qui contredit la proposition \ref{str_gen2}. D'autre part, l'extension
pr\'ec\'edente montre que \m{\ko_{p+q}} est une d\'eformation de
\m{\ko_{p}(-qC)\oplus\ko_{q}}, ce qui contredit la proposition \ref{str_gen3}.
\end{subsub}

\sepprop

\end{sub}

\sepsec

\section{Faisceaux quasi localement libres}\label{QLL}

On reprend les notations de \ref{QLL0}.

\sepsub

\Ssect{Modules quasi libres et faisceaux coh\'erents quasi localement
libres}{QLL-def2}

Soit $M$ un \m{\ko_{n,P}}-module de type fini.
On dit que $M$ est {\em quasi libre} s'il existe des entiers \m{m_1,\ldots,m_n}
non n\'egatifs et un isomorphisme
\m{M\simeq\bigoplus_{i=1}^n m_i\ko_{i,P}} .

\sepprop

\begin{subsub}{\bf Lemme : }\label{def_lem}
Les entiers \m{m_1,\ldots,m_n} sont uniquement d\'etermin\'es.
\end{subsub}

\begin{proof}
Posons \
\m{M=\bigoplus_{i=1}^nm_i\ko_{i,P}=
\bigoplus_{i=1}^nm_i\ko_{n,P}/(z^i)} .
Alors on a, pour \m{0\leq i< n},
\m{\dim_\C(z^iM/(z^i,x)M)=m_{i+1}+2m_{i+2}+\cdots+(n-i)m_n} ,
ce qui permet de retrouver les \m{m_i} \`a partir de $M$.
\end{proof}

\sepprop

On dit alors que $M$ est {\em de type} \m{(m_1,\ldots,m_n)}. On a \
\m{R(M)=\sigg_{i=1}^ni.m_i}  .

Soit $\ke$ un faisceau coh\'erent sur \m{C_n}.
On dit que $\ke$ est {\em quasi localement libre} en un point $P$ de
$C$ s'il existe un ouvert $U$ de \m{C_n} contenant
$P$ et des entiers non n\'egatifs \m{m_1,\ldots,m_n} tels que pour tout
point $Q$ de $U$, \m{\ke_{n,Q}} soit quasi localement libre de type
\m{m_1,\ldots,m_n}. Les entiers \m{m_1,\ldots,m_n} sont uniquement
d\'etermin\'es d'apr\`es le lemme pr\'ec\'edent, ne d\'ependent que de $\ke$,
et on dit que \m{(m_1,\ldots,m_n)} est le {\em type de } $\ke$.

On dit que $\ke$ est quasi localement libre s'il l'est en tout point de
\m{C_n}.

\sepsubsub

\begin{subsub} Extensions de modules quasi libres (ou de faisceaux quasi
localement libres) - \rm Il est \'evident qu'un somme directe finie de modules
quasi libres (resp. de faisceaux quasi localement libres) est quasi libre
(resp. quasi localement libre). Mais il est faux qu'une extension de modules
quasi libres (ou de faisceaux quasi localement libres) le soit. Par exemple
soit \m{\ki_P} l'id\'eal \m{(x,z)\subset\ko_{n,P}}. On a une suite exacte
\xmat{0\ar[r] & \ko_{n-1,P}\ar[r]^-i & \ki_P\ar[r]^-p & \ko_{C,P}\ar[r] & 0 ,}
o\`u $i$ est d\'efini par \m{i(1)=z}, et $p$ par \m{p(ax+bz)=\ov{a}},
\m{\ov{a}} d\'esignant l'image de $a$ dans \m{\ko_{C,P}}. Cependant \m{(x,z)}
n'est pas quasi libre, car \
\m{\ki_P\ot_{\ko_{n,P}}\ko_{C,P}\simeq \ko_{C,P}\oplus\C_P} \
n'est pas libre.
\end{subsub}

\sepprop

\begin{subsub}{\bf Th\'eor\`eme : }\label{str_gen0}
Soit $M$ un \m{\ko_{n,P}}-module de type fini. Alors $M$ est quasi libre si et
seulement si \m{{\rm Gr}(M)} est un \m{\ko_{C,P}}-module libre.
\end{subsub}

\begin{proof}
Il est clair que si $M$ est quasi libre, alors les termes de \m{{\rm Gr}(M)}
sont libres. R\'eciproquement, supposons que \m{{\rm Gr}(M)} soit libre. Alors
les termes de \m{{\rm Gr}(M)} sont libres, car sinon l'un d'entre eux aurait
un sous-module de torsion non nul et il en serait de m\^eme de \m{{\rm Gr}(M)}.
On d\'emontre que $M$ est quasi libre par r\'ecurrence sur $n$,
le r\'esultat \'etant trivial si \m{n=1}. Supposons donc que \m{n>1} et le
r\'esultat vrai pour \m{n-1}. Alors \m{M_2} est un \m{\ko_{n-1,P}}-module
de type fini dont le gradu\'e est \ \m{\oplus_{i=2}^{n}M_i/M_{i+1}}, et ses
termes sont donc libres. D'apr\`es l'hypoth\`ese de r\'ecurrence, \m{M_2} est
quasi libre, donc il existe des entiers \m{m_1,\ldots,m_{n-1}} tels que \
\m{M_2\simeq\bigoplus_{i=1}^{n-1}m_i\ko_{i,P}} .
D'autre part, il existe un entier $m$ tel que \m{M/M_2=m\ko_{C,P}}.
On a donc une extension
\[0\lra\bigoplus_{i=1}^{n-1}m_i\ko_{i,P}\lra M\lra m\ko_{C,P}\lra 0 .\]
Posons \ \m{M_{2C}=M_2\ot_{\ko_{n,P}}\ko_{C,P}}.

D'apr\`es la construction bien connue des extensions
 et la r\'esolution canonique de \m{\ko_{C,P}}
\[
\xymatrix{\cdots\ko_{n,P}\ar[r]^-{\times z^{n-1}} & \ko_{n,P}\ar[r]^-{\times z}
& \ko_{n,P}\ar[r] & \ko_{C,P}}\]
l'extension pr\'ec\'edente provient d'un morphisme \ \m{\phi:m\ko_{n,P}\to M_2}
\ s'annulant sur \m{z^{n-1}(m\ko_{n,P})}, ce qui est le cas de tous les
morphismes, puisque \m{M_2} est un \m{\ko_{n-1,P}}-module. Plus pr\'ecis\'ement
$M$ est isomorphe au conoyau de \ \m{\phi\oplus(\times z) : m\ko_{n,P}/(z^{n-1})
\to M_2\oplus m\ko_{n,P}}. 

Le morphisme $\phi$
est unique \`a un morphisme  \ \m{m.z\ko_{n,P}\to M_2} \ pr\`es, c'est-\`a-dire
que l'extension est en fait d\'etermin\'ee par l'\'el\'ement de \m{mM_{2C}}
induit par $\phi$.

Montrons que $\phi$ est surjectif. Soit \
\m{\eta=(\phi,z) : m\ko_{n,P}\to M_2\oplus m\ko_{n,P}} .
Alors on a \ \m{M\simeq\coker(\eta)}, la projection \ \m{p:M\to m\ko_{C,P}} \
\'etant induite par \ \m{M_2\oplus m\ko_{n,P}\to m\ko_{n,P}} \ et
l'inclusion \ \m{M_2\subset M} \ par \ \m{M_2\subset M_2\oplus m\ko_{n,P}}.
Puisque
\ \m{m\ko_{C,P}=M_{\mid C}}, on a \ \m{\ker(p)=zM}. Si \m{a\in M_2},
\m{b\in m\ko_{n,P}}, on note \m{\ov{(a,b)}} l'image de \m{(a,b)} dans $M$.
Alors on a \ \m{p(\ov{(a,b)})=0} si et seulement si $b$ est de la forme
\m{b=z\beta}. La condition \m{M_2=zM} implique
que pour tous \m{a\in M_2} et \m{\beta\in M_2} il existe \m{\gamma\in M_2}
et \m{\alpha\in m\ko_{n,P}} tels que \ \m{\ov{(a,z\beta)}=
\ov{(z\gamma,z\alpha)}}, et donc il existe \m{v\in m\ko_{n,P}} tel que \
\m{(a,z\beta)=(z\gamma+\phi(v),z(\alpha+v))} .
Il en d\'ecoule que \ \m{\imm(\phi)+zM_2=M_2}, donc
\[\imm(\phi)+z(\imm(\phi)+zM_2) \ = \ \imm(\phi)+z^2M_2 \ = \ M_2 .\]
On obtient de m\^eme \ \m{\imm(\phi)+z^rM_2=M_2} \ pour tout \m{r>0}, et en
prenant \m{r=n-1}, on voit que \ \m{\imm(\phi)=M_2}, c'est-\`a-dire que $\phi$
est surjectif.

Puisque l'extension est d\'etermin\'ee par l'\'el\'ement de \m{mM_{2C}} induit
par $\phi$, on peut supposer que $\phi$ provient d'applications lin\'eaires \
\m{\phi_i:\C^m\to\C^{m_i}} \ dont la somme \
\m{\C^m\to\oplus_{i=1}^{n-1}\C^{m_i}} est surjective. On se ram\`ene donc au
cas o\`u on a une d\'ecomposition \
\m{\C^m=\C^k\oplus\biggl(\oplus_{i=1}^{n-1}\C^{m_i}\biggr)} ,
et o\`u \m{\phi_i} est la projection. Il en d\'ecoule que $M$ est isomorphe \`a
la somme directe de \m{k\ko_{C,P}} et des extensions
\[0\lra m_i\ko_{i,P}\lra N_i\lra m_i\ko_{C,P}\lra 0\]
d\'efinies par l'identit\'e \m{\C^{m_i}\to\C^{m_i}}. Il est ais\'e de voir que
\m{N_i=m_i\ko_{i+1,P}}. On obtient finalement \
\m{M\simeq k\ko_{C,P}\oplus\biggl(\oplus_{i=1}^{n-1}m_i\ko_{i+1,P}\biggr)} ,
 ce qui d\'emontre le th\'eor\`eme.
\end{proof}

\sepprop

\begin{subsub}{\bf Corollaire : }\label{str_gen0b}
Soit $\ke$ un faisceau coh\'erent sur \m{C_n}. Alors $\ke$ est quasi localement
libre si et seulement si tous les termes de \m{{\rm Gr}(\ke)} sont localement
libres sur $C$.
\end{subsub}

\sepprop

\begin{subsub}{\bf Remarques : }\label{str_gen0c}\rm 1 -
Soit $\ke$ un faisceau coh\'erent sur \m{C_n}. Alors \m{{\rm Gr}_2(\ke)}, le
gradu\'e de la seconde filtration canonique de $\ke$, est localement libre si
$\ke$ est quasi localement libre. Mais la r\'eciproque est fausse : le faisceau
d'id\'eaux d'un point de $C$ n'est pas quasi localement libre, mais le gradu\'e
de sa seconde filtration canonique est localement libre.

2 - Il d\'ecoule ais\'ement de la d\'emonstration du th\'eor\`eme
\ref{str_gen0} que si $\ke$ est de type \m{(m_1,\ldots,m_n)} et quasi libre en
$P$, alors il existe un ouvert $U$ de \m{C_n} contenant $P$ tel que \
\m{\ke_{\mid U}\simeq\bigoplus_{i=1}^nm_i\ko_{i\mid U}}~.
\end{subsub}

\sepprop

Du th\'eor\`eme \ref{str_gen0} d\'ecoule imm\'ediatement le suivant, qui donne
la structure g\'en\'erique des faisceaux coh\'erents sur \m{C_n} :

\sepprop

\begin{subsub}{\bf Th\'eor\`eme : }\label{str_gen}
Soit $\ke$ un faisceau coh\'erent sur \m{C_n}. Alors il existe un ouvert non
vide $U$ tel que $\ke$ soit quasi localement libre en tout point de $U$.
\end{subsub}

\end{sub}

\sepsub

\Ssect{Morphismes de modules quasi libres et morphismes de faisceaux
coh\'erents quasi localement libres}{QLL-mor}

\begin{subsub}{\bf Th\'eor\`eme : }\label{str_mor}
1 - Soient $M$, $F$ des \m{\ko_{n,P}}-modules quasi libres et \m{f:M\to F} un
morphisme surjectif. Alors \m{\ker(f)} est quasi libre.

2 - Soient $\ke$, $\kf$ des faisceaux coh\'erents quasi localement libres sur
\m{C_n} et \m{\phi:\ke\to\kf} un morphisme surjectif. Alors  \m{\ker(\phi)} est
quasi localement libre.
\end{subsub}

\begin{proof}
On ne d\'emontrera que 1-, 2- s'en d\'eduisant ais\'ement.

{\em \'Etape 1 -}
On suppose d'abord que $F$ est un \m{\ko_{C,P}}-module libre.
On va se ramener au cas o\`u \ \m{z^{n-1}M=\nsp}.
Soit \m{N=\ker(f)}. Le
morphisme $f$ se factorise par \ \m{M_C=M\ot_{\ko_{n,P}}\ko_{C,P}} \ (vu comme
\m{\ko_{n,P}}-module) :
\xmat{f:M\flon[r] & M_C\flon[r]^-{\ov{f}} & F}
On peut donc supposer que $F$ est un quotient libre de \m{M_C}. Posons
\[M \ = \ \bigoplus_{i=1}^n M_i\ot_\C\ko_{i,P} , \ \ \ \
F \ = \ F_1\ot_\C\ko_{C,P} ,\]
o\`u \m{M_1,\ldots,M_n,F_1} sont des $\C$-espaces vectoriels de dimension
finie. Les \'el\'ements de $M$ sont donc des sommes \m{w_1+\cdots+w_n}, avec \
\m{w_i\in M_i\ot_\C\ko_{i,P}}. On note \m{\ov{w_i}} l'image de \m{w_i} dans
\m{M_i\ot_\C\ko_{C,P}}. Le morphisme \m{\ov{f}} provient d'une application
lin\'eaire \
\m{\ov{f} : \bigoplus_{i=1}^nM_i\to F_1} .
Soit \m{\ov{N}} son noyau. On a
\[N \ = \ \bigg\lbrace w_1+\cdots+w_n\in M \ ; \
\ov{w_1}+\cdots+\ov{w_n}\in\ov{N}\bigg\rbrace .\]
Soit \m{N_n\subset M_n} l'image de la projection \m{\ov{N}\to M_n}. D'apr\`es
la description pr\'ec\'edente on a
\[z^{n-1}N \ = \ N_n\ot_\C(z^{n-1}\ko_{n,P}) \ \subset \
M_n\ot_\C(z^{n-1}\ko_{n,P})=z^{n-1}M .\]
Soit \m{N'_n\subset M_n} un sous-espace vectoriel suppl\'ementaire de \m{N_n}.
On a
\[M/z^{n-1}N \ = \ (M_1\ot_\C\ko_{1,P})\oplus\cdots\oplus
(M_{n-2}\ot_\C\ko_{n-2,P})\oplus\big((M_{n-1}\oplus N_n)\ot_\C\ko_{n-1,P})\big)
\oplus(N'_n\ot_\C\ko_{n,P}) ,\]
donc \m{M/z^{n-1}N} est quasi libre. Supposons que l'on puisse prouver le
r\'esultat si \m{z^{n-1}N=\nsp}. Si \m{z^{n-1}N\not=\nsp}, on en d\'eduit que
\m{N/z^{n-1}N} est quasi libre : on voit ais\'ement que c'est le noyau d'un
morphisme surjectif de \m{M/z^{n-1}N} dans un \m{\ko_{C_P}}-module libre.
Donc d'apr\`es le th\'eor\`eme \ref{str_gen0} \m{Gr(N/z^{n-1}N)} est libre. On
a \ \m{Gr(N)=Gr(N/z^{n-1}N)\oplus z^{n-1}N}, donc \m{Gr(N)} est libre, et
d'apr\`es le th\'eor\`eme \ref{str_gen0} $N$ est quasi libre.

\medskip

{\em \'Etape 2 -}
On suppose toujours que $F$ est un \m{\ko_{C,P}}-module libre.
D'apr\`es l'\'etape 1 il suffit de traiter le cas o\`u \ \m{z^{n-1}N=\nsp}.
On d\'emontre que $N$ est quasi libre par r\'ecurrence sur $n$. C'est \'evident
si \m{n=1}. Supposons que \m{n>1} et que ce soit vrai pour \m{n-1}. On a
\ \m{\ov{N}\subset M_1\oplus\cdots M_{n-1}}, et on peut supposer que \
\m{F=\big((M_1\oplus\cdots M_{n-1})/\ov{N}\big)\oplus M_n} .
Soient
\[f':(M_1\ot_\C\ko_{1,P})\oplus\cdots\oplus(M_{n-1}\ot_\C\ko_{n-1,P})\lra
F'=\big((M_1\oplus\cdots M_{n-1})/\ov{N}\big)\ot_\C\ko_{C,P} \]
et \m{N'=\ker(f')}. Alors \m{f'} est surjectif et \ \m{N=N'\oplus(M_n\ot_\C
\ko_{n-1,P})}. D'apr\`es l'hypoth\`ese de r\'ecurrence \m{N'} est quasi libre,
donc il en est de m\^eme pour $N$. Le th\'eor\`eme \ref{str_mor} est donc
d\'emontr\'e dans les cas o\`u $F$ est un \m{\ko_{C,P}}-module.

\medskip

{\em \'Etape 3 -}
On traite maintenant le cas g\'en\'eral, par r\'ecurrence sur le plus petit
entier $m$ tel que $F$ soit un \m{\ko_{m,P}}-module, c'est-\`a-dire tel que
\m{z^mF=\nsp}. Le cas \m{m=1} a d\'ej\`a \'et\'e trait\'e.
Supposons que \m{m>1} et que ce soit vrai pour \m{m-1}.
Soient \m{N=\ker(f)}, $X$ le noyau du morphisme canonique \m{F\to
F_C=F\ot_{\ko_{n,P}}\ko_{C,P}} et \m{N_0} le noyau du morphisme compos\'e
surjectif
\xmat{M\flon[r]^f & F\flon[r] & F_C .}
On a un diagramme commutatif avec lignes et colonnes exactes
\[\xymatrix{
        &                 &                   & 0\ar[d]\\
        & 0\ar[d]         &                   & X\ar[d]\\
0\ar[r] & N\ar[r]\ar[d]   & M\fleq[d]\ar[r]^f & F\ar[r]\ar[d]   & 0\\
0\ar[r] & N_0\ar[r]\ar[d] & M\ar[r]           & F_C\ar[r]\ar[d] & 0\\
        & X\ar[d]         &                   & 0\\
	& 0
}\]
D'apr\`es l'\'etape 2, \m{N_0} est quasi libre, et d'apr\`es l'hypoth\`ese de
r\'ecurrence et la suite exacte verticale de gauche, $N$ est quasi libre. Ceci
ach\`eve la d\'emonstration du th\'eor\`eme \ref{str_mor}.
\end{proof}
\end{sub}

\sepsub

\Ssect{D\'eformations des $\ko_{n,P}$-modules quasi libres}
{QLL-def2b}

\begin{subsub}{\bf Conjecture : }\label{str_theo}
Soient $M$, $N$ des \m{\ko_{n,P}}-modules quasi libres. Alors
$N$ est une d\'eformation de $M$ si et seulement si on a \ \m{R(M)=R(N)} \ et
\ \m{F(M)\leq F(N)}.
\end{subsub}

Ce r\'esultat est d\'emontr\'e dans le cas o\`u \m{n=2} (cf. \ref{def_f_db2}).

\end{sub}

\sepsec

\section{Faisceaux coh\'erents sur les courbes doubles}\label{faisc_doub}

Soient $C_2$ une courbe projective double primitive et \
\m{C=C_1\subset C_2} \
sa filtration canonique. Soit $L$ le fibr\'e en droites associ\'e sur $C$, qui
n'est autre que \m{\ki_C}, le faisceau d'id\'eaux de $C$ dans \m{C_2}. On
notera \m{\ko_2=\ko_{C_2}}.

Si \m{P\in C}, on notera \m{z_P} (ou $z$ s'il n'y a pas d'ambiguit\'e) une
\'equation locale de $C$ dans \m{\ko_{2,P}}, et \m{x_P} (ou $x$ s'il n'y a pas
d'ambiguit\'e) un \'el\'ement de \m{\ko_{2,P}} au dessus d'un g\'en\'erateur
de l'id\'eal maximal de \m{\ko_{C,P}}.

\sepsub

\Ssect{Structure des faisceaux quasi-localement libres sur \m{C_2}}{QLL_C2}

\begin{subsub}\label{QLL_C2_0}\rm{\bf Premi\`ere filtration canonique - }
Soient $\kf$ un faisceau coh\'erent quasi localement libre sur \m{C_2},
\m{0\subset E \subset\kf} sa filtration canonique, et \m{F=\kf/E}. On a donc
une suite exacte
\[0\lra E\lra\kf\lra F\lra 0\]
et $E$, $F$ sont des fibr\'es vectoriels sur $C$. Soit \
\m{\Phi_\kf:\kf\ot\ki_C\lra\kf} \
le morphisme canonique. On a
\m{\kf\ot\ki_C\simeq F\ot L} , \m{\imm(\Phi_\kf)=E} ,
donc on peut voir $\Phi_\kf$ comme un morphisme surjectif \ \m{F\ot L\to E}.
On a une suite exacte canonique
\[0\lra H^1(\HHom(F,E))\lra\Ext^1_{\ko_2}(F,E)\lra H^0(\EExt^1_{\ko_2}(F,E))
\lra 0\]
d\'ecoulant de la suite spectrale des Ext. 
D'apr\`es le th\'eor\`eme \ref{pr5} Il existe un faisceau
localement libre $\F$ (resp. $\L$) sur \m{C_2} extension de $F$ (resp. $L$).
On a une r\'esolution localement libre de $F$
\xmat{\cdots\F\ot\L^2\ar[r]^-{f_2} & \F\ot\L\ar[r]^-{f_1} & \F\ar[r]^-{f_0}
& F\lra 0 ,}
qu'on construit en utilisant les suites exactes
\[0\lra F\ot L^{i+1}\lra\F\ot\L^i\lra F\ot L^i\lra 0 .\]
On en d\'eduit un isomorphisme canonique \
\m{\EExt^1_{\ko_2}(F,E) \ \simeq \ \HHom(F\ot L,E).} .
On a donc une suite exacte
\xmat{0\ar[r] & H^1(F^*\ot E)\ar[r] & \Ext^1_{\ko_2}(F,E)\ar[r]^-\delta &
\Hom(F\ot L,E)\ar[r] & 0 .}

Soit \m{\sigma\in\Ext^1_{\ko_2}(F,E)} l'\'el\'ement associ\'e \`a la suite
exacte \
\m{0\to E\to\kf\to F\to 0} .
Alors on a \ \m{\delta(\sigma)=\Phi_\kf} . Cela se voit ais\'ement en utilisant
la r\'esolution localement libre pr\'ec\'edente de $F$ ainsi que la
construction des extensions de
faisceaux d\'ej\`a utilis\'ee dans le d\'emonstration du th\'eor\`eme 
\ref{str_gen0}.
\end{subsub}

\sepprop

\begin{subsub}\label{QLL_C2_1}\rm{\bf Seconde filtration canonique - }
Soient $\Gamma$ le noyau de \
\m{\Phi_\kf\ot I_{L^*}:F\to E\ot L^*} \ et $G$ celui du morphisme
compos\'e
\xmat{\kf\ar[r] & F\ar[rr]^-{\Phi_\kf\ot I_{L^*}} & & E\ot L^* .}
On a donc des suites exactes
\[0\lra G\lra\kf\lra E\ot L^*\lra 0 , \ \ \ \
0\lra E\lra G\lra\Gamma\lra 0 .\]
Alors le faisceau $G$ est localement libre sur $C$. C'est le plus grand
sous-faisceau de $\kf$ de support $C$. Autrement dit, \m{0\subset G\subset\kf}
est la seconde filtration canonique de $\kf$.
\end{subsub}

\sepprop

\begin{subsub}\label{QLL_C2_2}\rm
{\bf R\'eciproque - Construction de faisceaux quasi localement libres - }
Soit \ \m{\Phi:F\ot L\to E} \ un morphisme surjectif.
Soient \m{\sigma\in\Ext^1_{\ko_2}(F,E)} au dessus de $\Phi$ et
\[0\lra E\lra\kf\lra F\lra 0\]
l'extension associ\'ee \`a $\sigma$. Alors il est facile de voir que
la filtration canonique de $\kf$ est \m{0\subset E\subset\kf}.
\end{subsub}

\sepprop

\begin{subsub}\label{QLL_C2_3}\rm{\bf Notations : }
\m{E_\kf=E} , \m{F_\kf=F}  , \m{G_\kf=G} , \m{\Gamma_\kf=\Gamma} .
\end{subsub}

\sepprop

\begin{subsub}\label{QLL_C2_4}\rm
{\bf Produits tensoriels - } Soient $\ke$, $\kf$ des faisceaux quasi localement
libres sur \m{C_2},
alors on a \ \m{E_{\ke\ot\kf}=E_\ke\ot E_\kf\ot L^*}, \m{F_{\ke\ot\kf}=
F_\ke\ot F_\kf}, \m{G_{\ke\ot\kf}=G_\ke\ot G_\kf\ot L^*} .
\end{subsub}
\end{sub}

\sepsub

\Ssect{ R\'esolutions localement libres de faisceaux quasi localement
libres}{QLL_C2B}

Soit $\ke$ un faisceau quasi localement libre sur \m{C_2} . Alors on a un
diagramme commutatif avec lignes exactes
\xmat{0\ar[r] & \Ext^1_{\ko_C}(F_\ke,F_\ke\ot L)\ar[r]^-a\flon[d]^p &
\Ext^1_{\ko_2}(F_\ke,F_\ke\ot L)\ar[r]^-b\ar[d]^\pi &
\End(F_\ke\ot L)\ar[r]\ar[d]^q & 0 \\
0\ar[r] & \Ext^1_{\ko_C}(F_\ke,E_\ke)\ar[r]^-\alpha &
\Ext^1_{\ko_2}(F_\ke,E_\ke)\ar[r]^-\beta & \Hom(F_\ke\ot L,E_\ke)\ar[r] & 0
}
les fl\`eches verticales \'etant induites par \m{\Phi_\ke} . 

\sepprop

\begin{subsub}\label{QLL_C2_4_0}{\bf Lemme : }
Soit \ \m{\sigma\in\Ext^1_{\ko_2}(F_\ke,E_\ke)} \ correspondant \`a \m{\ke}.
Alors on a \ \m{\sigma\in\imm(\pi)} . 
\end{subsub}

\begin{proof}
On a \ \m{\beta(\sigma)=\Phi_\ke=q(I_{F_\ke\ot L})} . Soit \ \m{\theta\in
\Ext^1_{\ko_2}(F_\ke,F_\ke\ot L)} tel que \Nligne
\m{b(\theta)=I_{F_ke\ot L}} . On a
\m{\beta\circ\pi(\theta)=\Phi_\ke} , donc \ \m{\sigma-\pi(\theta)\in
\imm(\alpha)} . Puisque $p$ est surjective, il existe \m{\lambda\in
\Ext^1_{\ko_C}(F_\ke,F_\ke\ot L)} \ tel que \ \m{\sigma-\pi(\theta)=\alpha\circ
p(\lambda)} . On a donc \ \m{\sigma=\pi(\theta+\alpha(\lambda))} .
\end{proof}

\sepprop

Il existe donc un faisceau localement libre $\kf$ et un diagramme commutatif
\xmat{
0\ar[r] & F_\ke\ot L\ar[r]\ar[d]^{\Phi_\ke} & \kf\ar[r]\ar[d]^f &
F_\ke\ar[r]\fleq[d] & 0 \\
0\ar[r] & E_\ke\ar[r] & \ke\ar[r] & F_\ke\ar[r] & 0
}
o\`u les suites exactes horizontales proviennent des premi\`eres filtrations
canoniques. Le morphisme $f$ est surjectif et on a \ \m{\ker(f)=\Gamma_\ke\ot
L} . On a d'autre part comme dans \ref{QLL_C2_0} une r\'esolution localement
libre de \m{\Gamma_\ke\ot L} :
\xmat{\cdots \ar[r] & {\bf \Gamma}\ot\L^2\ar[r]^-{f_2} & {\bf \Gamma}\ot\L
\ar[r]^-{f_1} & \Gamma_\ke\ot L\ar[r] & 0}
o\`u $\bf \Gamma$ est un fibr\'e vectoriel sur \m{C_2} extension de
\m{\Gamma_\ke} . On en d\'eduit une r\'esolution localement libre de $\ke$
\xmat{\cdots \ar[r] & {\bf \Gamma}\ot\L^2\ar[r]^-{f_2} & {\bf \Gamma}\ot\L
\ar[r]^-{f_1} & \kf\ar[r]^-f & \ke\ar[r] & 0}

\sepprop

\begin{subsub}\label{QLL_C2_4_1}{\bf Corollaire : }
On a, pour tout \m{i\geq 1} et tout faisceau quasi localement libre $\kg$ sur
\m{C_2} un isomorphisme canonique \
\m{\EExt^i_{\ko_2}(\ke,\kg)\simeq\HHom(\Gamma_\ke\ot L^i,\Gamma_\kg)} .
\end{subsub}

\end{sub}

\sepsub

\Ssect{Structure des faisceaux sans torsion sur \m{C_2}}{ST_C2}

\begin{subsub}\label{C2_x2b2}Faisceaux de torsion sur $C$ - \rm Soient
\m{P\in C} et $N$ un \m{\ko_{C,P}}-module de torsion. Alors il existe des
entiers \m{n_1,\ldots,n_p} uniques, tels que \m{n_1\geq\cdots\geq n_p} et que
\[N \ \simeq \ \bigoplus_{i=1}^p\ko_{C,P}/(x_P^{n_i}) .\]
Soit $T$ un faisceau coh\'erent de torsion sur $C$. On pose
\m{\wT=\EExt^1_{\ko_C}(T,\ko_C)} .
Alors on a \m{\wT\simeq T} (isomorphisme non canonique), mais on a un
isomorphisme canonique \m{{\widetilde\wT}\simeq T}. Soit $E$ est un fibr\'e
vectoriel sur $C$. Alors on a \
\m{\Ext^i_{\ko_C}(T,E)=\nsp \quad {\rm si \ } i\not=1} \
et \ \m{\Ext^1_{\ko_C}(T,E)\simeq\Hom(E^*,\wT)}.
D'apr\`es la proposition \ref{pr3} et le fait que
\m{\Tor^1_{\ko_2}(E,\ko_C)} est de torsion, on a un isomorphisme canonique \
\m{\Ext^1_{\ko_C}(T,E)\simeq\Ext^1_{\ko_2}(T,E)} .
Si $T$, \m{T'} sont des faisceaux coh\'erents de torsion sur $C$ et
\m{f:T\to T'} un morphisme, on note \m{\widetilde f} le morphisme induit
\m{{\widetilde {T'}}\to\wT}.
\end{subsub}

\sepprop

\begin{subsub}{\bf Lemme : }\label{C2_x3_0}
Soient $\F$ un fibr\'e vectoriel sur \m{C_2} et \m{F=\F_{\mid C}}. Alors

1 - Le morphisme canonique \
\m{\Ext^1_{\ko_2}(T,\F)\to\Ext^1_{\ko_2}(T,F)} \
est nul. On a donc un isomorphisme canonique \
\m{\Ext^1_{\ko_2}(T,\F)\simeq\Ext^1_{\ko_2}(T,F\ot L)} .

2 - Pour tout \m{j\geq 2} on a \ \m{\Ext^j_{\ko_2}(T,\F)=\nsp} , et si
\m{j\geq 1} on a
\[\Ext^j_{\ko_2}(T,F) \ \simeq \ \Ext^1_{\ko_2}(T,F\ot (L^*)^{j-1}) \ \simeq \
\Hom(E^*\ot L^{j-1},\wT) \ .\]
\end{subsub}

\begin{proof}
La conclusion de 1- d\'ecoule de la suite exacte \
\m{0\to F\ot L\to\F\to F\to 0} .
Puisque \m{\HHom(T,\F)=\HHom(T,F)=0}, on a des isomorphismes canoniques
\[\Ext^1_{\ko_2}(T,\F)\simeq H^0(\EExt^1_{\ko_2}(T,\F)) , \ \ \ \
\Ext^1_{\ko_2}(T,F)\simeq H^0(\EExt^1_{\ko_2}(T,F)) ,\]
et il suffit, pour prouver 1-,  de montrer que le morphisme canonique \
\m{\EExt^1_{\ko_2}(T,\F)\to\EExt^1_{\ko_2}(T,F)} \
est nul. On se ram\`ene donc au probl\`eme suivant : soient \m{i>0} un entier
et \m{P\in C}, il faut montrer que le morphisme canonique
\[\Phi: \Ext^1_{\ko_{2,P}}(\ko_{C,P}/(x_P^i),\ko_{2,P})\lra
\Ext^1_{\ko_{2,P}}(\ko_{C,P}/(x_P^i),\ko_{C,P})\]
est nul. On utilise la r\'esolution libre suivante de \m{\ko_{C,P}/(x_P^i)} (en
tant que \m{\ko_{2,P}}-module) :
\xmat{
\cdots 2\ko_{2,P}\ar[r]^-f & 2\ko_{2,P}\ar[r]^-h & \ko_{2,P}\flon[r] &
\ko_{C,P}/(x_P^i)}
avec \m{f=\begin{pmatrix}z_P & -x_P^i\\ 0 & z_P\end{pmatrix}} ,
\m{h=(z_P,x_P^i)} . On en d\'eduit le diagramme commutatif
\xmat{
\ko_{2,P}\ar[rr]^-{^th}\ar[d] & & 2\ko_{2,P}\ar[rr]^-{^tf}\ar[d]^\phi & &
2\ko_{2,P}\ar[d]\\
\ko_{C,P}\ar[rr]^-{^th\ot I_{\ko_{C,P}}} & &
2\ko_{C,P}\ar[rr]^-{^tf\ot I_{\ko_{C,P}}} & &
2\ko_{C,P}
}
On a
\[\Ext^1_{\ko_{2,P}}(\ko_{C,P}/(x_P^i),\ko_{2,P})\simeq\ker(^tf)/\imm(^th) ,\]
\[\Ext^1_{\ko_{2,P}}(\ko_{C,P}/(x_P^i),\ko_{C,P})\simeq
\ker(^tf\ot I_{\ko_{C,P}})/\imm(^th\ot I_{\ko_{C,P}}) ,\]
et $\Phi$ est induit par $\phi$.
Un calcul simple montre que
\[\ker(^tf)/\imm(^th) \ \simeq \ ((z_P,x_P^i),(0,z_P))/((z_P,x_P^i)) ,\]
est donc engendr\'e par la classe de \m{(0,z_P)}, dont l'image dans
\m{2\ko_{C,P}} est nulle. On a donc \m{\Phi=0}. Ceci d\'emontre 1-.

Pour d\'emontrer que \ \m{\Ext^j_{\ko_2}(T,\F)=\nsp} \ on utilise encore la
r\'esolution pr\'ec\'edente de \m{\ko_{C,P}/(x_P^i)} et les autres assertions
de 2- d\'ecoulent alors de la suite exacte \ \m{0\to F\ot L\to\F\to F\to 0} .
\end{proof}

\sepsubsub

\begin{subsub}\label{C2_x3_2}Faisceaux sans torsion sur \m{C_2} - \rm
Soient $\ke$ un faisceau coh\'erent sans torsion sur \m{C_2}, \m{0\subset E
\subset\ke} sa filtration canonique, \m{F\oplus T=\ke/E}, $F$ \'etant
localement libre et $T$ de torsion sur $C$. On a donc une suite exacte
\[0\lra E\lra\ke\lra F\oplus T\lra 0 .\]
Soit $\kf$ le noyau du morphisme \m{\ke\to T}. On a un diagramme commutatif
avec lignes et colonnes exactes
\[(D1) \quad\quad\quad \xymatrix{
        &                 & 0\ar[d]              & 0\ar[d] \\
0\ar[r] & E\ar[r]\fleq[d] & \kf\ar[r]\ar[d]      & F\ar[r]\ar[d]         & 0 \\
0\ar[r] & E\ar[r]         & \ke\ar[r]^-q\ar[d]^p & F\oplus T\ar[r]\ar[d] & 0 \\
        &                 & T\fleq[r]\ar[d]      & T\ar[d] \\
        &                 & 0                    & 0
}\]
\end{subsub}
Notons que $\kf$ n'est pas uniquement d\'etermin\'e, mais $E$, $F$ et $T$ le
sont. On peut donc noter \m{E=E_\ke}, \m{F=F_\ke}, \m{T=T_\ke}.

\sepprop

\begin{subsub}{\bf Lemme : }\label{C2_x3}
Le faisceau $\kf$ est quasi localement libre et sa filtration canonique est
\m{0\subset E\subset\kf}.
\end{subsub}

\begin{proof}
Il suffit de d\'emontrer la seconde assertion, qui entra\^ine le premi\`ere
d'apr\`es le th\'eor\`eme \ref{str_gen0}. Soit \m{P\in C}. Il suffit de montrer
que \m{E_P\subset z_P\kf_P}. Soit \m{e\in E_P}. Alors il existe
\m{\epsilon\in\ke_P} tel que \m{e=z\epsilon} (car \m{0\subset E\subset\ke} est
la filtration canonique de $\ke$). Soit \ \m{q(\epsilon)=(f,t)}. Il existe
\m{u\in\ke_P} tel que \m{q(u)=(0,t)}, et un entier \m{k>0} tel que
\m{x_P^kt=0}. On a donc \m{q(x_P^ku)=0}, c'est-\`a-dire que \m{x_P^ku\in E_P}.
On a donc \m{z_Px_P^ku=0}, d'o\`u \m{zu=0}, car $\ke$ est sans torsion. Soit
\m{\epsilon'=\epsilon-u}. Alors on a \m{z\epsilon'=e}, et \m{p(\epsilon')=0},
c'est-\`a-dire \m{\epsilon'\in\kf_P}. On a donc \m{e\in z\kf_P}.
\end{proof}

\sepprop

Le morphisme canonique \m{\Phi_\kf:F\ot L\to E} ne d\'epend pas de $\kf$. En
fait il est induit par \m{\ke\ot\ki_C\to\ke}. On peut donc le noter
\m{\Phi_\ke}. On notera \m{\Gamma(\ke)=\ker(\Phi_\ke)}.

On note \m{V_\ke} le noyau du morphisme surjectif \m{\ke\to F}. 
Il est concentr\'e sur $C$ et localement libre sur $C$ (pour le voir on remarque
qu'il est contenu dans le premier terme de la seconde filtration canonique de
$\ke$ d\'ecrite plus loin).

On note \m{\sigma_\ke} l'\'el\'ement de \m{\Ext^1_{\ko_2}(T,\kf)} associ\'e \`a
la suite exacte \
\m{0\to\kf\to\ke\to T\to 0} .
D'apr\`es le diagramme commutatif pr\'ec\'edent l'image de \m{\sigma_\ke} dans
\m{\Ext^1_{\ko_2}(T,F)} est nulle. On peut donc voir \m{\sigma_\ke} comme un
\'el\'ement de \m{\Ext^1_{\ko_2}(T,E)}, qui est ind\'ependant du choix de $\kf$
et d\'efini \`a un automorphisme de $T$ pr\`es. Rappelons qu'on a un
isomorphisme canonique
\[\Ext^1_{\ko_2}(T,E) \ \simeq \ \Hom(E^*,\wT) .\]

\sepprop

\begin{subsub}{\bf Lemme : }\label{C2_x4}
Le morphisme \ \m{\sigma_\ke:E^*\to\wT} \ est surjectif.
\end{subsub}

\begin{proof}
On d\'eduit de la suite exacte \ \m{0\to\kf\to\ke\to T\to 0} \ une injection
\Nligne \m{\End(T)\lra\Ext^1_{\ko_2}(T,\kf)} .
On voit ais\'ement qu'elle associe \`a un automorphisme $\tau$ de $T$ le
morphisme \ \m{\widetilde{\tau}\circ\sigma_\ke:E^*\to\wT}. Si \m{\sigma_\ke}
n'est pas surjectif, il d\'ecoule ais\'ement de la description des faisceaux de
torsion sur $C$ (cf. \ref{C2_x2b2}) qu'il existe \m{\tau\in\End(T)} tel que
\m{\widetilde{\tau}\circ\sigma_\ke=0}.
\end{proof}

\sepsubsub

\begin{subsub}\label{C2_x5_01}{\bf Remarques : }\rm 1 - On montre ais\'ement que
\m{\sigma_\ke} est l'\'el\'ement de \m{\Ext^1_{\ko_C}(T,E)} associ\'e \`a la
suite exacte \
\m{0\to E\to V_\ke\to T\to 0} .
Ceci permet de red\'emontrer le lemme \ref{C2_x4}.

2 - {\em R\'eciproque - } Soient $\kf$ un faisceau coh\'erent quasi localement
libre sur \m{C_2}, \m{0\subset E\subset\kf} sa filtration canonique,
\m{F=\kf/E}, et \ \m{\sigma : E^*\to\wT} \ un morphisme surjectif. On peut voir
$\sigma$ comme un \'el\'ement de \m{\Ext^1_{\ko_2}(T,E)}. On en d\'eduit un
\'el\'ement \m{\sigma_0} de \m{\Ext^1_{\ko_2}(\kf,T)} et l'extension associ\'ee
\ \m{0\to\kf\to\ke\to T\to 0} .
Puisque l'image de \m{\sigma_0} dans \m{\Ext^1_{\ko_2}(F,T)} est nulle, on
obtient le diagramme commutatif de \ref{C2_x3_2}. Soit \m{P\in C}. Puisque \
\m{E_P\subset z_P\kf_P\subset z_P\ke_P}, \m{0\subset E\subset\ke} \ est la
filtration canonique de $\ke$. On montre ais\'ement comme dans le lemme
\ref{C2_x4} que la surjectivit\'e de $\sigma$ entra\^ine que $\ke$ est sans
torsion. Enfin, on a \m{\sigma_\ke=\sigma}.
\end{subsub}

\sepsubsub

\begin{subsub}\label{C2_x5_0}Index d'un faisceau sans torsion - \rm
On pose \ \m{i(\ke)=h^0(T)} , qu'on appelle l'{\em index } de $\ke$.
\end{subsub}

\sepsubsub

\begin{subsub}\label{C2_x5_1}Seconde filtration canonique - \rm Soit \m{G_\ke}
le noyau du morphisme compos\'e
\xmat{\ke\ar[r] & F\ar[rr]^-{\Phi_\ke\ot I_{L^*}} & & E\ot L^* .}
On a donc un diagramme commutatif avec lignes et colonnes exactes
\[(D) \quad\quad\quad \xymatrix{
        & 0\ar[d]                & 0\ar[d]\\
0\ar[r] & G_\kf\ar[r]\ar[d]     & G_\ke\ar[r]\ar[d] & T\ar[r]\fleq[d] & 0\\
0\ar[r] & \kf\ar[r]\ar[d]        & \ke\ar[r]\ar[d]    & T\ar[r]         & 0\\
        & E\ot L^*\fleq[r]\ar[d] & E\ot L^*\ar[d]\\
	& 0                      & 0
}\]
\end{subsub}

Il est ais\'e de voir que
le faisceau \m{G_\ke} est concentr\'e sur $C$ et localement libre sur $C$.
De plus, \m{0\subset G_\ke\subset\ke} est la seconde filtration canonique de
$\ke$.

\sepprop

\begin{subsub}\label{C2_x5}Familles de faisceaux sans torsion - \rm
Soit $S$ une vari\'et\'e alg\'ebrique. On appelle {\em famille de faisceaux
sans torsion sur \m{C_2} param\'etr\'ee par $S$} un faisceau coh\'erent $\kk$
sur \m{S\times C_2}, plat sur $S$ et tel que pour tout point ferm\'e
\m{s\in S}, \m{\kk_s} soit sans torsion sur \m{C_2}.
\end{subsub}

\sepprop

\begin{subsub}\label{C2_x5b}Exemples - \rm
Soient \m{P\in C} \m{k\geq 1} un entier et \m{\beta\in\ko_{C,P}}. On pose
\m{x=x_P}, \m{z=z_P}. On note
\m{\kl_{k,P,x,\beta}} le faisceau d'id\'eaux sur \m{C_2} tel que pour tout
\m{Q\in C\backslash\lbrace P\rbrace}, on ait \
\m{(\kl_{k,P,x,\beta})_Q=\ko_{2,P}}~,
et \ \m{(\kl_{k,P,x,\beta})_P = (x^k+\beta z)} . C'est un faisceau localement
libre, donc un fibr\'e en droites sur \m{C_2}. On a \
\m{\kl_{k,P,x,\beta\mid C}=\ko_C(-kP)} .
Pour tout entier \m{k\geq 1} on note \m{\ki_{k,P}} le faisceau d'id\'eaux du
sous-sch\'ema \m{kP} de $C$. On a donc \ \m{\ki_{k,P,x}=(x^k,z)} .
On a un isomorphisme
\[(*) \quad\quad\ki_{k,P\mid C} \ \simeq \ \ko_C(-kP)\oplus T_k ,\]
o\`u \m{T_k} est le faisceau de torsion nul en dehors de $P$ et tel que
\m{T_{k,P}=\ko_{C,P}/(x^k)} . L'isomorphisme \m{(*)} au point $P$ est
\xmat{\alpha x^k+\gamma z\fmaps[r] & (\pi(\alpha x^k),p(\gamma)) }
(o\`u \m{\pi:\ko_{C_2,P}\to\ko_{C,P}} et \m{p:\ko_{C_2,P}\to\ko_{C,P}/(x^k)}
sont les projections), donc il d\'epend du choix de $x$.\end{subsub}
\end{sub}

\sepsub

\Ssect{Plongements dans un faisceau quasi localement libre}{C2_x6_0b}

On consid\`ere le faisceau sans torsion $\ke$ de \ref{C2_x3_2}. On conserve les
notations de \ref{ST_C2}.

\sepprop

\begin{subsub}\label{C2_x7}{\bf Proposition : } Soit
\m{\mu\in\Ext^1_{\ko_2}(E\ot L^*,G_\ke)} l'\'el\'ement correspondant \`a la
suite exacte \
\m{0\to G_\ke\to\ke\to E\ot L^*\to 0} \ du diagramme \m{(D)}. Alors

1 - Il existe un \m{\nu\in\Ext^1_{\ko_2}(V_\ke\ot L^*,G_\ke)} tel que $\mu$
soit l'image de $\nu$ par l'application
\[\Ext^1_{\ko_2}(V_\ke\ot L^*,G_\ke)\lra\Ext^1_{\ko_2}(E\ot L^*,G_\ke)\]
induite par l'inclusion \m{E\subset V_\ke}.

2 - Soit \ \m{0\lra G_\ke\lra\ku\lra V_\ke\ot L^*\lra 0} \ la suite exacte
correspondant \`a $\nu$. Alors $\ku$ est quasi localement libre sur \m{C_2}.
\end{subsub}

\begin{proof}
On consid\`ere la suite exacte \ \m{0\lra E\lra V_\ke\lra T\lra 0} . On en
d\'eduit le morphisme
\[\delta : \Ext^1_{\ko_2}(E\ot L^*,G_\ke)\lra\Ext^2_{\ko_2}(T\ot L^*,G_\ke).
\]
Il suffit, pour obtenir 1-, de prouver que \m{\delta(\mu)=0}.

Les d\'efinitions de \m{V_\ke}, \m{G_\ke} sont {\em locales} : on peut
d\'efinir, pour tout \m{P\in C} et tout \m{\ko_{2,P}}-module sans torsion de
type fini $M$, les \m{\ko_{2,P}}-modules \m{V_M}, \m{G_M}, de telle sorte que
\m{V_{\ke P}=V_{\ke_P}}, \m{G_{\ke P}=G_{\ke_P}}. On a \
\m{\Ext^2_{\ko_2}(T\ot L^*,G_\ke)=H^0(\\Ext^2_{\ko_2}(T\ot L^*,G_\ke))} ,
et on a, pour tout \m{P\in C} un diagramme commutatif \'evident
\xmat{
\Ext^1_{\ko_2}(E\ot L^*,G_\ke)\ar[rr]^-{ev}\ar[dd]^-\delta & &
\Ext^1_{\ko_{2,P}}(E_P\ot L^*_P,G_{\ke_P})\ar[dd]^-{ev_P}\\
\\
\Ext^2_{\ko_2}(T\ot L^*,G_\ke)\ar[rr]^-{ev} & &
\Ext^2_{\ko_{2,P}}(T_P\ot L^*_P,G_{\ke_P})
}
Pour montrer que \m{\delta(\mu)=0}, il suffit donc de montrer que pour tout
\m{P\in C} on a \m{\delta_P(ev(\mu))=0}.

Les constructions de \ref{ST_C2} sont {\em additives} : on a
\m{V_{\ke\oplus\ke'}=V_\ke\oplus V_{\ke'}},\Nligne
\m{G_{\ke\oplus\ke'}=G_\ke\oplus G_{\ke'}}, etc... C'est aussi le cas des
constructions
locales. Il suffit donc, d'apr\`es le corollaire \ref{DUAL_2b}, de traiter le
cas o\`u \m{(\ki_{k,P})_P}. On est donc ramen\'e au cas o\`u \m{\ke=\ki_{k,P}}.
L'assertion 1- d\'ecoule donc du diagramme commutatif avec lignes exactes
\xmat{
0\ar[r] & L=G_{\ki_{k,P}}\ar[r]\fleq[d] & \ki_{k,P}\ar[r]\flinc[d] &
\ko_C(-kP)=E\ar[r]\flinc[d] & 0\\
0\ar[r] & L\ar[r] & \ko_2\ar[r] & \ko_C=V_{\ki_{k,P}}\ar[r] & 0
}
et de la proposition 4.3.1 de \cite{dr1}.

L'assertion 2- est aussi locale. On peut donc encore se ramener au cas o\`u
\m{\ke=\ki_{k,P}}. Il faut montrer que pour tout diagramme commutatif avec
lignes exactes de \m{\ko_{2,P}}-modules
\xmat{
0\ar[r] & L_P\ar[r]\fleq[d] & (\ki_{k,P})_P\ar[r]\flinc[d] &
\ko_C(-kP)_P\ar[r]\flinc[d] & 0\\
0\ar[r] & L_P\ar[r] & M\ar[r] & \ko_{C,P}\ar[r] & 0
}
on a \m{M\simeq\ko_{2,P}}. On a \ \m{L_P\simeq\ko_C(-kP)_P\simeq\ko_{C,P}}, et
\m{\Ext^1_{\ko_{2,P}}(\ko_{C,P},\ko_{C,P})\simeq\ko_{C,P}}
et on v\'erifie ais\'ement que l'\'el\'ement de
\m{\Ext^1_{\ko_{2,P}}(\ko_{C,P},\ko_{C,P})} associ\'e \`a la suite exacte du
haut est de la forme \m{\alpha x^k_P} (avec \m{\alpha\in\ko_{C,P}} inversible).
D'autre part le morphisme \ \m{\ko_{C,P}=\ko_C(-kP)_P\hookrightarrow\ko_{C,P}}
\ est la multiplication par \m{x_P^k}. Donc l'\'el\'ement de
\m{\Ext^1_{\ko_{2,P}}(\ko_{C,P},\ko_{C,P})} associ\'e \`a la suite exacte du
bas est $\alpha$, d'o\`u \m{M\simeq\ko_{2,P}}. Ceci  d\'emontre 2- .
\end{proof}

\sepprop

\begin{subsub}\label{C2_x7b}{\bf Corollaire : } Il existe un faisceau quasi
localement libre $\kv$ sur \m{C_2} et un diagramme commutatif avec lignes et
colonnes exactes
\xmat{
 & & 0\ar[d] & 0 \ar[d]\\
0\ar[r] & G_\ke\ar[r]\fleq[d] & \ke\ar[r]\ar[d] & E\ot L^*\ar[r]\ar[d] & 0\\
0\ar[r] & G_\ke\ar[r] & \kv\ar[r]\ar[d] & V_\ke\ot L^*\ar[r]\ar[d] & 0\\
 & & T\ar[d]\fleq[r] & T\ar[d]\\
 & & 0 & 0
}
En particulier, $\ke$ est isomorphe au noyau d'un morphisme surjectif d'un
faisceau quasi localement libre dans $T$.
\end{subsub}
\end{sub}

\sepsub

\Ssect{Dualit\'e}{DUAL_C2}

Soient \m{P\in C} et $M$ un  \m{\ko_{2,P}}-module de type fini. On note
\m{M^\vee} le {\em dual} de $M$ :\Nligne
\m{M^\vee=\Hom(M,\ko_{2,P})}~.
Si $N$ est un \m{\ko_{C,P}}-module, on note \m{N^*} le dual de $N$ :
\Nligne\m{N^*=\Hom(N,\ko_{C,P})}.

Soit $\ke$ un faisceau coh\'erent sur \m{C_n}. On note \m{\ke^\vee} le {\em
dual} de $\ke$ : \
\m{\ke^\vee\simeq\HHom(\ke,\ko_2)} .
Si $E$ est un faisceau coh\'erent sur $C$, on note \m{E^*} le dual de $E$ :
\m{E^*=\Hom(E,\ko_C)}.
Ces notations sont justifi\'ees par le fait que \m{E^\vee\not=E^*} :

\sepprop

\begin{subsub}{\bf Lemme : }\label{DUAL_1}
1 - Soit $N$ un \m{\ko_{C,P}}-module de type fini. Alors on a \ \m{N^\vee\simeq
N^*}. Si $N$ est libre on a \ \m{\Ext^i_{\ko_{2,P}}(N,\ko_{2,P})=\nsp} \
si \m{i\geq 1}.

2 - Soit $E$ un faisceau coh\'erent sur $C$. Alors on a \ \m{E^\vee\simeq
E^*\ot L}. Si $E$ est localement libre on a \ \m{\Ext^i_{\ko_2}(E,\ko_2)=\nsp} \
si \m{i\geq 1}.
\end{subsub}

\begin{proof}
On ne d\'emontrera que la seconde assertion, la premi\`ere \'etant analogue.
Soient $\F$ un fibr\'e vectoriel sur \m{C_2} et \m{F=\F_{\mid C}}. On a donc
une suite exacte canonique \
\m{0\to F\ot L\to\F\to F\to 0} .
En examinant ce qui se passe en chaque point de $C$ on en d\'eduit ais\'ement
la suite exacte duale \
\m{0\to F^\vee\to\F^\vee\to (F\ot L)^\vee\to 0} .
On a un morphisme canonique \ \m{\F^\vee\lra F^*} \ \'evident. En examinant ce
qui se passe en chaque point de $C$ on voit ais\'ement que ce morphisme est
surjectif et que son noyau est exactement \m{F^\vee}. On a donc \
\m{(F\ot L)^\vee=F^*},
d'o\`u \ \m{F^\vee=F^*\ot L}. Pour d\'emontrer 2- on peut se limiter au cas
o\`u $E$ est localement libre. D'apr\`es le th\'eor\`eme \ref{pr5} il existe un
prolongement de $E$ en un fibr\'e vectoriel $\E$ sur \m{C_2}. D'apr\`es ce
qui pr\'ec\`ede on a bien \m{E^\vee=E^*\ot L}.

Soit $\L$ une extension de $L$ \`a \m{C_2}. On a une r\'esolution localement
libre de $E$, en tant que faisceau sur \m{C_2} :
\[\cdots\lra\F\ot\L^2\lra\F\ot\L\lra\F\lra F\lra 0 ,\]
obtenue en juxtaposant les suites exactes \
\m{0\to F\ot L^{i+1}\to\F\ot\L^i\lra F\ot L^i\to 0} .
On en d\'eduit imm\'ediatement que \ \m{\Ext^i_{\ko_2}(E,\ko_2)=0} \ si
\m{i\geq 1}.
\end{proof}

\sepsub

\begin{subsub}{\bf Dualit\'e des faisceaux quasi localement libres - }
\label{DUAL_QLL}\rm Les propri\'et\'es suivantes se d\'emontrent ais\'ement
ce qui se passe en chaque point de $C$. 
Si $\kf$ est un faisceau coh\'erent quasi localement libre, alors on a
\[E_{\kf^\vee} \ \simeq \ E_\kf^*\ot L^2 , \ \ \ \
F_{\kf^\vee} \ \simeq \ G_\kf^*\ot L , \ \ \ \
G_{\kf^\vee} \ \simeq \ F_\kf^*\ot L .\]
Le morphisme \ \m{\Phi_{\kf^\vee}:F_{\kf^\vee}\ot L\to E_{\kf^\vee}} \ est le
morphisme \ \m{G_\kf^*\ot L^2\to E_\kf^*\ot L^2} \ provenant de la suite
exacte  \ \m{0\to E_\kf\to G_\kf\to\Gamma_\kf\to 0} .

Plus g\'en\'eralement, si $\ke$, $\kf$ sont des faisceaux quasi localement
libres sur \m{C_2}, et si\Nligne \m{\kh= \HHom(\ke,\kf)} , alors on a \
\m{E_\kh=\HHom(E_\ke,E_\kf\ot L)} , \m{G_\kh=\HHom(F_\ke,G_\kf)} ,
et des suites exactes
\[0\lra\HHom(\Gamma_\ke,E_\kf)\oplus\HHom(E_\ke,\Gamma_\kf\ot L)\lra F_\kh\lra
\HHom(\Gamma_\ke,\Gamma_\kf)\oplus\HHom(E_\ke,E_\kf)\lra 0 ,\]
\[0\lra\HHom(\Gamma_\ke,E_\kf)\oplus\HHom(E_\ke,\Gamma_\kf\ot L)\lra
\Gamma_\kh\lra\HHom(\Gamma_\ke,\Gamma_\kf)\lra 0 .\]

Il d\'ecoule du lemme \ref{DUAL_1} que si \ \m{0\to\ke'\to\ke\to\ke''\to 0} \
est une suite exacte de faisceaux quasi localement libres sur \m{C_2}, alors
la suite transpos\'ee \ \m{0\to{\ke'}^\vee\to\ke^\vee\to{\ke''}^\vee\to 0} \
est aussi exacte.

Si $\kf$, $\kf'$ sont des faisceaux quasi localement libres sur \m{C_2}, le
morphisme canonique \Nligne
 \m{\kf^\vee\ot\kf'\to\HHom(\kf,\kf')} \ n'est pas en
g\'en\'eral un isomorphisme (par exemple \m{\ko_C^\vee\ot\ko_C=L},
\m{\HHom(\ko_C,\ko_C)=\ko_C}). On a en fait une suite exacte canonique
\[0\lra\Gamma(\kf)^*\ot\Gamma(\kf')\ot L\lra\kf^\vee\ot\kf\lra
\HHom(\kf,\kf')\lra\Gamma(\kf)^*\ot\Gamma(\kf')\lra 0 . \]
\end{subsub}

\sepsubsub

\begin{subsub}{\bf Dualit\'e des faisceaux sans torsion - }\label{DUAL_ST}\rm
Si $M$ est un \m{\ko_{2,P}}-module de type fini et $\ke$ un faisceau coh\'erent
sans torsion dur \m{C_2}, on a des morphismes canoniques \
\m{i_M:M\to M^{\vee\vee}} , \m{i_\ke:\ke\lra\ke^{\vee\vee}} .

On dit que $M$ (resp. $\ke$) est {\em r\'eflexif} si \m{i_M} (resp. \m{i_\ke})
est un isomorphisme.
\end{subsub}

\sepprop

\begin{subsub}{\bf Proposition : }\label{DUAL_2}
1 -  Un \m{\ko_{2,P}}-module de type fini est r\'eflexif si et seulement si il
est sans torsion.

2 - Un faisceau coh\'erent sur \m{C_2} est r\'eflexif si et seulement si il est
sans torsion.
\end{subsub}

\begin{proof}
Il est clair que 2- d\'ecoule de 1-. Soit $M$ un \m{\ko_{2,P}}-module de type
fini. Il est clair que si $M$ est r\'eflexif, il est sans torsion. Pour
d\'emontrer la r\'eciproque on peut utiliser une \'etude de la structure des
\m{\ko_{2,P}}-modules sans torsion analogue \`a \ref{C2_x3_2}. On trouve qu'il
existe un sous-\m{\ko_{2,P}}-module quasi libre $N$ de $M$ tel que \m{T=M/N}
soit de torsion. On en d\'eduit une suite exacte \
\m{0\to M^\vee\to N^\vee\lra\wT\to 0} .
Il est imm\'ediat que \m{M^\vee} est sans torsion. Montrons que \m{M\subset
M^{\vee\vee}} . Il faut prouver que pour tout \m{m\in M} il
existe \m{\phi\in M^\vee} tel que \m{\phi(m)\not=0}. C'est imm\'ediat si $M$
est quasi libre. En g\'en\'eral, il existe un entier \m{k\geq 0} tel que
\m{x^k_Pm\in N}. Puisque $N$ est quasi libre, il existe \m{\psi\in N^\vee} tel
que \m{\psi(x^k_Pm)\not=0}. D'autre part il existe un entier \m{p\geq 0} tel
que \m{x_P^p\psi\in M^\vee}. Il suffit donc de prendre \m{\phi=x_P^{k+p}\psi}.
On a maintenant un diagramme commutatif avec lignes exactes
\xmat{
0\ar[r] & N\fleq[d]\ar[r] & M^{\vee\vee}\ar[r] & T\fleq[d]\\
0\ar[r] & N\ar[r]         & M\flinc[u]\ar[r]   & T\ar[r] & 0
}
On en d\'eduit imm\'ediatement que \m{M^{\vee\vee}=M}.
\end{proof}

\sepprop

Soit \m{m\geq 1} un entier. On note \
\m{I_{m,P}=(x^m,z)} .

\sepprop

\begin{subsub}{\bf Corollaire : }\label{DUAL_2b}
Soit $M$ un \m{\ko_{2,P}}-module sans torsion. Alors il existe des entiers $m$,
$q$ et une suite d'entiers \m{n_1,\ldots,n_p} tels que
\[M \ \simeq \ \biggl(\bigoplus_{i=1}^pI_{n_i,P}\biggr)\oplus m\ko_{2,P} 
\oplus q\ko_{C,P}.\]
\end{subsub}

\begin{proof}
Comme dans la d\'emonstration de la proposition \ref{DUAL_2} on montre qu'il
existe un sous-\m{\ko_{2,P}}-module quasi libre $N$ de \m{M^\vee} tel que
\m{M^\vee/N} soit un \m{\ko_{C,P}}-module de torsion. Il en d\'ecoule que \
\m{M^{\vee\vee}=M\subset N^\vee} \ et que \m{N^\vee/M} est un
\m{\ko_{C,P}}-module de torsion. Le r\'esultat d\'ecoule alors du
lemme suivant.
\end{proof}

\sepprop

\begin{subsub}{\bf lemme : }\label{C2_x3_1}
Soient $M$ un \m{\ko_{2,P}}-module quasi libre, $T$ un \m{\ko_{C,P}}-module de
torsion et \m{\pi:M\to T} un morphisme. Alors il des entiers $p$, $m$, $q$ et
une suite d'entiers\Nligne \m{n_1\geq n_2\geq\cdots\geq n_p\geq 1} tels que
\[\ker(\pi) \ \simeq \ \biggl(\bigoplus_{i=1}^pI_{n_i,P}\biggr)\oplus
m\ko_{2,P}\oplus q\ko_{C,P} .\]
\end{subsub}

\begin{proof}
Soient \m{N=\ker(\pi)}, \m{M_C=M\ot_{\ko_{2,P}}\ko_{C,P}}, et $E$ le noyau de
la restriction \m{M\to M_C}. Le morphisme $\pi$ se factorise par \m{M_C} :
\xmat{M\ar[r] & M_C\ar[r]^-{\ov{\pi}} & T}
Soit \m{F=\ker(\ov{\pi})} . Alors $E$ et $F$ sont des \m{\ko_{C,P}}-modules
libres. On a un diagramme commutatif avec lignes et colonnes exactes :
\xmat{
        & 0\ar[d]         & 0\ar[d] \\
        & E\ar[d]\fleq[r] & E\ar[d] \\
0\ar[r] & N\ar[r]\ar[d]   & M\ar[r]\ar[d]   & T\ar[r]\fleq[d] & 0\\
0\ar[r] & F\ar[r]\ar[d]   & M_C\ar[r]\ar[d] & T\ar[r]         & 0\\
        & 0               & 0
}
Posons \m{E=r\ko_{C,P}}, \m{F=s\ko_{C,P}}. On a donc une suite exacte \
\m{0\to r\ko_{C,P}\to N\to s\ko_{C,P}\to 0} .
On a \ \m{\Ext^1_{\ko_P}(s\ko_{C,P},r\ko_{C,P})=L(\C^s,\C^r)\ot\ko_{C,P}} . La
suite exacte pr\'ec\'edente est donc associ\'ee \`a une matrice \
\m{A=(\phi_{ij})_{1\leq i\leq r,1\leq j\leq s}} \ d\'el\'ements de
\m{\ko_{C,P}}. On peut supposer que tous les \m{\phi_{ij}} ne sont pas nuls. Si
\m{\phi_{ij}\not=0}, on pose \ \m{\phi_{ij}=x^{m_{ij}}\eta_{ij}} , avec
\m{\eta_{ij}} inversible. On peut supposer que \m{\phi_{11}\not=0} et que
\m{m_{11}} est minimal. Dans ce cas en ajoutant des multiples de la premi\`ere
colonne de $A$ aux autres colonnes et en manipulant aussi les lignes, on se
ram\`ene au cas o\`u \m{\phi_{1j}=0} si \m{j\not=1} et \m{\phi_{i1}=0} si
\m{i\not=1}. On fait la m\^eme chose avec la sous-matrice \
\m{A'=(\phi_{ij})_{2\leq i\leq r,2\leq j\leq s}}, et se ram\`ene au cas
suivant : $A$ poss\`ede une sous-matrice carr\'ee $B$ diagonale, et tous ses
autres termes sont nuls. Le r\'esultat d\'ecoule alors imm\'ediatement du fait
que si \m{\beta\in\ko_{C,P}} est de la forme \ \m{\beta=x^m\alpha}, avec
$\alpha$ inversible, et si \
\m{0\to\ko_{C,P}\to E\to\ko_{C,P}\to 0} \
l'extension correspondant \`a $\beta$, alors on a \ \m{E\simeq I_{m,P}} .
\end{proof}

\sepprop

Le r\'esultat suivant d\'ecoule de \ref{C2_x6_0b}. On en donne une autre
d\'emonstration.

\sepprop

\begin{subsub}{\bf Corollaire : }\label{DUAL_2bb}
1 - Soit $M$ un \m{\ko_{2,P}}-module de type fini et sans torsion. Alors il
existe un module quasi libre $N$, un \m{\ko_{C,P}}-module de torsion $T$ et un
morphisme surjectif \m{\phi:N\to T} tels que \m{M\simeq\ker(\phi)}.

2 - Soit $\ke$ un faisceau coh\'erent sans torsion sur \m{C_2}. Alors il existe
un faisceau quasi localement libre $\kg$ sur \m{C_2}, un faisceau de torsion
$T$ sur $C$ et un morphisme surjectif \m{\phi:\kg\to T} tels que
\m{\ke\simeq\ker(\phi)}.
\end{subsub}

\begin{proof}
On ne d\'emontrera que 2-, 1- \'etant analogue. Le faisceau \m{\ke^\vee} est
sans torsion. On a vu dans \ref{ST_C2} qu'il existait un faisceau quasi
localement libre $\kf$ sur \m{C_2}, un faisceau de torsion \m{T_0} sur $C$ et
une suite exacte \
\m{0\to\kf\to\ke^\vee\to T_0\to 0} . 
En dualisant cette suite exacte on obtient la suivante
\xmat{0\ar[r] & \ke^{\vee\vee}\ar[r] & \kf^\vee\ar[r] &
\EExt^1(T_0,\ko_2)\ar[r]^-\phi & \EExt^1(\ke^\vee,\ko_2)\ar[r] & 0\ .}
On a \m{\ke^{\vee\vee}=\ke} d'apr\`es la proposition \ref{DUAL_2}, et
\m{\EExt^1(T_0,\ko_2)=L\ot\widetilde{T_0}} . Il suffit donc de prendre
\m{\kg=\kf^\vee} et \m{T=\ker(\phi)}.
On a en fait \ \m{\EExt^1(\ke^\vee,\ko_2)=0} : cela se d\'emontre en utilisant
le corollaire \ref{DUAL_2b} et des r\'esolutions libres des modules \m{I_{n,P}}
(cf. la d\'emonstration de la proposition \ref{C2_x6_1}).
\end{proof}
\end{sub}

\sepsec

\section{D\'eformations des faisceaux coh\'erents sur les courbes multiples}
\label{def_f_db}

\Ssect{Sch\'emas de Hilbert ponctuels sur \m{C_n}}{TORS}

Soit \m{d\geq 1} un entier. Soit \m{U\subset C_n} un ouvert affine tel qu'il
existe un plongement de $U$ dans une surface affine lisse $S$. On note
\m{\Hilb^d(C_n)_{\mid U}} la sous-vari\'et\'e de \m{\Hilb^d(S)}
constitu\'ee des sous-sch\'emas finis de longueur $d$ contenus dans $U$,
et $\kt$ la restriction \`a \ \m{\Hilb^d(C_n)_{\mid U}\times U} \ du sch\'ema
universel \m{\kt_0} sur \ \m{\Hilb^d(S)\times S}.
La structure de vari\'et\'e de \m{\Hilb^d(C_n)_{\mid U}} est d\'efinie de la
fa\c con suivante : soit \m{p_0:\Hilb^d(S)\times S\to\Hilb^d(S)} la
projection et \m{V=p_{0*}(\kt_0)}, qui est un fibr\'e vectoriel de rang $d$.
Tout \'el\'ement de l'id\'eal de $U$ dans $S$ induit une section de $V$, et
\m{\Hilb^d(C_n)_{\mid U}} est pr\'ecis\'ement le lieu des z\'eros de ces
sections. Le sch\'ema
\m{\Hilb^q(C_n)_{\mid U}} muni de \m{\kt} est le {\em sch\'ema de Hilbert} des
sous-sch\'emas finis de longueur $q$ de $U$.

\sepprop

\begin{subsub}{\bf Proposition : }\label{C2_x6_1}
Soit \m{Z\subset U} un sous-sch\'ema fini de longueur $d>0$ de support $P$,
\m{T=\ko_Z}. Alors si $Z$ est contenu dans $C$ le morphisme canonique \
\m{\Ext^1_{\ko_n}(T,T)\to\Ext^1_{\ko_S}(T,T)} \
est un isomorphisme. Autrement dit l'inclusion \
\m{\Hilb^d(C_n)_{\mid U}\subset\Hilb^d(S)} \
induit un isomorphisme entre les espaces tangents en $Z$.
\end{subsub}

\begin{proof} 
Cela d\'ecoule du fait que l'\'equation de $U$ dans $S$ est une puissance 
$n$-i\`eme de celle de $C$.
\end{proof}

\sepprop

Le r\'esultat pr\'ec\'edent implique que toute d\'eformation infinit\'esimale
(c'est-\`a-dire param\'etr\'ee par \m{\spec(\C[t]/(t^2))}) de $T$ en tant que
faisceau sur $S$ en est une de $T$ en tant que faisceau sur \m{C_n},
c'est-\`a-dire que c'est une famille de faisceaux dont le support est contenu
dans \m{C_n}.

Les d\'eformations de $Z$ en tant que sous-sch\'ema de $C$ induisent des
d\'eformations de $T$ : elles correspondent au sous-espace \
\m{\Ext^1_{\ko_C}(T,T)\subset\Ext^1_{\ko_n}(T,T)} \
de codimension $d$.

Quand $n=2$ on peut d\'ecrire enti\`erement les sous-sch\'emas $T$ de $C_2$ de
support $P$ :
Soit \m{Z} un sous-sch\'ema fini de longueur $d>0$ de support $P$,
\m{T=\ko_Z}. Alors le sch\'ema $Z$ est une intersection compl\`ete dans $U$.
Deux cas peuvent se produire :
\begin{itemize}
\item[-] $I=(x^{p+q},zx^q)$, avec \m{p,q\geq 0}, \m{p+q>0}, \m{d=p+2q} .
\item[-] $I=(x^{q+m+p}+zx^q\alpha,zx^{q+m})$, avec \m{m>0}, \m{p,q\geq 0},
\m{p+2q+2m=d}, \m{\alpha\in\ko_{C,P}} inversible.
\end{itemize}

\end{sub}

\sepsub

\Ssect{D\'eformations des faisceaux localement libres sur $C$}{def_f_db1}

\begin{subsub}{\bf Proposition : }\label{DUAL_0}
On suppose que $C_2$ est plong\'ee dans une surface lisse $S$. Soit $E$ un
fibr\'e vectoriel sur $C$. Alors le morphisme canonique \
\m{\Ext^1_{\ko_2}(E,E)\to\Ext^1_{\ko_S}(E,E)} \
est un isomorphisme.
\end{subsub}

\begin{proof}
D'apr\`es le th\'eor\`eme \ref{pr5} il existe un fibr\'e vectoriel $\E$ (resp.
un fibr\'e en droites $\L$) sur \m{C_2} dont la restriction \`a $C$ est $E$
(resp. $L$). On a alors une r\'esolution canonique de $E$ sur \m{C_2}
\[\cdots\lra\E\ot\L^2\lra\E\ot\L\lra\E\lra E \ .\]
On en d\'eduit un isomorphisme \ \m{\EExt^1_{\ko_2}(E,E)\simeq\HHom(E\ot L,E)}.
D'autre part, en utilisant la r\'esolution localement libre de \m{\ko_C} sur
$S$ :
\m{0\to\ko_S(-C)\to\ko_S\to\ko_C\to 0} \
on obtient \ \m{\Tor^1_{\ko_S}(E,\ko_C)\simeq E\ot L}. On en d\'eduit, en
utilisant la suite spectrale des Ext et la proposition \ref{pr3} le diagramme
commutatif avec lignes et colonnes exactes
\xmat{
& 0\ar[d]\\
0\ar[r] & \Ext^1_{\ko_C}(E,E)\ar[r]\ar[d] & \Ext^1_{\ko_S}(E,E)\ar[r]\fleq[d] &
\Hom(E\ot L,E)\ar[d]^\alpha\\
0\ar[r] & \Ext^1_{\ko_{C_2}}(E,E)\ar[r]\ar[d] & \Ext^1_{\ko_S}(E,E)\ar[r]^-\phi
& \Hom(\Tor^1_{\ko_S}(E,\ko_{C_2}),E)\ar[d]\\
& \Hom(E\ot L,E)\ar[d] & & 0\\ & 0\\
}
On en d\'eduit que \ \m{\ker(\alpha)\simeq\Hom(E\ot L,E)}, d'o\`u \m{\alpha=0}
et \m{\phi=0}.
\end{proof}

\sepprop

\begin{subsub}{\bf Remarque : }\label{DUAL_0b}\rm
Ce r\'esultat se g\'en\'eralise ais\'ement de la fa\c con suivante : soient
$n$, $k$ des entiers avec \m{n\geq 2}, \m{1\leq k\leq n}, $Y$
une courbe multiple plong\'ee dans une surface lisse $S$,\Nligne
\m{C_1=C\subset C_2\subset\cdots\subset C_n=Y} \
sa filtration canonique, et $E$ un fibr\'e vectoriel sur \m{C_k}. Alors, si
\m{2k\leq n}, le morphisme canonique \
\m{\Ext^1_{\ko_{C_n}}(E,E)\to\Ext^1_{\ko_S}(E,E)} \
est un isomorphisme.
\end{subsub}

\sepprop

\begin{subsub}{\bf Th\'eor\`eme : }\label{DUAL_0c}
Soient \m{p\geq 0}, $d$ des entiers et $E$ un fibr\'e vectoriel de rang 2 et de
degr\'e $d$ sur $C$. Alors

1 - Si $E$ se d\'eforme en faisceaux sans torsion sur
\m{C_2}, non concentr\'es sur $C$ et d'index $p$, alors il existe un fibr\'e en
droites $V$ sur $C$ de degr\'e \m{\frac{1}{2}(d+\deg(L)+p)} et un morphisme non
nul \m{\alpha:V\to E} tels que \ \m{\Hom((E/\imm(\alpha))\ot L,V)\not=\nsp} .

2 - Si $E$ poss\`ede un sous-fibr\'e en droites $V$ de degr\'e
\m{\frac{1}{2}(d+\deg(L)+p)} tel que\Nligne
\m{\Hom((E/V)\ot L,V)\not=\nsp}, alors $E$ se
d\'eforme en faisceaux sans torsion sur \m{C_2}, non concentr\'es sur $C$ et
d'index $p$.
\end{subsub}

\begin{proof}
D\'emontrons d'abord 2- .
 Soit
\m{D=E/V}. On peut supposer que \m{D\ot L\subset V}, et il existe un
unique sous-sch\'ema \m{Z\subset C} de dimension 0 tel que \m{V/(D\ot
L)\simeq\ko_Z}. Alors il existe des faisceaux $\ke$ sans torsion sur \m{C_2},
non concentr\'es sur $C$, d'index $p$, de premi\`ere filtration canonique
\m{0\subset D\ot L\subset\ke}, tels que \m{V_\ke=V} et que le morphisme
canonique \m{D\ot L\to V_\ke} soit l'inclusion \m{D\ot L\subset V}. On a donc
une suite exacte \
\m{0\to V\to\ke\to D\to 0} .
On a aussi une suite exacte \ \m{0\lra V\lra E\lra D\lra 0}. La famille des
extensions de $D$ par $V$ (param\'etr\'ee par \m{\Ext^1_{\ko_2}(D,V)}) est donc
une d\'eformation de $E$ dont le faisceau g\'en\'erique est sans torsion, non
concentr\'e sur $C$ et d'index $p$.

On d\'emontre maintenant 1-.
Supposons que $E$ soit une d\'eformation de faisceaux sans
torsion non concentr\'es sur $C$ et d'index $p$. Il existe donc un germe de
courbe lisse \m{(Y,0)}, une famille plate $\kf$ de faisceaux coh\'erents sur
\m{C_2} param\'etr\'ee par $Y$ telle que si \m{y\not=0}, \m{\kf_y} soit de rang
g\'en\'eralis\'e 2, sans torsion, non concentr\'e sur $C$, d'index $p$, et que
\m{\kf_0\simeq E}. La famille de fibr\'es en droites sur $C$
\m{(V_{\kf_y})_{y\in Y\backslash\lbrace 0\rbrace}} se prolonge en $0$ et la
restriction \m{\Pic(C_2)\to\Pic(C)} admet des sections locales. Il existe donc
un fibr\'e en droites $\kv$ sur \m{C_2\times Y} tel que pour tout \m{y\in
Y\backslash\lbrace 0\rbrace} on ait \ \m{\kv_{y\mid C}\simeq V_{\kf_y}}. En
rempla\c cant $\kf$ par \m{\kf\ot\kv^{-1}} on se ram\`ene au cas o\`u pour tout
\m{y\in Y\backslash\lbrace 0\rbrace} on a \ \m{V_{\kf_y}=\ko_C}.

On va d\'emontrer 1-, avec \m{V=\ko_C}. Soient $t$ un
g\'en\'erateur de l'id\'eal maximal de \m{\ko_{Y,0}} et \m{p_2:C_2\times Y\to
C_2} \ la projection. Il existe un morphisme
\[\pi : p_2^*(\ko_C)_{\mid Y\backslash\lbrace 0\rbrace}\lra
\kf_{\mid Y\backslash\lbrace 0\rbrace}\]
tel que pour tout \m{y\in Y\backslash\lbrace 0\rbrace} on ait \m{\imm(\pi_y)=
V_{\kf_y}}. On va montrer qu'il existe un entier $k$ tel que \m{t^k\pi} se
prolonge en un morphisme \m{\gamma:p_2^*(\ko_C)\to\kf} non nul en 0. Il suffit
de montrer que c'est vrai au voisinage de tout point de $C$. Soient \m{P\in C},
et \
\m{\sigma:\ko_{C_2\times(Y\backslash\lbrace 0\rbrace)}\lra
\kf_{\mid Y\backslash\lbrace 0\rbrace}} \
le compos\'e de $\pi$ et de la projection
\m{\ko_{C_2\times(Y\backslash\lbrace 0\rbrace)}\to
p_2^*(\ko_C)_{\mid Y\backslash\lbrace 0\rbrace}},
qui est donc une section de \m{\kf_{\mid Y\backslash\lbrace 0\rbrace}}.
Alors \m{\kf_{(P,0)}} est un \m{(\ko_{2,P}\ot\ko_{Y,0})}-module, et on d\'eduit
de $\sigma$ un \'el\'ement \m{\sigma_P} du localis\'e \m{\kf_{(P,0),(t)}}.
Il existe un entier $p$ tel que \m{t^p\sigma_P} se prolonge en un \'el\'ement
\m{\gamma_P} de \m{\kf_{(P,0)}}. On peut supposer que $p$ est minimal. Il reste
\`a montrer que \m{(\gamma_P)_0}, l'image de \m{\gamma_P} dans le
\m{\ko_{2,P}}-module \m{(\kf_0)_P}, est non nul. Si \m{(\gamma_P)_0=0}, on a
\m{\gamma_P\in t\kf_{(P,0)}}. On a donc \m{\gamma_P=t\nu}, avec \m{\nu\in
\kf_{(P,0)}}, d'o\`u \m{t^p\sigma_P=t\nu}. Mais, puisque $\kf$ est plat sur
$Y$, la multiplication par $t$ est injective. Donc \m{t^{p-1}\sigma_P=\nu}, ce
qui contredit la minimalit\'e de $p$. On a donc \m{(\gamma_P)_0\not=0}.

Il en d\'ecoule que \m{\ko_C\subset E} et
\m{\kd=\coker(\gamma)} est une famille plate de faisceaux
de rang 1 sur $C$. Pour tout \m{y\in Y\backslash\lbrace 0\rbrace}, on a \
\m{\Hom(\kd_y\ot L,\ko_C)\not=\nsp} \ (car \m{\kd_y=E_{\kf_y}} et
\m{\ko_C=V_{\kf_y}}). Donc par semi-continuit\'e on a \ \m{\Hom((E/\ko_C)\ot
L,\ko_C)\not=\nsp}.
\end{proof}
\end{sub}

\sepsub

\Ssect{D\'eformations des faisceaux quasi localement libres sur les courbes
doubles}{def_f_db2}

Soient $\kf$ un faisceau quasi localement libre sur \m{C_2} et \m{0\subset
E\subset\kf} sa premi\`ere filtration canonique. On pose \
\m{r_0(\kf)=rg(E)} .
On a toujours \
\m{R(\kf)\geq 2r_0(\kf)} ,
avec \'egalit\'e si et seulement si $\kf$ est localement libre.

Si \m{P\in C}, on d\'efini de m\^eme l'entier \m{r_0(M)} pour tout
\m{\ko_{2,P}}-module quasi libre $M$.

\sepprop

\begin{subsub}{\bf Proposition : }\label{Def_1}
Soient \m{P\in C}, $M$ un \m{\ko_{2,P}}-module quasi libre, et \m{r_0} un
entier tel que \ \m{0<2r_0\leq R(M)}. Alors $M$ se d\'eforme en
\m{\ko_{2,P}}-modules quasi libres $N$ tels que \ \m{r_0(N)=r_0} \ si et
seulement si on a \ \m{r_0\geq r_0(M)}.
\end{subsub}

\begin{proof}
On a \ \m{rg(M)=R(M)-r_0(M)}. Les modules $N$ d\'eformations de $M$ sont de
rang g\'en\'eralis\'e \m{R(M)} et de rang
au plus \m{rg(M)}, donc \ \m{r_0(N)\geq r_0(M)}. R\'eciproquement, supposons
que \ \m{0<2r_0\leq R(M)} \ et \ \m{r_0>r_0(M)}. Soit \m{s=R(M)-2r_0(M)}, donc
\Nligne \m{M\simeq r_0(M)\ko_{2,P}\oplus s\ko_{C,P}} .
On doit montrer que $M$ se d\'eforme en modules isomorphes \`a \
\m{r_0\ko_{2,P}\oplus (s-2(r_0-r_0(M)))\ko_{C,P}} .
Il suffit de montrer que \m{2(r_0-r0(M))\ko_{C,P}} se d\'eforme en modules
isomorphes \`a
\m{(r_0-r_0(M))\ko_{2,P}}. Cela revient \`a montrer que
\m{2\ko_{C,P}} se d\'eforme en modules isomorphes \`a \m{\ko_{2,P}}. Cela se
voit ais\'ement en consid\'erant des extensions de \m{\ko_{2,P}}-modules
\ \m{0\to\ko_{C,P}\to A\to\ko_{C,P}\to 0} .
\end{proof}

\sepprop

Il en d\'ecoule que si un faisceau quasi localement libre $\ke$ sur \m{C_2} se
d\'eforme en faisceaux quasi localement libres de rang \m{r_0}, alors on a \
\m{r_0\geq r_0(\ke)}. La r\'eciproque est fausse : si par exemple \m{R(\ke)}
est pair, et si $\ke$ se d\'eforme en faisceaux localement libres, on doit
avoir \ \m{\Deg(\ke)\equiv\frac{R(\ke)}{2}\deg(L)\ ({\rm mod}\ 2)}. J'ignore si
cette condition est suffisante pour que $\ke$ se d\'eforme en faisceaux
localement libres.

\sepsubsub

\begin{subsub}\label{def_f_db2_0}{\bf D\'eformations de modules - }\rm
Soit \m{z\in\ko_{2P}} un g\'en\'erateur de l'id\'eal de $C$.
On a \ \m{\Ext_{\ko_{2P}}^i(\ko_{CP},\ko_{CP})\simeq\ko_{CP}} \
si \m{i\geq 1}. Cela se d\'emontre en utilisant la r\'esolution libre de
\m{\ko_{CP}} :
\xmat{(K^\bullet)\quad\quad\cdots\ko_{2P}\ar[r]^-{\times z} &
\ko_{2P}\ar[r]^-{\times z} & \ko_{2P}\ar[r]^-{\times z} & \ko_{2P}}
Si \m{i,j\leq 1}, le produit \ \m{\Ext_{\ko_{2P}}^i(\ko_{CP},\ko_{CP})\times
\Ext_{\ko_{2P}}^j(\ko_{CP},\ko_{CP})\to\Ext_{\ko_{2P}}^{i+j}(\ko_{CP},\ko_{CP})}
\ est la multiplication \ \m{\ko_{CP}\times\ko_{CP}\to\ko_{CP}} . Cela se voit
en interpr\'etant les \'el\'ements de \m{\Ext_{\ko_{2P}}^k(\ko_{CP},\ko_{CP})}
comme des morphismes de degr\'e \m{-k} de \m{(K^\bullet)} dans lui-m\^eme.

Soient $k$ un entier positif et \m{A_k=\C\lbrack t \rbrack/(t^k)} .
On appelle {\em d\'eformation d'ordre $k$} d'un \m{\ko_{2,P}}-module $M$ un
\m{(\ko_{2P}\ot_\C A_k)}-module $\bf M$ plat sur \m{A_k} tel que
\m{{\bf M}/t{\bf M}=M} . On appelle {\em d\'eformation d'ordre infini} de $M$ un
\m{(\ko_{2P}\ot_\C\C\lbrack t\rbrack)}-module $\bf M$ plat sur
\m{\C\lbrack t\rbrack} tel que \m{{\bf M}/t{\bf M}=M} . Les d\'eformations
d'ordre 2 de $M$ sont param\'etr\'ees naturellement par
\m{\Ext^1_{\ko_{2P}}(M,M)} . Si \m{\sigma\in\Ext^1_{\ko_{2P}}(M,M)} et si
\m{M_\sigma} est le \m{(\ko_{2P}\ot_\C A_2)}-module correspondant, alors
\m{M_\sigma} s'\'etend en une d\'eformation d'ordre 3 si et seulement si
\m{\sigma^2=0} dans \m{\Ext^1_{\ko_{2P}}(M,M)} (cf. \cite{dr1}, 3) . Soient $r$
un entier positif, \m{\sigma\in\Ext^1_{\ko_{2P}}(r\ko_{CP},r\ko_{CP})\simeq
\End(r\ko_{CP})}, et \m{M_\sigma} le \m{(\ko_{2P}\ot_\C A_2)}-module
correspondant. Alors \m{M_\sigma} s'\'etend en une d\'eformation d'ordre 3 si et
seulement si \m{\sigma\circ\sigma\in\End(r\ko_{CP})} est nul.

Soit \m{\phi\in\End(r\ko_{CP})} tel que \m{\phi^2=0} . On montre ais\'ement
qu'il existe une base de \m{r\ko_{CP}} dans laquelle la matrice de $\phi$ est
une matrice diagonale de matrices du type 0 ou \m{\begin{pmatrix}0 & 0\\ 0 &
\alpha\end{pmatrix}}. Il en d\'ecoule d'abord qu'une d\'eformation d'ordre 2 de
\m{r\ko_{CP}} pouvant s'\'etendre en une d\'eformation d'ordre 3 se d\'ecompose
en une somme directe de d\'eformations de \m{2\ko_{CP}} et de d\'eformations
triviales de \m{\ko_{CP}}.

Soit \m{x\in\ko_{CP}} un g\'en\'erateur de l'id\'eal maximal. Soit \m{\phi\in
\End(2\ko_{CP})} tel que \m{\phi^2=0}. On peut supposer que la matrice de $\phi$
est \m{\begin{pmatrix}0 & 0\\ 0 & x^k\end{pmatrix}} , \m{k\geq 0}. Une telle
d\'eformation s'obtient en consid\'erant le faisceau d'id\'eaux \m{\ki_{kP}} de
$kP$. On a une suite exacte canonique
\[0\lra L\lra\ki_{kP}\lra\ko_C(kP)\lra 0 .\]
Soit \m{\sigma\in\Ext^1_{\ko_2}(\ko_C(kP),L)} correspondant \`a l'extension
pr\'ec\'edente. On consid\`ere la famille d'extensions param\'etr\'ee par $\C$
\[0\lra p_{\ko_2}^*(L)\lra\ke\lra p_{\ko_2}^*(\ko_C(kP))\lra 0\]
(o\`u \m{p_{\ko_2}} est la projection \m{C_2\times\C\to C_2}) telle qu'en tout
\m{t\in\C} l'extension \Nligne
\m{0\to L\to\ke_t\to\ko_C(kP)\to 0} \ correspond \`a
\m{t\sigma}. Alors \m{(\ke_{tP})_{t\in\C}} est une d\'eformation d'ordre infini
de \m{2\ko_{CP}} qui est une extension de l'extension d'ordre 2 d\'efinie par
$\phi$. Les d\'eformations en modules libres sont celles pour lesquelles  
\m{k=0} . 

Il d\'ecoule de ce qui pr\'ec\`ede qu'une extension d'ordre 2 de \m{r\ko_{CP}}
s'\'etend en une d\'eformation d'ordre 3 si et seulement si elle s'\'etend en
une d\'eformation d'ordre infini.

Soit \ \m{M=a\ko_{2P}\oplus b\ko_{CP}} \ un \m{\ko_{2P}}-module quasi libre.
Alors on a \Nligne
 \m{\Ext^1_{\ko_{2P}}(M,M)=\Ext^1_{\ko_{2P}}(b\ko_{CP},b\ko_{CP})} .
Les d\'eformations d'ordre 2 de $M$ se r\'eduisent donc \`a celles de
\m{b\ko_{CP}} .
\end{subsub}

\sepsubsub

\begin{subsub}\label{def_f_db2_1}{\bf D\'eformations de faisceaux - }\rm
Soit $\ke$ un faisceau coh\'erent quasi localement libre sur \m{C_2}. Alors on a
d'apr\`es la suite spectrale des Ext et \ref{QLL_C2_4_1} une suite exacte
\xmat{0\ar[r] & H^1(\EEnd(\ke))\ar[r] & \Ext^1_{\ko_2}(\ke,\ke)\ar[r]^-\pi &
H^0(\EExt^1_{\ko_2}(\ke,\ke))=\HHom(\Gamma_\ke\ot L,\Gamma_\ke)\ar[r] & 0
\quad .}
Soit $\kf$ une d\'eformation d'ordre 2 de $\ke$ correspondant \`a \m{\sigma\in
\Ext^1_{\ko_2}(\ke,\ke)} . Il d\'ecoule de \ref{def_f_db2_0} que $\kf$ s'\'etend
en une d\'eformation d'ordre 3 si et seulement si le morphisme \
\m{\pi(\sigma)^2:\Gamma_\ke\ot L^2\to\Gamma_\ke} \ est nul.

{\bf Conjecture : } $\kf$ s'\'etend en une d\'eformation param\'etr\'ee par
une courbe lisse si et seulement si \m{\pi(\sigma)^2=0} .
\end{subsub}

\end{sub}

\sepsec

\section{Faisceaux d'id\'eaux de points}
\label{GEN_FIB}

Soient $C_n$ une courbe projective multiple primitive de multiplicit\'e $n>1$,
plong\'ee dans une surface projective lisse $S$, \
\m{C=C_1\subset C_2\subset\cdots\subset C_n} \
la filtration canonique et  $L$ le fibr\'e en droites sur $C$ associ\'e
(cf. \ref{cour_mul}).
On pose, pour \m{1\leq i\leq n}, \m{\ko_i=\ko_{C_i}} .  Soient $\L$ le faisceau
d'id\'eaux de $C$ dans \m{C_n} (qui est un fibr\'e en droites sur \m{C_{n-1}} et
\m{L=\L_{\mid C}} .

Si \m{P\in C}, on note \m{x_P} un \'el\'ement de \m{\ko_{n,P}} au dessus d'un
g\'en\'erateur de l'id\'eal maximal de \m{\ko_{C,P}}, et \m{z_P} un
g\'en\'erateur de l'id\'eal de $C$ dans \m{\ko_{n,P}}. L'id\'eal maximal de
\m{\ko_{n,P}} est donc \m{(x_P,z_P)}. On note \m{\kl_P} le faisceau d'id\'eaux
sur \m{C_n} \'egal \`a \m{(x_P)} en $P$ et \`a \m{\ko_{n,P'}} en \m{P'\not=P} . 
C'est un faisceau inversible. Notons que contrairement \`a ce qu'indique la
notation, \m{\kl_P} d\'epend du choix de \m{x_P}. Soit \m{\kj_{n,P}} le faisceau
d'id\'eaux sur \m{C_n} \'egal \`a \m{(x_P,z_�^{n-1})} en $P$ et \`a
\m{\ko_{n,P'}} en \m{P'\not=P} . Il d\'epend de \m{x_P} mais pas de \m{z_P} . 

\sepsub

\Ssect{D\'eformations des faisceaux d'id\'eaux de points}{def_id}

Soit $Z$ un ensemble fini non vide de $p$ points de $C$. On note \m{\ki_{n,Z}}
le faisceau d'id\'eaux de $Z$ sur \m{C_n}. On note \ \m{\kl_Z=\bigotimes_{P\in
Z}\kl_P} , \m{\kj_{n,Z}=\bigotimes_{P\in Z}\kj_{n,P}} .

\sepprop

\begin{subsub}\label{propX1}
{\bf Proposition : } Soient $\kd$ un fibr\'e en droites sur \m{C_n},
\m{\kd_{n-1}=\kd_{\mid C_{n-1}}}, \m{D=\kd_{\mid C}}. On suppose que la
restriction \ \m{H^0(\kd)\to H^0(\kd_{n-1})} \ est surjective.
Alors on a une suite exacte canonique \
\m{0\to H^0(\ko_C(Z)\ot L^{n-1}\ot D)\to\Hom(\ki_{n,Z},\ki_{n,Z}\ot\kd)\to
H^0(\kd_{n-1})\to 0 } .
\end{subsub}

\begin{proof}
On a \m{\L\subset\ki_{n,Z}} car \m{Z\subset C}, et \ \m{\ki_{n,Z}/\L\simeq
\ko_C(-Z)} . Soit \Nligne
\m{\phi\in\Hom(\ki_{n,Z},\ki_{n,Z}\ot\kd)} . Alors on a
\ \m{\phi(\L)\subset\L\ot\kd} : pour le voir on utilise la caract\'erisation
suivante de $\L$ : en $P\in C$, si \m{z\in\ko_{nP}} engendre l'id\'eal de $C$,
\m{\L_P\subset\ki_{n,Z,P}} est l'annulateur de \m{z^{n-1}}. Le morphisme
\ \m{\phi:\L\to\L\ot\kd} \ \'equivaut \`a une section $s$ de \m{\kd_{n-1}} .
Soit $\sigma$ le morphisme \ \m{\ki_{n,Z}\to\ki_{n,Z}\ot\kd} \ induit par
une section de $\kd$ au dessus de $s$. Alors on a \ \m{(\phi-s)(\L)=0} , donc
\m{\phi-\sigma} est induit par un morphisme \ \m{\ko_C(-Z)\to\ki_{n,Z}\ot\kd} .
L'image de ce morphisme est contenue dans le sous-faisceau constitu\'e en $P$
de l'annulateur de $z$, ce sous-faisceau est isomorphe \`a \m{L^{n-1}\ot D}.
La proposition \ref{propX1} en d\'ecoule imm\'ediatement.
\end{proof}

\sepprop

En particulier, si \m{\kd=\ko_n}, on obtient \ 
\m{\End(\ki_{n,Z})\simeq\C\oplus H^0(\ko_C(Z)\ot L^{n-1})} , et si\Nligne
\m{H^0(D)=H^0(L)=\nsp}, \m{\Hom(\ki_{n,Z},\ki_{n,Z}\ot\kd)\simeq
H^0(\ko_C(Z)\ot L^{n-1}\ot D)} .

\sepprop

\begin{subsub}\label{propX2}
{\bf Proposition : } On a \ \m{\EEnd(\ki_{n,Z})\simeq
\kj_{n,Z}\ot\kl_Z^{-1}} .
\end{subsub}

\begin{proof}
Cela d\'ecoule ais\'ement du fait que pour tout \m{P\in Z}, les endomorphismes
de \m{\ki_{n,P}} sont du type
\xmat{ax_P+bz_P\fmaps[r] & \alpha(ax_P+bz_P)+\beta az^{n-1} ,}
avec \m{\alpha,\beta\in\ko_{n,P}} .
\end{proof}

\sepprop

\begin{subsub}\label{propX3}
{\bf Proposition : } On suppose que \ \m{h^0(L)=0} \ et que \ 
\m{h^0(K_{S\mid C})=0} ou \ \m{K_{S\mid C_n}\simeq\ko_n} . Alors on
a
\[\dim(\Ext^1_{\ko_{C_n}}(\ki_{n,Z},\ki_{n,Z})) \ = \ 1 + \frac{n^2}{2}C^2 +
\frac{n}{2}K_SC + p + h^0(\ko_C(Z)\ot L^{n-1}) ,\]
\[\dim(\Ext^1_{\ko_{S}}(\ki_{n,Z},\ki_{n,Z})) \ = \ 1 + n^2C^2 +
h^0(\ko_C(Z)\ot L^{n-1}) + h^0(\ko_C(Z)\ot L^{n-1}\ot K_S) +
h^0(K_{S\mid C_n}) .\]
\end{subsub}

\begin{proof}
D\'emontrons d'abord la premi\`ere \'egalit\'e. Soit \ \m{T=\ko_n/\kj_{n,Z}} .
C'est un faisceau de support $Z$ tel que \ \m{h^0(T)=(n-1)p} . On a une suite
exacte, d'apr\`es la proposition \ref{propX2}
\[0\lra\EEnd(\ki_{n,Z})\lra\kl_Z^{-1}\lra T\ot\kl_Z^{-1}\lra 0 ,\]
d'o\`u on d\'eduit la suite exacte
\[0\lra \End(\ki_{n,Z})\lra H^0(\kl_Z^{-1})\lra H^0(T)\lra
H^1(\EEnd(\ki_{n,Z}))\lra H^1(\kl_Z^{-1})\lra 0 .\]
D'apr\`es la proposition \ref{propX1} on a \ \m{\End(\ki_{n,Z})\simeq
\C\oplus H^0(\ko_C(Z)\ot L^{n-1})} , donc
\[h^1(\EEnd(\ki_{n,Z})) \ = \ -\chi(\kl_Z^{-1})+(n-1)p+h^0(\ko_C(Z)\ot L^{n-1})
+ 1 .\]
On a \ \m{h^0(\kl_Z^{-1}/\ko_n)=np} , donc \ \m{\chi(\kl_Z^{-1})=np-
\frac{n^2}{2}C^2-\frac{n}{2}K_SC} . D'apr\`es la suite spectrale des Ext on a \
\m{\dim(\Ext^1_{\ko_{C_n}}(\ki_{n,Z},\ki_{n,Z}))=h^1(\EEnd(\ki_{n,Z}))+
h^0(\EExt^1_{\ko_{C_n}}(\ki_{n,Z},\ki_{n,Z}))} . On a \
\m{\EExt^1_{\ko_{C_n}}(\ki_{n,Z},\ki_{n,Z}))\simeq 2\ko_Z} . Cela se d\'emontre
ais\'ement en utilisant en chaque point $P$ de $Z$ la r\'esolution suivante de
\m{(\ki_{n,Z})_P} :
\xmat{\ko_{n,P}^2\ar[rrr]^{\begin{pmatrix}0 & z_P^{n-1}\\ z_P &
-x_P\end{pmatrix}}
& & & {\ko_{n,P}^2}\ar[rrr]^{\begin{pmatrix}x_P & z_P^{n-1}\\
z_P & 0\end{pmatrix}} & & &
{\ko_{n,P}^2}\ar[rr]^{(z_P,x_P)} & & \ko_{n,P}}
On en d\'eduit imm\'ediatement la premi\`ere \'egalit\'e.

D\'emontrons maintenant la seconde. On a, sur $S$,
\m{\chi(\ki_{n,Z},\ki_{n,Z})=-n^2C^2} . D'autre part, d'apr\`es la proposition
\ref{propX1} et la dualit\'e de Serre
\[\dim(\Ext^2_{\ko_{S}}(\ki_{n,Z},\ki_{n,Z})) \ = \
\dim(\Hom(\ki_{n,Z},\ki_{n,Z}\ot K_S)) \ = \ h^0(\ko_C(Z)\ot L^{n-1}\ot K_S)
+ h^0(K_{S\mid C_n}) ,\]
d'o\`u le r\'esultat.
\end{proof}

\sepprop

\begin{subsub}\label{XXX1}\rm
Sous les hypoth\`eses de la proposition \ref{propX3} on a donc
\begin{eqnarray*}
\dim(\Ext^1_{\ko_{S}}(\ki_{n,Z},\ki_{n,Z})/
\Ext^1_{\ko_{C_n}}(\ki_{n,Z},\ki_{n,Z})) & = & \frac{(n-1)^2}{2}C^2-
\frac{n-1}{2}K_SC\\ & & + \ h^0(\ko_C(-Z)\ot L^{-n}) + h^0(K_{S\mid C_n}) .
\end{eqnarray*}
\end{subsub}

\sepsubsub

\begin{subsub}\label{def_id2}{\bf Les cas du plan projectif et des surfaces 
K3 - }\rm Si $S$ est \m{\P_2} ou une surface K3, on a \ \m{H^0(\kn_{C_n})=
H^0(\ko_S(nC))/\langle\sigma^n\rangle} , o\`u \m{\kn_C} d\'esigne le fibr\'e
normal de $C$ dans $S$ et \m{\sigma\in H^0(\ko_S(C))} une \'equation de $C$.

Soit $\kf$ une d\'eformation plate de \m{\ki_{n,Z}} (en tant que faisceau sur
$S$) param\'etr\'ee par une vari\'et\'e alg\'ebrique (ou analytique)
connexe $U$, \m{x\in U} un point ferm\'e tel que \m{\kf_x\simeq\ki_{n,Z}} .
On en d\'eduit une d\'eformation de \m{C_n}, d'o\`u l'application canonique
\[\Theta^\kf_x : T_xU\lra H^0(\ko_S(nC))/\langle\sigma^n\rangle .\]
Soit \m{\ke_0} une d\'eformation semi-universelle de \m{\ki_{n,Z}} (en tant que
faisceau sur $S$) param\'etr\'ee par un germe de vari\'et\'e analytique
\m{(X_0,x_0)}. Alors \m{\Theta^\kf_x} se factorise \`a travers
\[\Theta^{\ke_0}_{x_0} : T_{x_0}X_0=\Ext^1_{\ko_S}(\ki_{n,Z},\ki_{n,Z})
\lra H^0(\ko_S(nC))/\langle\sigma^n\rangle .\]
En particulier, l'image de \m{\Theta^\kf_x} est contenue dans celle de
\m{\Theta^{\ke_0}_{x_0}} .
\end{subsub}

\sepprop

\begin{subsub}\label{theoX1}
{\bf Th\'eor\`eme : } L'application \m{\Theta^{\ke_0}_{x_0}} induit
un isomorphisme
\[\Ext^1_{\ko_S}(\ki_{n,Z},\ki_{n,Z})/
\Ext^1_{\ko_{C_n}}(\ki_{n,Z},\ki_{n,Z})\ \simeq \ V/\langle\sigma^n\rangle ,\]
o\`u \m{V\subset H^0(\ko_S(nC))} est l'espace des courbes passant par tous
les points de $Z$.
\end{subsub}

\begin{proof}
Le noyau de \m{\Theta^{\ke_0}_{x_0}} contient \m{\Ext^1_{\ko_{C_n}}(\ki_{n,Z},
\ki_{n,Z})}, qui correspond aux faisceaux de support \m{C_n}. On a d'apr\`es
\ref{XXX1}
\[\dim(V/\langle\sigma^n\rangle) \ = \
\dim(\Ext^1_{\ko_S}(\ki_{n,Z},\ki_{n,Z})/
\Ext^1_{\ko_{C_n}}(\ki_{n,Z},\ki_{n,Z})) ,\]
donc il suffit de prouver que \ \m{V/\langle\sigma^n\rangle\subset
\imm(\Theta^{\ke_0}_{x_0})} . Pour cela il suffit de construire une
d\'eformation $\kf$ de \m{\ki_{n,Z}} comme pr\'ec\'edemment telle que
\m{V/\langle\sigma^n\rangle\subset\imm(\Theta^\kf_x)} . On prend \m{U=V} et
$\kf$ est la famille des faisceaux d'id\'eaux de $Z$.
\end{proof}

\sepprop

\begin{subsub}\label{exX1}\rm
{\bf Exemple : } Soient $x$, $y$, $z$ des coordonn\'ees de \m{\P_2}. Supposons
que $C$ soit la droite d'\'equation \m{z=0}, et que \ \m{Z=\lbrace (x_i,y_i,0);
1\leq i\leq p\rbrace} . Pour \m{1\leq i\leq p}, soit \m{z_i} un nombre complexe
tel que \ \m{z_i^n=x_i^n+y_i^n}. Pour \m{\alpha\in\C}, soit \m{\kc_\alpha} la
courbe d'\'equation \ \m{z^n=\alpha^nx^n+\alpha^ny^n} , qui est lisse si
\m{\alpha\not=0}. Notons que \m{C_n=\kc_0}. On a \
\m{Z_\alpha=\lbrace(x_i,y_i,\alpha z_i);1\leq i\leq p\rbrace\subset\kc_\alpha} .
Soit \m{\ki_\alpha} le faisceau d'id\'eaux de \m{Z_\alpha} dans \m{\kc_\alpha},
vu comme faisceau sur \m{\P_2}. Alors \m{\ka=(\ki_\alpha)_{\alpha\in\C}} est une
d\'eformation de \m{\ki_{n,Z}} . L'ensemble $\Gamma$ des courbes \m{\kc_\alpha}
est une courbe de \m{\P(H^0(\ko_{\P_2}(n)))} passant par \m{\langle z^n\rangle}.
La tangente \`a cette courbe en \m{\langle z^n\rangle} n'est pas contenue dans
\m{V/\langle z^n\rangle} . Ceci n'est pas contradictoire avec le th\'eor\`eme
\ref{theoX1} car on voit ais\'ement que \
\m{d\Theta^\ke_{\langle z^n\rangle}=0} .
\end{subsub}

\end{sub}

\sepsub

\Ssect{Faisceaux de rang g\'en\'eralis\'e 2 sur les courbes doubles}{XXX2}

On s'int\'eresse ici aux faisceaux sans torsion de rang g\'en\'eralis\'e 2 sur
\m{C_2}. Les plus simples sont les fibr\'es en droites sur \m{C_2} et les
fibr\'es vectoriels de rang 2 sur $C$. Il reste \`a \'etudier les faisceaux
sans torsion de rang g\'en\'eralis\'e 2, non concentr\'es sur $C$ et d'index
positif.

Soient $\kd$ un fibr\'e en droites sur \m{C_2} de degr\'e $d$, et
\m{D=\kd_{\mid C}}. On a donc une suite exacte \
\m{0\to D\ot L\to\kd\to D\to 0} , et \
\m{\deg(D) \ = \ (d-\deg(L))/2} .
Si maintenant $\ke$ est un faisceau coh\'erent sans torsion de rang
g\'en\'eralis\'e 2 et de degr\'e $d$ sur \m{C_2}, non concentr\'e sur $C$, et
si \ \m{0\subset E\ot L\subset\ke} \ est sa filtration canonique, alors $E$ est
un fibr\'e en droites sur $C$. On a donc une suite exacte
\[0\lra E\ot L\lra\ke\lra E\oplus T\lra 0 ,\]
$T$ \'etant un faisceau de torsion sur $C$. On a donc
\[\deg(E) \ = \ \frac{d-\deg(L)-i(\ke)}{2} .\]
Les seuls autres faisceaux de rang g\'en\'eralis\'e 2, de degr\'e $d$ et sans
torsion sur \m{C_2} sont les fibr\'es vectoriels de rang $2$ et de degr\'e $d$
sur $C$.

\sepprop

Soit $\ke$ un faisceau coh\'erent sans torsion de rang g\'en\'eralis\'e 2 sur
\m{C_2}, non concentr\'e sur $C$. Alors il existe d'apr\`es le corollaire
\ref{DUAL_2bb} un fibr\'e en droites $\kf$ sur \m{C_2}, un faisceau de torsion
$T$ sur $C$ et une suite exacte \
\m{0\to\ke\to\kf\to T\to 0} . On en d\'eduit imm\'ediatement la

\sepprop

\begin{subsub}{\bf Proposition : }\label{C2_x7c}
Soit $\ke$ un faisceau coh\'erent sans torsion de rang g\'en\'eralis\'e 2 sur
\m{C_2}, non concentr\'e sur $C$. Alors il existe un fibr\'e en droites $\kf$
sur \m{C_2} et un unique sous-sch\'ema fini \m{Z\subset C} tels que \
\m{\ke\simeq\ki_Z\ot\kf} ,
\m{\ki_Z} d\'esignant le faisceau d'id\'eaux de $Z$ dans \m{C_2}.
\end{subsub}

\sepprop

Le fibr\'e en droites $\kf$ n'est pas en g\'en\'eral unique.

\sepprop

\begin{subsub}Dualit\'e - \rm
On utilise les notations de \ref{faisc_doub}. On utilise les notations de
\ref{faisc_doub}. Si \m{n\geq 1} est un entier, on note \m{\kl(n,P,x_P)} le
faisceau d'id\'eaux du sous-sch\'ema de \m{C_2} de support $P$ engendr\'e par
\m{x_P^n} en $P$. Comme la notation le sugg\`ere il d\'epend du choix de
\m{x_P}. C'est un fibr\'e en droites sur \m{C_2}. Soit $Z$ un un sous-sch\'ema
fini \m{Z\subset C}. On va d\'ecrire le faisceau dual de \m{\ki_Z}. Posons
\[Z \ = \ \sigg_i n_iP_i ,\]
avec \m{n_i\geq 1}, les \m{P_i} \'etant des points distincts de \m{C}. Alors on
a
\[\ki_Z^\vee \ \simeq\ki_Z\ot\bigg(\psigg_i \kl(n_i,P_i,x_{P_i})\bigg)^{-1}
\ .\]
\end{subsub}

\sepsubsub

\begin{subsub}{D\'eformations infinit\'esimales des faisceaux de rang
g\'en\'eralis\'e 2 - }\label{def_r2}\rm
On rappelle maintenant et on pr\'ecise les constructions de \ref{ST_C2}.
Soient $\ke$ un faisceau coh\'erent sans torsion
de rang g\'en\'eralis\'e 2 sur \m{C_2}, non concentr\'e sur $C$, et \m{D\ot
L\subset\ke} sa premi\`ere filtration canonique, $D$ \'etant un fibr\'e en
droites sur $C$. Alors on a \ \m{\ke_{\mid C}\simeq D\oplus T}, $T$ \'etant un
faisceau de torsion sur $C$. Fixons un isomorphisme
\[\eta:\ke_{\mid C} \ \simeq \ D\oplus T .\]
On en d\'eduit la suite exacte \
\m{0\to D\ot L\to\ke\to D\oplus T\to 0} .
Soit \Nligne
\m{(\nu_\eta,\sigma_\eta)\in\Ext_{\ko_2}^1(D,D\ot L)\oplus
\Ext^1_{\ko_2}(T,D\ot L)} \
l'\'el\'ement associ\'e. Le noyau \m{\kf_\eta} du morphisme compos\'e \
\m{\ke\to D\oplus T\to T} \ est un fibr\'e en droites sur \m{C_2}, et on a une
suite exacte
\[0\lra D\ot L\lra\kf_\eta\lra D\lra 0\]
qui est associ\'ee \`a \m{\nu_\eta}.

Les automorphismes de \m{D\oplus T} sont repr\'esent\'es par des matrices
\m{\tau=\begin{pmatrix}\alpha & 0\\\phi & \beta\end{pmatrix}}, avec
\m{\alpha\in\C^*}, \m{\beta\in\Aut(T)} et \m{\phi:D\to T}. On a
\m{v_{\tau\circ\eta} = \alpha v_\eta+\phi\sigma_\eta} , 
\m{\sigma_{\tau\circ\eta} = \beta\sigma_\eta} .\Nligne
Puisque \ \m{\Ext^1_{\ko_2}(T,D\ot L)=\Ext^1_{\ko_C}(T,D\ot L)}, la
multiplication par \Nligne
\m{\sigma_\eta : \Hom(D,T)\to\Ext^1_{\ko_2}(D,D\ot L)} \ est
\`a valeurs dans \m{\Ext^1_{\ko_C}(D,D\ot L)=H^1(L)}. On a \ \m{H^1(L)\subset
H^1(\ko_2)}. On peut donc voir la multiplication par \m{\sigma_\eta} comme une
application
\[\beta_\eta:\Hom(D,T)\lra H^1(\ko_2) .\]
Le rang de cette application ne d\'epend pas du choix de $\eta$.
\end{subsub}

\sepprop

Un calcul \'el\'ementaire mais fastidieux permet de d\'emontrer la

\sepprop

\begin{subsub}{\bf Proposition : }\label{C2_x7d}
On a des suites exactes canoniques
\[0\lra H^0(\ko_2)\lra\End(\ke)\lra\ker(\beta_\eta)\lra 0 ,\]
\[0\lra\Ext^1_{\ko_2}(T,T)\lra\Ext^1_{\ko_2}(\ke,\ke)\lra\coker(\beta_\eta)\lra
0.\]
La seconde suite exacte est scind\'ee.
\end{subsub}

\sepprop

{\bf Cons\'equences : } On peut montrer en utilisant la proposition
\ref{C2_x7d} que les seules d\'eformations de $\ke$ proviennent de
d\'eformations de $T$ et de \m{\kf_\eta}. Les d\'eformations de $T$ sont
d\'ecrites en \ref{TORS}.
\end{sub}

\sepsec

\section{Faisceaux de rang g\'en\'eralis\'e 3 sur les courbes doubles}
\label{GEN_FIB2}

Soient $S$ une surface projective lisse irr\'eductible et \m{C\subset S} une
courbe projective lisse irr\'eductible. Soient \m{C_2\subset S} la courbe
double associ\'ee, \m{L=\ko_C(-C)} et \m{l=-\deg(L)}. On suppose que
\m{l=C^2\geq 1}. Le genre de $C$ est \ \m{g=\frac{1}{2}(C^2+K_SC)+1} .

Soit $\ke$ un faisceau quasi localement libre de rang g\'en\'eralis\'e 3 sur
\m{C_2}, non concentr\'e sur $C$. Il est donc localement du type
\m{\ko_2\oplus\ko_C}. Pour fixer les id\'ees, rappelons qu'on a un diagramme
commutatif avec lignes et colonnes exactes
\xmat{& & & 0\ar[d]\\
& 0\ar[d] & & \Gamma_\ke\ar[d]\\
0\ar[r] & E_\ke\ar[r]\ar[d] & \ke\fleq[d]\ar[r] & F_\ke\ar[r]\ar[d] & 0\\
0\ar[r] & G_\ke\ar[r]\ar[d] & \ke\ar[r] & E_\ke\ot L^*\ar[r]\ar[d] & 0\\
& \Gamma_\ke\ar[d] & & 0\\
& 0}
la premi\`ere ligne exacte provenant de la premi\`ere filtration canonique de
$\ke$ et la seconde ligne exacte de la seconde filtration canonique de $\ke$.
On a \ \m{rg(E_\ke)=rg(\Gamma_\ke)=1}, \m{rg(F_\ke)=rg(G_\ke)=2}.

Le rang et le degr\'e des fibr\'es vectoriels \m{E_\ke}, \m{F_\ke}, \m{G_\ke},
\m{\Gamma_\ke} sont invariants par d\'eformation. Cette propri\'et\'e se 
g\'en\'eralise aux faisceaux quasi localement libres de type \m{(m-1,m)} ou 
\m{(m,m)}.

\sepsub

\Ssect{Faisceaux (semi-)stables}{RGEN3_2}

\begin{subsub}\label{stab_def}{\bf D\'efinition : } Soit $\ke$ un faisceau
coh\'erent sur \m{C_2} pur de dimension 1. On dit que $\ke$ est 
{\em semi-stable} (resp. stable)
s'il est sans torsion et si pour tous sous-faisceau propre \m{\kf\subset\ke} on
a \ \m{\mu(\kf)\leq\mu(\ke)} (resp. $<$).
\end{subsub}

\sepprop

Il d\'ecoule du th\'eor\`eme de Riemann-Roch \ref{RR0} et du calcul de
polyn\^omes de Hilbert effectu\'e dans \ref{hilb_pol} que cette notion de
(semi-)stabilit\'e est \'equivalente \`a celle de Simpson (cf. \cite{si}). Il
existe donc des vari\'et\'es de modules de faisceaux semi-stables de polyn\^ome
de Hilbert donn\'e.

Soit $\ke$ un faisceau quasi localement libre de rang g\'en\'eralis\'e 3 sur
\m{C_2}, non concentr\'e sur $C$. On pose \ \m{\epsilon=\deg(E_\ke)}, \
\m{\gamma=\deg(\Gamma_\ke)} . Si $\ke$ est semi-stable (resp. stable) on
obtient, en consid\'erant les sous-faisceaux \m{E_\ke}, \m{G_\ke}, que \
\m{\gamma-2l\leq\epsilon\leq l+\gamma} (resp $<$). Le r\'esultat suivant montre
que la (semi-)stabilit\'e de $\ke$ d\'epend uniquement des propri\'et\'es des
fibr\'es de rang 2 \m{F_\ke} et \m{G_\ke} :

\sepprop

\begin{subsub}\label{stab_lem}{\bf Lemme : } Soit $\ke$ un faisceau
quasi localement libre de rang g\'en\'eralis\'e 3 sur \m{C_2}, non concentr\'e
sur $C$. Alors $\ke$ est semi-stable (resp. stable) si et seulement si les deux
propri\'et\'es suivantes sont v\'erifi\'ees :
\begin{enumerate}
\item[(i)] Pour tout sous-fibr\'e en droites $D'$ de $G_\ke$ on a \
$\deg(D')\leq\mu(\ke)$ (resp. $<$).
\item[(ii)] Pour tout fibr\'e en droites quotient $D''$ de $F_\ke$ on a \
$\mu(\ke)\leq\deg(D'')$ (resp. $<$).
\end{enumerate}
\end{subsub}

\begin{proof}
On ne traitera que le cas de la semi-stabilit\'e, l'autre cas \'etant analogue.
Les conditions sont \'evidemment n\'ecessaires. Supposons les v\'erifi\'ees.
Soit \m{\kf\subset\ke} un sous-faisceau de rang 1 ou 2. Il faut montrer que
\m{\mu(\kf)\leq\mu(\ke)} . Si $\kf$ est de rang 1,
\m{E_\kf\subset E_\ke} est sans torsion, donc $\kf$ est un faisceau sur $C$,
c'est donc un sous-faisceau de \m{G_\ke}, et il est localement libre de rang 1.
Soit \m{F\subset G_\ke} l'image inverse du sous-faisceau de torsion de
\m{G_\ke/\kf}, c'est un sous-fibr\'e de \m{G_\ke}. On a \m{\mu(F)\leq\mu(\ke)}
d'apr\`es (i), et \m{\mu(\kf)\leq\mu(F)}, donc \m{\mu(\kf)\leq\mu(\ke)}.
Supposons $\kf$ de rang 2. Alors \m{\ke/\kf} est de rang 1, donc \m{E_{\ke/\kf}}
est de torsion. Soit $\kg$ son image inverse dans $\ke$. Alors on a
\m{\mu(\kf)\leq\mu(\kg)} et il suffit de prouver que \m{\mu(\kg)\leq\mu(\ke)} .
Cela d\'ecoule du fait que \m{\ke/\kg} est un quotient de \m{F_\ke}, car
\m{\ke/\kg} est un faisceau sur $C$, et de (ii).
\end{proof}

On en d\'eduit que si \m{G_\ke} et \m{F_\ke} sont semi-stables (resp. stables)
il en est de m\^eme de $\ke$.

\end{sub}

\Ssect{Vari\'et\'es de modules}{RGEN3_3}

On suppose que \m{\gamma-l<\epsilon<\gamma}, ce qui \'equivaut \`a
\m{\mu(E_\ke)<\mu(G_\ke)} et \m{\mu(\Gamma_\ke)<\mu(F_\ke)} . 

On note \m{\km_s(\epsilon,\gamma)} la vari\'et\'e de modules des faisceaux quasi
localement $\ke$ libres stables de rang g\'en\'eralis\'e 3 tels que
\m{\deg(E_\ke)=\epsilon}, \m{\deg(\Gamma_\ke)=\gamma} et que \m{F_\ke},
\m{G_\ke} soient stables. C'est un ouvert de la vari\'et\'e de modules
\m{M(3,2\epsilon+\gamma+l)} des faisceaux semi-stables de rang g\'en\'eralis\'e
3 et de degr\'e g\'en\'eralis\'e \m{2\epsilon+\gamma+l}.

\sepprop

\begin{subsub}\label{stab_prop}{\bf Proposition : } La vari\'et\'e
\m{\km_s(\epsilon,\gamma)} est irr\'eductible de dimension \m{5g+2l-4} . La
sous-vari\'et\'e r\'eduite associ\'ee est lisse.
\end{subsub}

\begin{proof}
Soient \m{\iota\in W^{2,\epsilon-\gamma+l}_0} (cf. \ref{BN}), correspondant \`a
un morphisme \m{\ko\to F_0}. Soient \m{\Gamma\in J^\gamma}, \m{F=F_0\ot\Gamma},
\m{E=(F/\Gamma)\ot L\in J^\epsilon}. On a donc une suite exacte
\xmat{0\ar[r] & \gamma\ar[r]^i & F\ar[r]^-p & E\ot L^*\ar[r] & 0 \ .}
Rappelons qu'on aune suite exacte
\xmat{0\ar[r] & \Ext^1_{\ko_C}(F,E)\ar[r] & \Ext^1_{\ko_2}(F,E)\ar[r]^-\pi &
\Hom(F\ot L,E)\ar[r] & 0 \ .}
Soit \m{\sigma\in\pi^{-1}(p)} . Soit \ \m{0\to E\to\ke_\sigma\to F\to 0} \
l'extension associ\'ee, o\`u \m{\ke_\sigma} est quasi localement libre de rang
g\'en\'eralis\'e 3, \m{E_{\ke_\sigma}=E}, \m{\Gamma_{\ke_\sigma}=\Gamma},
\m{F_{\ke_\sigma}=F} .

\begin{subsub}\label{stab_prop2}{\bf Lemme : } Soit \m{\eta\in
\Ext^1_{\ko_C}(\Gamma,E)\subset\Ext^1_{\ko_2}(\Gamma,E)} correspondant \`a la
suite exac-\break te \ \m{0\lra E\lra G_{\ke_\sigma}\lra\Gamma\lra 0} . Alors
$\eta$ est l'image de $\sigma$ par l'application \Nligne
\m{\theta:\Ext^1_{\ko_2}(F,E)\to
\Ext^1_{\ko_2}(\Gamma,E)} \ d\'eduite de l'inclusion \m{\Gamma\subset F}.
\end{subsub}

\begin{proof}
Cela se d\'emontre en utilisant des r\'esolutions libres ad\'equates de
$\Gamma$, $E$. On en d\'eduit des r\'esolutions libres de $F$, \m{\ke_\sigma},
\m{G_{\ke_\sigma}} permettant de repr\'esenter les \m{\Ext^1} et de prouver le
lemme.
\end{proof}

\begin{subsub}\label{stab_prop3}{\bf Lemme : } Pour tout \m{\eta_0\in
\Ext^1_{\ko_C}(\Gamma,E)} il existe \m{\sigma_0\in\pi^{-1}(p)} tel que
\m{\theta(\sigma_0)=\eta_0} .
\end{subsub}

\begin{proof}
Soit \m{\theta_C:\Ext^1_{\ko_C}(F,E)\to\Ext^1_{\ko_C}(\Gamma,E)} \ l'application
d\'eduite de l'inclusion \m{\Gamma\subset F}. Alors \m{\theta_C} est surjective.
Soit \m{\alpha\in\Ext^1_{\ko_C}(F,E)} tel que \m{\theta_C(\alpha)=\eta_0-\eta} .
Il suffit de prendre \m{\sigma_0=\sigma+\alpha} .
\end{proof}

{\em Fin de la d\'emonstration de la proposition \ref{stab_prop} - } D'apr\`es
le lemme \ref{stab_prop3} on peut choisir $\sigma$ de telle sorte que
\m{G_{\ke_\sigma}} soit stable. Dans ce cas
\m{\ke_\sigma\in\km_s(\epsilon,\gamma)} . Le reste de la proposition
\ref{stab_prop} se d\'emontre ais\'ement.
\end{proof}

\end{sub}


\vskip 2cm

\begin{thebibliography}{99}
\bibitem{ar} Artin, M. {\em On deformations of singularities.} Tata Institute
of Fundamental Research Lect. Notes. 54 (1976).
\bibitem{ba_fo} B\u anic\u a, C., Forster, O. {\em Multiple structures on plane
curves.} In Contemporary Mathematics 58, Proc. of Lefschetz Centennial Conf.
(1986), AMS, 47-64.
\bibitem{be_fr} Beorchia, V., Franco, D. {\em On the Moduli Space of 't Hooft
Bundles.} Ann. Univ. Ferrara, sez. VII, Sc. Mat. Vol. XLVII (2001), 253-268.
\bibitem{bho} Bhosle Usha N. {\em Generalized parabolic bundles and applications
to torsion free sheaves on nodal curves.} Arkiv for Matematik 30 (1992),
187-215.
\bibitem{bho2} Bhosle Usha N. {\em Picard groups of the moduli spaces of vector
bundles.} Math. Ann. 314 (1999) 245-263.
\bibitem{bo} Boratynski, M., {\em Locally complete intersection multiple
structures on smooth algebraic curves.} Proc. of the Amer. Math. Soc. 115
(1992), 877-879.
\bibitem{BMNO} Brambila-Paz, L., Mercat, V., Newstead, P.E., Ongay, F. {\em
Nonemptiness of Brill-Noether loci.} Intern. J. Math. 11 (2000), 737-760.
\bibitem{dr} Dr\'ezet, J.-M. {\em Vari\'et\'es de modules alternatives}. Ann.
de l'Inst. Fourier 49 (1999), 57-139.
\bibitem{dr1}Dr\'ezet, J.-M. {\em D\'eformations des extensions larges de
faisceaux} . Preprint (2002), math.AG/0201125.
\bibitem{go} Godement, R. {\em Th\'eorie des faisceaux.} Actualit\'es
scientifiques et industrielles 1252, Hermann, Paris (1964).
\bibitem{ha} Hartshorne, R. {\em Algebraic Geometry.} GTM 52, Springer-Verlag
(1977).
\bibitem{hu_le} Huybrechts, D., Lehn, M. {\em The Geometry of Moduli Spaces of
Sheaves.} Aspect of Math. E31, Vieweg (1997).
\bibitem{in} Inaba, M.-A. {\em On the moduli of stable sheaves on some
nonreduced projective schemes.} Journ. of Alg. Geom. 13 (2004), 1-27.
\bibitem{in2} Inaba, M.-A. {\em On the moduli of stable sheaves on a reducible
projective scheme and examples on a reducible quadric surface.} Nagoya Math.
J. (2002), 135-181.
\bibitem{lp}Le Potier, J. {\em Syst\`emes coh\'erents et structures de niveau.}
Ast\'erisque 214 (1993).
\bibitem{lp4} Le Potier, J. {\em Faisceaux semi-stables de dimension 1 sur le
plan projectif.} Revue roumaine de math. pures et appliqu\'ees 38 (1993),
635-678.
\bibitem{man0} Manolache, N. {\em Cohen-Macaulay Nilpotent Structures.} Revue
Roumaine Math. pures et appl. 31 (1986), 563-575.
\bibitem{man} Manolache, N. {\em Multiple Structures on Smooth Support.} Math.
Nachr. 167 (1994), 157-202.
\bibitem{man0b} Manolache, N. {\em Double rational normal curves with linear
syzygies.} Manusc. Math. 104 (2001), 503-517.
\bibitem{man1} Manolache, N. {\em Cohen-Macaulay Nilpotent Schemes.} Preprint
(2003), math.AG/0312514.
\bibitem{ma1} Maruyama, M. {\em Moduli of stable sheaves I.} J. Math. Kyoto
Univ. 17 (1977), 91-126.
\bibitem{ma2} Maruyama, M. {\em Moduli of stable sheaves II.} J. Math. Kyoto
Univ. 18 (1978), 577-614.
\bibitem{ma_tr}Maruyama, M., Trautmann, G. {\em Limits of instantons.}
Intern. Journ. of Math. 3 (1992), 213-276.
\bibitem{na_ra} Narasimhan, M.S., Ramanan, S. {\em Moduli of vector bundles on
a compact Riemann surface.} Ann. of Math. 89 (1969), 14-51.
\bibitem{nu_tr}N\"ussler, T., Trautmann, G. {\em Multiple Koszul structures on
lines and instanton bundles.} Intern. Journ. of Math. 5 (1994), 373-388.
\bibitem{ra-vi} Raghavendra, N., Vishwanath, P.A. {\em Moduli of pairs and
generalized theta divisors.} Tohoku Math. J. 46 (1994), 321-340.
\bibitem{rego} Rego, C.J. {\em Deformation of modules on curves and surfaces.}
Singularities, Representations of Algebras and Vector Bundles, Proceedings
Lambrecht 1985, Lect. Notes in Math. 1273, Springer-Verlag (1987), 157-167.
\bibitem{ro} Rosenlicht, M. {\em Generalized Jacobian varieties.} Ann. of Math.
59 (1954), 505-530.
\bibitem{se} Serre, J.-P. {\em Groupes alg\'ebriques et corps de classes.}
Hermann, Paris (1959).
\bibitem{ses} Seshadri, C.S. {\em Fibr\'es vectoriels sur les courbes
alg\'ebriques.} Ast\'erisque 96 (1982).
\bibitem{si} Simpson, C.T. {\em Moduli of representations of the fundamental
group of a smooth projective variety I.} Publ. Math. IHES 79 (1994), 47-129.
\bibitem{si_tr} Siu Y., Trautmann, G. {\em Deformations of coherent
analytic sheaves with compact supports.} Memoirs of the Amer.
Math. Soc., Vol. 29, N. 238 (1981).
\bibitem{sun0} Sun, X. {\em Degeneration of moduli spaces and generalized
theta functions.} Journ. of Alg. Geom. 9 (2000), 459-527.
\bibitem{sun1} Sun X., {\em Degeneration of SL(n)-bundles on a reducible curve.}
Proc. of Algebraic Geometry in East Asia (Japan, 2001) and math.AG/0112072.
\bibitem{sun} Sundaram, N. {\em Special divisors and vector bundles.} Tohoku
Math. J. 39 (1987), 175-213.
\bibitem{va} Vatne, J.E. {\em Multiple structures.} Thesis. Preprint (2002),
math.AG/0210042.

\end{thebibliography}
\end{document}